%% file: Frittelli_et_al.tex
\documentclass[letterpaper, 11pt]{article}
\usepackage{amsmath}
\usepackage{amssymb}
\usepackage{amsthm}
\usepackage[english]{babel}
\usepackage{caption}
\usepackage{color}
\usepackage{enumitem}
\usepackage[T1]{fontenc}
\usepackage[left=2.5cm, right=2.5cm, top=2.5cm, bottom = 2.5cm]{geometry}
\usepackage{graphicx}
\usepackage[hidelinks]{hyperref}
\usepackage[utf8]{inputenc}
\usepackage{makecell}
\usepackage{mathtools}
\usepackage[numbers,sort]{natbib}
\usepackage{placeins}
\usepackage[center]{subfigure}
\usepackage{tikz}
\usepackage[normalem]{ulem} 

\usetikzlibrary{calc}

\definecolor{lightGreen}{rgb}{0.6, 1, 0.6}
\definecolor{darkGreen}{rgb}{0.0, 0.5, 0.0}
\definecolor{lightBlue}{rgb}{0.8, 0.8, 1}
\definecolor{lightRed}{rgb}{1, 0.8, 0.8}
\definecolor{purple}{rgb}{0.55,0.4,1}

\newcommand{\red}[1]{\textcolor{red}{#1}}
\newcommand{\blue}[1]{\textcolor{blue}{#1}}
\newcommand{\green}[1]{\textcolor{darkGreen}{#1}}

\newcommand{\bolda}{\boldsymbol{a}}
\newcommand{\boldf}{\boldsymbol{f}}
\newcommand{\boldg}{\boldsymbol{g}}
\newcommand{\boldn}{\boldsymbol{n}}
\newcommand{\boldx}{\boldsymbol{x}}
\newcommand{\boldy}{\boldsymbol{y}}

\newcommand{\boldv}{\boldsymbol{v}}
\newcommand{\boldw}{\boldsymbol{w}}
\newcommand{\boldalpha}{\boldsymbol{\alpha}}
\newcommand{\boldnu}{\boldsymbol{\nu}}

\newcommand{\boldxi}{\boldsymbol{\xi}}
\newcommand{\boldeta}{\boldsymbol{\eta}}

\newcommand{\intomega}{\int_{\Omega}}
\newcommand{\intomegah}{\int_{\Omega_h}}
\newcommand{\intgamma}{\int_{\Gamma}}
\newcommand{\intgammah}{\int_{\Gamma_h}}

\newcommand{\ut}{\dfrac{\partial u}{\partial t}}
\newcommand{\vt}{\dfrac{\partial v}{\partial t}}
\newcommand{\Ut}{\dfrac{\partial U}{\partial t}}
\newcommand{\Vt}{\dfrac{\partial V}{\partial t}}
\newcommand{\nablagamma}{\nabla_{\Gamma}}
\newcommand{\nablagammah}{\nabla_{\Gamma_h}}

\DeclareMathOperator{\area}{area}
\DeclareMathOperator{\ess}{ess}
\DeclareMathOperator{\Tr}{Tr}

\newtheorem{assumption}{Assumption}
\newtheorem{corollary}{Corollary}
\newtheorem{definition}{Definition}
\newtheorem{lemma}{Lemma}
\newtheorem{proposition}{Proposition}
\newtheorem{remark}{Remark}
\newtheorem{theorem}{Theorem}

\title{Bulk-surface virtual element method for systems of PDEs in two-space dimensions}

\author{Massimo Frittelli\footnote{University of Salento, Department of Innovation Engineering, Via per Monteroni, 73100 Lecce, Italy. Email: \texttt{massimo.frittelli@unisalento.it}}, Anotida Madzvamuse\footnote{University of Sussex, School of Mathematical and Physical Sciences, Department of Mathematics, BN1 9QH, Brighton, UK, \& University of Johannesburg, Department of Mathematics, South Africa and Universit\'{a} degli Studi di Bari Aldo Moro, Bari, Italy. Email: \texttt{A.Madzvamuse@sussex.ac.uk}}, Ivonne Sgura\footnote{University of Salento, Department of Mathematics and Physics ``E. De Giorgi'', Via per Arnesano, 73100 Lecce, Italy. Email: \texttt{ivonne.sgura@unisalento.it}}}

\date{}

\begin{document}

\maketitle

\begin{abstract}
In this paper we consider a coupled bulk-surface PDE in two space dimensions. The model consists of a PDE in the bulk that is coupled to another PDE on the surface through general nonlinear boundary conditions. For such a system we propose a novel method, based on coupling a virtual element method \cite{beirao2013basic} in the bulk domain to a surface finite element method \cite{dziuk2013finite} on the surface. {The proposed method, which we coin the {\it Bulk-Surface Virtual Element Method (BSVEM)} includes, as a special case, the bulk-surface finite element method (BSFEM) on triangular meshes \cite{Madzvamuse_2016}. The method exhibits} second-order convergence in space, provided the exact solution is $H^{2+1/4}$ in the bulk and $H^2$ on the surface, where the additional $\frac{1}{4}$ is required only {in the simultaneous presence} of surface curvature and non-triangular elements. Two novel techniques introduced in our analysis are (i) an $L^2$\emph{-preserving inverse trace operator} for the analysis of boundary conditions and (ii) the \emph{Sobolev extension} as a replacement of the lifting operator \cite{elliottranner2013finite} for sufficiently smooth exact solutions. The generality of the polygonal mesh can be exploited to optimize the computational time of matrix assembly. The method takes an optimised matrix-vector form that also simplifies the known special case of BSFEM on triangular meshes \cite{Madzvamuse_2016}. Three numerical examples illustrate our findings.
\end{abstract}

\section*{Keywords}
Bulk-surface PDEs, Bulk-surface finite elements, Bulk-surface virtual elements, Bulk-surface reaction-diffusion systems, Virtual elements

\section*{Mathematics Subject Classification}
65M12, 65M15, 65M20, 65M60, 65N30

\section{Introduction}
\label{sec:introduction}

{An interesting class of PDE problems that is recently drawing attention in the literature is that of \emph{coupled bulk-surface partial differential equations} (BSPDEs). Given a number $d\in\mathbb{N}$ of space dimensions, a BSPDE is a system of $m\in\mathbb{N}$ PDEs posed in the \emph{bulk} $\Omega \subset\mathbb{R}^d$ coupled with $n\in\mathbb{N}$ PDEs posed on the surface $\Gamma := \partial \Omega$ through either linear or non-linear coupling, see for instance \cite{Madzvamuse_2016}.} {For ease of presentation, we confine the exposition to the case $m=n=1$, even if the whole study applies to the case of arbitrary $m$ and $n$. In the time-independent case,} {let $u(\boldx)$ and $v(\boldx)$ be the bulk and surface variables obeying the following \emph{elliptic bulk-surface (BS) linear problem}:
\begin{equation}
\label{elliptic_model}
\begin{cases}
- \Delta u (\boldx) + u(\boldx) = f(\boldx), \qquad \boldx \in \Omega;\\
- \Delta_\Gamma v (\boldx) + v(\boldx) + \dfrac{\partial u}{\partial \boldnu}(\boldx)  = g(\boldx), \qquad \boldx \in \Gamma;\\
\dfrac{\partial}{\partial \boldnu} u(\boldx) = -\alpha u(\boldx) + \beta v(\boldx), \qquad \boldx \in \Gamma,\\
\end{cases}
\end{equation}
where $\alpha, \beta > 0$, while $\Delta$ and $\Delta_{\Gamma}$ denote the Laplace and Laplace-Beltrami operators respectively, and $\boldnu$ denotes the outward unit normal vector field on $\Gamma$ (see Appendix A for notations and definitions). The above model \eqref{elliptic_model} was considered in \cite{elliottranner2013finite}. }

{For a nonlinear time-dependent generalisation, let $u(\boldx,t)$ and $v(\boldx,t)$ be the bulk and surface variables obeying the following \emph{bulk-surface reaction-diffusion system} (BSRDS):
\begin{equation}
\label{parabolic_problem}
\begin{cases}
\vspace{1mm}
\dfrac{\partial u}{\partial t} (\boldx,t) - d_u\Delta u (\boldx,t) = q(u(\boldx,t)), \quad \boldx \in \Omega, \ t \in [0,T];\\
\vspace{1mm}
\dfrac{\partial v}{\partial t} (\boldx,t) - d_v\Delta_\Gamma v (\boldx,t) + \dfrac{\partial u}{\partial \boldnu} (\boldx,t) = r(u(\boldx,t), v(\boldx,t)), \quad \boldx \in \Gamma, \ t \in [0,T];\\
\dfrac{\partial u}{\partial \boldnu} (\boldx,t) = s(u(\boldx,t), v(\boldx,t)), \quad \boldx \in \Gamma, \ t \in [0,T],
\end{cases}
\end{equation}
where $T$ is the final time and $q(u), r(u,v), s(u,v)$ are possibly nonlinear functions. The model comprises several time-dependent BSPDE models existing in the literature, see for example \cite{Ratz_2014, Madzvamuse_2016, Elliott_2017, cusseddu2019}.}

The quickly growing interest toward BSPDEs arises from the numerous applications of such PDE problems in different areas, such as cellular biological systems \cite{Mackenzie_2016, Elliott_2017, cusseddu2019, Paquin_Lefebvre_2019} or fluid dynamics \cite{Bianco_2013, Lee_2015, burman2016cut}, among many other applications. Among the various state-of-the art numerical methods for the spatial discretisation of BSPDEs existing in the literature we mention finite elements \cite{elliottranner2013finite, Madzvamuse_2015, Madzvamuse_2016, Kovacs_2016}, trace finite elements \cite{gross2015trace}, cut finite elements \cite{burman2016cut} and discontinuous Galerkin methods \cite{chernyshenko2018hybrid}. 

{On the other hand, the Virtual Element Method (VEM) for bulk-only PDEs} is a recent extension of the well-known finite element method (FEM) for the numerical approximation of several classes of PDEs on flat domains \cite{beirao2013basic} or surfaces \cite{frittelli2018virtual}. The key feature of VEM is that of being a \emph{polygonal} method, i.e.  it handles elements of a quite general polygonal shape, rather than just of triangular shape \cite{beirao2013basic}. The success of virtual elements is due to several advantages arising from non-polygonal mesh generality. A non-exhaustive list of such advantages includes: (i) computationally cheap mesh pasting \cite{benkemoun2012anisotropic, chen2014memory, frittelli2018virtual}, (ii) efficient adaptive algorithms \cite{cangiani2014adaptive}, flexible approximation of the domain and in particular of its boundary \cite{dai2007n}, and the possibility of enforcing higher regularity to the numerical solution \cite{antonietti2016c, da2013arbitrary, brezzi2013virtual}, just to mention a few. Thanks to these advantages, several extensions of the original VEM for the Poisson equation \cite{beirao2013basic} were developed for numerous PDE problems, such as heat \cite{vacca2015virtual} and wave equations \cite{vacca2016virtual}, reaction-diffusion systems \cite{adak2019convergence}, Cahn-Hilliard equation \cite{antonietti2016c}, Stokes equation \cite{da_veiga_2017_stokes}, linear elasticity \cite{da2013virtual}, plate bending \cite{brezzi2013virtual}, fracture problems \cite{benedetto2016globally}, eigenvalue problems \cite{mora2015virtual} and many more.

The purpose of the present paper is to introduce a \emph{bulk-surface virtual element method} (BSVEM) for the spatial discretisation of a coupled system of BSPDEs in two space dimensions. The proposed method combines the VEM for the bulk equation(s) with the surface finite element method (SFEM) for the surface equation(s). We apply the proposed method to {(i) the linear elliptic BS Poisson problem \eqref{elliptic_model} and (ii) the BSRDS \eqref{parabolic_problem}}.

The main novelty in the present study is devoted to error analysis. In fact, the simultaneous presence of non-triangular elements and boundary approximation error (which cannot be neglected in the context of BSPDEs, {because the boundary is itself the domain of a surface PDE}) provide new numerical analysis challenges. We prove that the proposed method possesses optimal second-order convergence in space, provided the exact solution is $H^{2+1/4}$ in the bulk and $H^2$ on the surface. $H$ denotes the Hilbert space. This is slightly more than the usual requirement of $H^2$ both in the bulk and on the surface \cite{Kovacs_2016}. However, our analysis requires this slightly higher regularity assumption only \emph{in the simultaneous presence of a curved boundary $\Gamma$ and non-triangular elements close to the boundary}, which is a novel case. Otherwise, our results fall back to the known cases in the literature, see for instance \cite{Kovacs_2016} for the case of triangular BSFEM and \cite{adak2019convergence} for the case of polygonal VEM in the absence of curvature in the boundary $\Gamma$. {The extension to higher order cases would require the usage of curved elements, because otherwise the geometric error arising from the boundary approximation would dominate over the numerical error, see \cite{demlow2009higher}. This challenge will be addressed in future studies.} The proposed analysis has three by-products:

\begin{enumerate}
\item The bulk-VEM of lowest polynomial order $k=1$ possesses optimal convergence also \emph{in the presence of curved boundaries}. A first work in this direction is found in \cite{beirao2019curved}, in which the authors consider polygonal elements with a curved boundary that match the exact domain in order to {avoid errors arising from boundary approximation}. Subsequently, in \cite{bertoluzza2019}, the need of matching the exact domain was removed by introducing suitable corrections {in the method}. In the present work, instead, we obtain similar results in the low polynomial case $k=1$ by {harnessing in a novel way the geometric error estimates of BS polyhedral domains in} \cite{elliottranner2013finite}.
\item A potential alternative approach to the \emph{lifting operator} used in the analysis of SPDEs \cite{dziuk2013finite} and BSPDEs \cite{elliottranner2013finite} is the Sobolev extension operator. In fact, we prove that, if a function is $H^{2+1/4}$ instead of $H^2$ in the two-dimensional bulk, the Sobolev extension retains optimal approximation properties of lifting. Moreover, the Sobolev extension of a function has the property of preserving its $W^{m,p}$ class, while its lift does not {because it is not $\mathcal{C}^k$ for any positive integer $k$}. This property, which is crucial in our analysis, is potentially beneficial to the {error analysis of bulk-only or BSVEM approximation of more general PDEs, where the boundary curvature was not accounted for, see for instance \cite{beirao2013basic, vacca2015virtual, mora2015virtual, vacca2016virtual, da_veiga_2017_stokes, frittelli2018virtual}.}
\item We construct a special \emph{inverse trace operator} that we use for accounting for general boundary conditions. Like the standard inverse trace operator, our inverse trace maps a function $v \in H^1(\Gamma)$ to a function $v_B \in H^1(\Omega)$ such that $\|v_B\|_{H^1(\Omega)} \leq C\|v\|_{H^1(\Gamma)}$, with $C$ depending on the domain $\Omega$. Our inverse trace has the stronger property that $\|v_B\|_{L^2(\Omega)} \leq C\|v\|_{L^2(\Gamma)}$ and $|v_B|_{H^1(\Omega)} \leq C\|v\|_{H^1(\Gamma)}$. This means that the proposed operator preserves both $L^2$ and $H^1$ norms up to the same multiplicative constants, i.e. it is an \emph{$L^2$-preserving inverse trace}.
\end{enumerate}
The proposed method has all the benefits of polygonal meshes, two of which will be illustrated in the present work. First, the usage of suitable polygons drastically reduces the computational complexity of matrix assembly on equal meshsize in comparison to the triangular BSFEM. Similar results are obtained in the literature through other methods, such as trace FEMs \cite{gross2015trace} or cut FEMs \cite{burman2016cut}. Second, a curved portion of the boundary can be approximated with a single element with many edges. The BSVEM lends itself to other advantages due to its polygonal nature, such as efficient algorithms for adaptivity or mesh pasting, see for instance \cite{cangiani2014adaptive}. These aspects will be addressed in future studies.

Hence, the structure of our paper is as follows.
In Section \ref{sec:bspdes}, {we elaborate on the weak formulations, existence and regularity for problems \eqref{elliptic_model} and \eqref{parabolic_problem}}.
In Section \ref{sec:bsvem}, we introduce polygonal BS meshes, analyse geometric error, define suitable function spaces, analyse their approximation properties and present the spatial discretisation of the considered BSPDE problems.
In Section \ref{sec:convergence_analysis}, we carry out the convergence error analysis for the parabolic case, the main result being optimal second-order spatial convergence in the $L^2$ norm, both in the bulk and on the surface.
In Section \ref{sec:time_disc} we present an IMEX-Euler time discretisation of the parabolic problem.
In Section \ref{sec:mesh_advantage} we show that polygonal meshes can be exploited to significantly reduce the computational time of the matrix assembly.
In Section \ref{sec:numerical_examples} we illustrate our findings through {three numerical examples. The first example shows (i) optimised matrix assembly and (ii) optimal convergence for the elliptic problem \eqref{elliptic_model}. The second example shows optimal convergence in space and time for the parabolic problem \eqref{parabolic_problem} for the linear case. The third example compares the BSFEM- and BSVEM-IMEX Euler approximations of the wave pinning RDS considered in \cite{cusseddu2019}.}
In Appendix A we provide some basic definitions and results required in our analysis.
In Appendix B we report some lengthy proofs involved in Section \ref{sec:convergence_analysis}.

\section{Weak formulations of BSPDEs, existence and regularity}
\label{sec:bspdes}
In this section we state the weak formulations of the elliptic \eqref{elliptic_model} and parabolic problems \eqref{parabolic_problem}, respectively. For problem \eqref{elliptic_model}, we multiply the first two equations of \eqref{elliptic_model} by two test functions $\varphi \in H^1(\Omega)$ and $\psi \in H^1(\Gamma)$, respectively, then we apply Green's formula in the bulk $\Omega$ and on the one-dimensional manifold $\Gamma$ \cite{dziuk2013finite}. We obtain the following formulation: find $u \in H^1(\Omega)$ and $v \in H^1(\Gamma)$ such that
\begin{equation}
\label{elliptic_problem_weak_form_intermediate}
\begin{cases}
\vspace{2mm}
\displaystyle \intomega \nabla u \cdot \nabla \varphi + \intomega u\varphi = \intomega f\varphi + \intgamma \dfrac{\partial u}{\partial\boldnu}\varphi;\\
\displaystyle \intgamma \nablagamma v\cdot\nablagamma \psi + \intgamma v\psi + \intgamma \dfrac{\partial u}{\partial \boldnu} \psi = \intgamma g\psi,
\end{cases}
\end{equation}
for all $\varphi \in H^1(\Omega)$ and $\psi \in H^1(\Gamma)$. By using the third equation of \eqref{elliptic_model} in \eqref{elliptic_problem_weak_form_intermediate}, we obtain the following weak formulation: find $u\in H^1(\Omega)$ and $v\in H^1(\Gamma)$ such that 
\begin{equation}
\label{elliptic_problem_weak_form}
\begin{cases}
\vspace{2mm}
\displaystyle \intomega \nabla u \cdot \nabla \varphi + \intomega u\varphi + \intgamma (\alpha u -\beta v)\varphi = \intomega f\varphi;\\
\displaystyle \intgamma \nablagamma v\cdot\nablagamma \psi + \intgamma (-\alpha u + (\beta+1)v) \psi = \intgamma g\psi,
\end{cases}
\end{equation}
for all $\varphi \in H^1(\Omega)$ and $\psi \in H^1(\Gamma)$. 

By reasoning similarly to the elliptic problem \eqref{elliptic_model}, we obtain the following weak formulation of the parabolic problem \eqref{parabolic_problem}: find $u \in L^\infty([0,T]; H^1(\Omega))$ and $v\in L^\infty([0,T]; H^1(\Gamma))$ such that
\begin{equation}
\label{weak_formulation}
\begin{cases}
\vspace{2mm}
\displaystyle\intomega \ut \varphi +  d_u\intomega \nabla u \cdot \nabla \varphi = \intomega q(u)\varphi + \intgamma s(u,v)\varphi;\\
\displaystyle\intgamma \vt \psi +  d_v\intgamma \nablagamma v \cdot \nablagamma \psi + \intgamma s(u,v)\psi = \intgamma r(u,v)\psi,
\end{cases}
\end{equation} 
for all $\varphi \in L^\infty([0,T]; H^1(\Omega))$ and $\psi \in L^\infty([0,T]; H^1(\Gamma))$. The following theorem contains existence, uniqueness and regularity results for the weak problem \eqref{elliptic_problem_weak_form}.

\begin{theorem}[Existence, uniqueness and regularity for problem \eqref{elliptic_problem_weak_form}]
\label{thm:existence_elliptic}
If $\Gamma$ is a $\mathcal{C}^3$ surface, $f\in L^2(\Omega)$ and $g\in L^2(\Gamma)$, then the weak elliptic problem \eqref{elliptic_problem_weak_form} has a unique solution $(u,v) \in H^2(\Omega) \times H^2(\Gamma)$ that fulfils
\begin{align}
\label{regularity_elliptic_u_v}
\|(u,v)\|_{H^2(\Omega) \times H^2(\Gamma)} \leq C \|(f,g)\|_{L^2(\Omega) \times L^2(\Gamma)};\\
\label{regularity_elliptic_tr_u}
\|\Tr(u)\|_{H^{3/2}(\Omega)} \leq C \|(f,g)\|_{L^2(\Omega) \times L^2(\Gamma)},
\end{align}
where $C>0$ is a constant that depends on $\alpha, \beta$ and $\Omega$.
\begin{proof}
See \cite{elliottranner2013finite} for estimate \eqref{regularity_elliptic_u_v}. Estimate \eqref{regularity_elliptic_tr_u} follows from \eqref{regularity_elliptic_u_v} by using the trace inequality \eqref{trace_inequality}.
\end{proof}
\end{theorem}

The problem of existence, uniqueness and regularity for the parabolic problem \eqref{weak_formulation} is much more complicated and strictly depends on the nature of the kinetics $q(\cdot)$, $r(\cdot)$ and of the coupling kinetics $s(\cdot)$. For the remainder of this work we will adopt the following set of assumptions.

\begin{assumption}[Existence, uniqueness and regularity for problem \eqref{weak_formulation}]
\label{thm:assumptions_parabolic}
We assume that:
\begin{itemize}
\item $\Gamma$ is a $\mathcal{C}^3$ surface, $q,r,s$ are $\mathcal{C}^2$ and globally Lipschitz functions.
\item The initial datum $(u_0,v_0)$ fulfils $u_0\in H^2(\Omega)$, $\Tr(u_0) \in H^2(\Gamma)$ and $v_0\in H^2(\Gamma)$.
\item There exists a unique solution $(u,v)$ that fulfils
\begin{align}
\|(u,\dot{u})\|_{L^\infty([0,T]; H^{2+1/4}(\Omega))} + &\|(\Tr(u),\Tr(\dot{u}),v,\dot{v})\|_{L^\infty([0,T]; H^2(\Gamma))}\\
\notag
\leq &C \exp(T)\left(\|u_0\|_{H^{2+1/4}(\Omega)} + \|(\Tr(u_0),v_0)\|_{H^2(\Gamma)}\right),
\end{align}
where $T>0$ is the final time and $C>0$ depends on $\Omega$, $\|q\|_{\mathcal{C}^2(\mathbb{R})}$, $\|r\|_{\mathcal{C}^2(\mathbb{R}^2)}$ and $\|s\|_{\mathcal{C}^2(\mathbb{R}^2)}$.
\end{itemize}
\end{assumption}

In many applications, assuming globally Lipschitz kinetics is too restrictive and an ad-hoc analysis is required, see for instance \cite{Elliott_2017}. However, there are notable examples of BSRDS with globally Lipschitz kinetics, such as the wave pinning model studied in \cite{cusseddu2019} and considered in the numerical example in Section \ref{sec:experiment_parabolic}. From here onwards, we shall assume that the weak parabolic problem \eqref{weak_formulation} has a unique and sufficiently regular solution.

\section{The Bulk-Surface Virtual Element Method}
\label{sec:bsvem}
\noindent
In this section we will introduce the Bulk-Surface Virtual Element Method (BSVEM) through the following steps:
\begin{itemize}
\item describe the polygonal BS meshes that will be used in the method (Subsection \ref{sec:bsmeshes});
\item quantify the geometric error arising from polygonal approximation of BS domains (Subsection \ref{sec:variationalcrime});
\item introduce the discrete function spaces and bilinear forms used in the method and their approximation properties (Subsections \ref{sec:vem_space_bulk}-\ref{sec:ingredients_surface});
\item present the spatially discrete formulations of the elliptic- and parabolic problems \eqref{elliptic_model} and \eqref{parabolic_problem}, respectively (Subsection \ref{sec:spatially_discrete_formulations}).
\end{itemize}

\subsection{Polygonal bulk-surface meshes}
\label{sec:bsmeshes}
\begin{figure}
\begin{center}
\subfigure[Illustration of the \emph{bulk} $\Omega$, enclosed by the \emph{surface} $\Gamma$, the \emph{narrow band} $U_\delta$ and the Fermi stripe $U$.]{\input{exact_domain.tex}\label{fig:exact_domain}}
\hspace*{1cm}
\subfigure[Illustration of the \emph{discrete bulk} $\Omega_h$, enclosed by the \emph{discrete surface} $\Gamma_h$, the \emph{discrete narrow band} $\Omega_B$ and the Fermi stripe $U$.]{\input{mesh_portions.tex}\label{fig:discrete_domain}}
\end{center}
\caption{Illustration of continuous domain, discrete domain and related notations.}
\label{fig:mesh_illustration}
\end{figure}
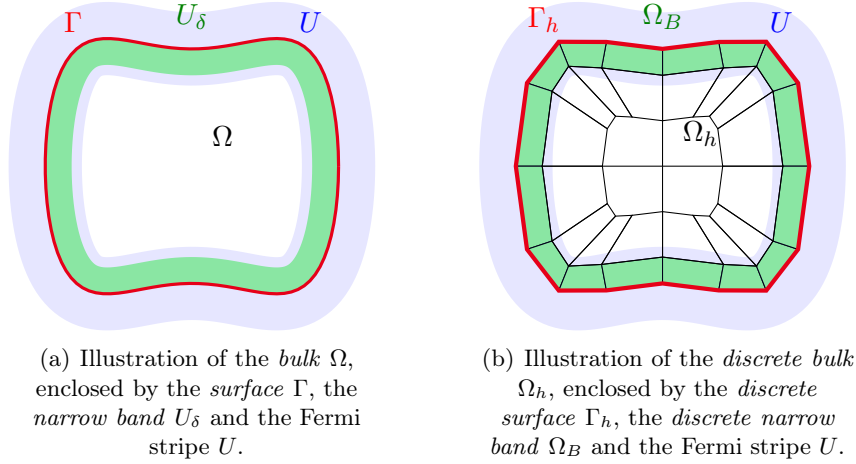

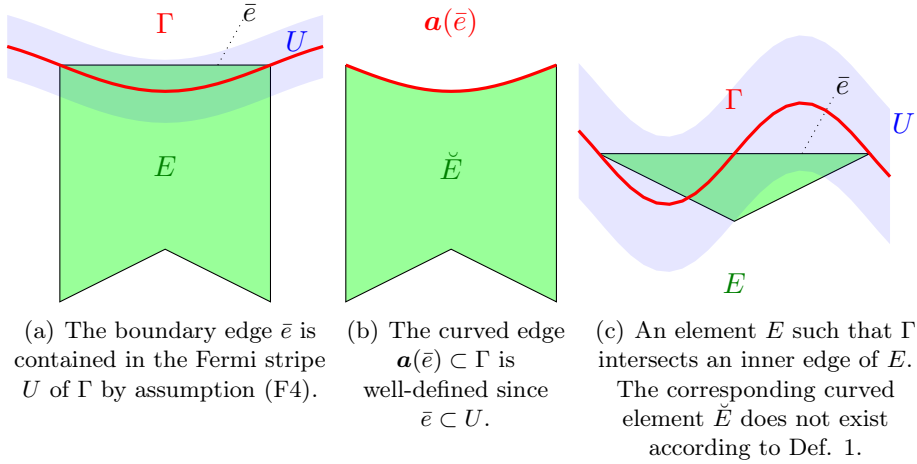
\begin{figure}
\begin{center}
\subfigure[The boundary edge $\bar{e}$ is contained in the Fermi stripe $U$ of $\Gamma$ by assumption \ref{B4}.]{\input{element.tex}}
\subfigure[The curved edge $\bolda(\bar{e}) \subset \Gamma$ is well-defined since $\bar{e}\subset U$.]{\input{exact_element.tex}}
\subfigure[An element $E$ such that $\Gamma$ intersects an inner edge of $E$. The corresponding curved element $\breve{E}$ does not exist according to Def. \ref{def:boundary_elements}.]{\input{element_intersection.tex}}
\end{center}
\caption{Construction of the exact element $\breve{E}$ corresponding to the polygonal element $E$ according to Def. \ref{def:boundary_elements}. {If $h$ is sufficiently small (depending on the curvature of $\Gamma$ and the mesh regularity parameter $\gamma_1$) then $\Gamma$ cannot intersect any inner edges and curved elements are well-defined (subplots (a)-(b)). Otherwise, if $\Gamma$ intersects an inner edge of $E$, $\breve{E}$ is not well-defined (subplot (c)).}}
\label{fig:exact_element}
\end{figure}

Let $h>0$ be a positive number called \emph{meshsize} and let $\Omega_h = \cup_{E\in\mathcal{E}_h} E$ be a polygonal approximation of the bulk $\Omega$, where $\mathcal{E}_h$ is a set of non-degenerate polygons. The polygonal bulk $\Omega_h$ automatically induces a piecewise linear approximation $\Gamma_h$ of $\Gamma$, defined by $\Gamma_h = \partial \Omega_h$, exactly as in the case of triangular meshes, see for example \cite{elliottranner2013finite}. Notice that we can write $\Gamma_h = \cup_{F \in \mathcal{F}_h} F$, where $\mathcal{F}_h$ is the set of the edges of $\Omega_h$ that lie on $\Gamma_h$. We assume that:
\begin{enumerate}[label=\textnormal{(F\arabic*)}]
\item the diameter of each element $E\in\mathcal{E}_h$ does not exceed $h$;\label{B1}
\item for any two distinct elements $E_1$, $E_2 \in \mathcal{E}_h$, their intersection $E_1\cap E_2$ is either empty, or a common vertex, or a common edge.
\item all nodes of $\Gamma_h$ lie on $\Gamma$;
\item every edge $F\in\mathcal{F}_h$ is contained in the Fermi stripe $U$ of $\Gamma$ (see Fig. \ref{fig:mesh_illustration} for an illustration).\label{B4}
\end{enumerate}
\begin{enumerate}[label=\textnormal{(V\arabic*)}]
\item there exists $\gamma_1 > 0$ such that every $E \in \mathcal{E}_h$ is star-shaped with respect to a ball of radius $\gamma_1 h_E$, where $h_E$ is the diameter of $E$;\label{A1} 
\item there exists $\gamma_2 > 0$ such that for all $E \in \mathcal{E}_h$, the distance between any two nodes of $E$ is at least $\gamma_2 h_E$.\label{A2}
\end{enumerate}
Assumptions \ref{B1}-\ref{B4} are standard in the SFEM literature, see for instance \cite{dziuk2013finite}, while assumptions \ref{A1}-\ref{A2} are standard in the VEM literature, see for instance \cite{beirao2013basic}. The combined assumptions \ref{B1}-\ref{A2} will prove sufficient in our BS setting. In the following definitions and results we provide the necessary theory for estimating the geometric error arising from the boundary approximation.
\begin{definition}[Essentials of polygonal BS meshes]
\label{def:boundary_elements}
An edge $\bar{e}$ of any element $E\in\mathcal{E}_h$ is called a \emph{boundary edge} if $\bar{e} \subset \Gamma_h$, otherwise $\bar{e}$ is called an \emph{inner edge}. Let $\mathcal{BE}(E)$ and $\mathcal{IE}(E)$ be the sets of boundary and inner edges of $E$, respectively. An element $E\in\mathcal{E}_h$ is called an \emph{external element} if it has at least one boundary edge, otherwise $E$ is called an \emph{internal element}.
Let $\Omega_B$ be the \emph{discrete narrow band} defined as the union of the external elements of $\Omega_h$ as illustrated in Fig. \ref{fig:discrete_domain}.
From Assumption \ref{B4}, for any boundary edge $\bar{e}$, we have that $\bolda(\bar{e}) \subset \Gamma$, where $\bolda$ is the normal projection defined in Lemma \ref{lmm:fermi}. Hence, {for sufficiently small $h$ and for all $E\in\mathcal{E}_h$, it is possible to} define the exact element $\breve{E}$ as the compact set enclosed by the edges
$\mathcal{IE}(E) \cup \{\bolda(\bar{e}) | \bar{e} \in \mathcal{BE}(E)\}$,
see Fig. \ref{fig:exact_element} for an illustration.
\end{definition}

\begin{remark}[Properties of polygonal BS meshes]
For any BS mesh $(\Omega_h, \Gamma_h)$ of meshsize $h>0$ it holds that:
\begin{itemize}
\item for sufficiently small $h>0$, the discrete narrow band $\Omega_B$ is contained in the Fermi stripe $U$ as shown in Fig. \ref{fig:discrete_domain} (blue colour);
\item the collection of all exact elements is a coverage of the exact bulk $\Omega$, that is $\cup_{E\in\mathcal{E}_h} \breve{E} = \Omega$.
\end{itemize}
\end{remark}

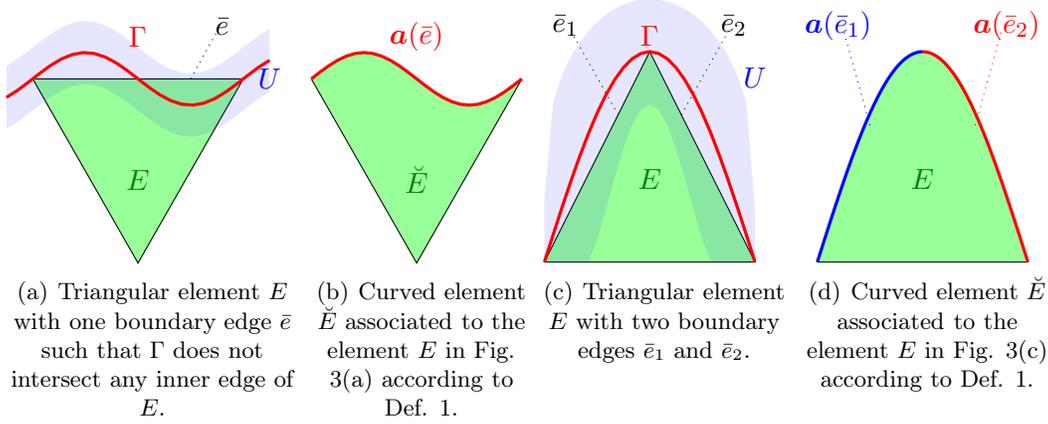
\begin{figure}
\begin{center}
\subfigure[Triangular element $E$ with one boundary edge $\bar{e}$ such that $\Gamma$ does not intersect any inner edge of $E$.]{\input{element_smooth.tex}}
\subfigure[Curved element $\breve{E}$ associated to the element $E$ in {Fig. 3(a)} according to Def. \ref{def:boundary_elements}.]{\input{element_smooth_exact.tex}}
\subfigure[Triangular element $E$ with two boundary edges $\bar{e}_1$ and $\bar{e}_2$.]{\input{element_corner.tex}}
\subfigure[Curved element $\breve{E}$ associated to the element $E$ in {Fig. 3(c)} according to Def. \ref{def:boundary_elements}.]{\input{element_corner_exact.tex}}\\
\end{center}
\caption{{On the mapping $G_E$ between a triangular element $E$ and its curved counterpart $\breve{E}$. In the case of one boundary edge (subplots (a)-(b)), a $\mathcal{C}^2$ homeomorphism $G_E: E\rightarrow\breve{E}$ is known to exist, see \cite{elliottranner2013finite}. With two adjacent boundary edges (subplots (c)-(d)) the mapping $G_E$ cannot be smooth, because $G_E$ maps the non-smooth curve $\bar{e}_1 \cup \bar{e}_2$ onto the smooth curve $\bolda(\bar{e}_1) \cup \bolda(\bar{e}_2)$. The approach proposed in Lemma \ref{lmm:homeomorphism} solves this issue.}}
\label{fig:element_corner}
\end{figure}

Let $N \in\mathbb{N}$ and let $\boldx_i$, $i=1,\dots,N$, be the nodes of $\Omega_h$, which can be ordered in an arbitrary way. However, if $\Omega$ has a rectangular shape and the nodes are ordered along a Cartesian grid, the matrices associated with the method will have a block-tridiagonal structure. Let $M \in\mathbb{N}$, $M < N$ and assume that the nodes of $\Gamma_h$ are $\boldx_k$, $k=1,\dots,M$, i.e. the first $M$ nodes of $\Omega_h$. Throughout the paper we need the following \emph{reduction matrix} $R \in\mathbb{R}^{N\times M}$ defined as $R := [I_M; 0]$, where $I_M$ is the $M\times M$ identity matrix. Equivalently, $R = (r_{ik})$ is defined as
\begin{equation}
\label{reduction_matrix}
r_{ik} = 
\begin{cases}
\delta_{ik} \qquad &\text{if}\ i=1,\dots,M;\\
0 \qquad &\text{if}\ i=M+1,\dots,N,
\end{cases}
\end{equation}
for all $k=1,\dots,M$, where $\delta_{ik}$ is the Kronecker symbol. The reduction matrix $R$ fulfils the following two properties:
\begin{itemize}
\item For $\boldv \in \mathbb{R}^{N}$, $R^T\boldv \in \mathbb{R}^M$ is the vector with the first $M$ entries of $\boldv$;
\item For $\boldw \in \mathbb{R}^M$, $R\boldw \in \mathbb{R}^N$ is the vector whose first $M$ entries are those of $\boldw$ and the other $N-M$ entries are $0$.
\end{itemize}
In what follows, we will use the matrix $R$ for an optimised implementation of the BSVEM.

\subsection{Variational crime}
\label{sec:variationalcrime}
We now consider the geometric error due to the boundary approximation. Since the surface variational crime in surface finite elements is well-understood \cite{dziuk2013finite}, we will mainly focus on the variational crime in the bulk. To this end, it is useful to analyse the relation between any element $E\in\mathcal{E}_h$ and its exact counterpart $\breve{E}$, see Def. \ref{def:boundary_elements}. For the special case of triangular meshes with at most one boundary edge per element (Fig. \ref{fig:element_corner}(a)-(b)), there exists {a $\mathcal{C}^2$ homeomorphism $G_E: E \rightarrow \breve{E}$ that is quadratically} close to the identity with respect to the meshsize, see \cite{elliottranner2013finite}. 
{Instead, when an element $E$ has adjacent, non-collinear boundary edges (which can occur even in triangular meshes, see Fig. \ref{fig:element_corner}(c)-(d)), such a smooth homeomorphism does not exist. In fact, a smooth mapping cannot map the adjacent boundary edges of $E$ - a curve with corners - onto a portion of the smooth curve $\Gamma$.} However, in the following result we show the existence of a {homeomorphism} between $E$ and $\breve{E}$ with slightly weaker regularity, which is sufficient for our purposes.

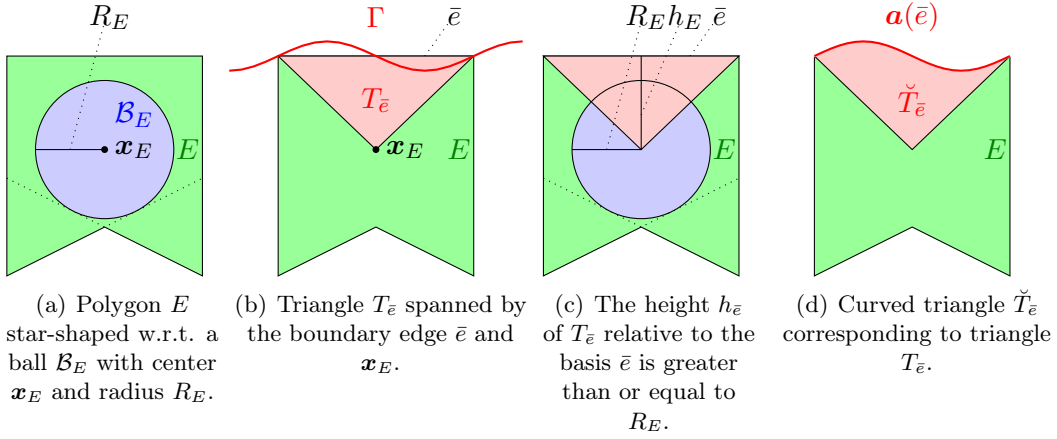
\begin{figure}
\begin{center}
\subfigure[Polygon $E$ star-shaped w.r.t. a ball $\mathcal{B}_E$ with center $\boldx_E$ and radius $R_E$.]{\label{fig:proof_step_1}\input{lemma_proof_1.tex}}
\subfigure[Triangle $T_{\bar{e}}$ spanned by the boundary edge $\bar{e}$ and $\boldx_E$.]{\label{fig:proof_step_2}\input{lemma_proof_2.tex}}
\subfigure[The height $h_{\bar{e}}$ of $T_{\bar{e}}$ relative to the basis $\bar{e}$ is greater than or equal to $R_E$.]{\label{fig:proof_step_3}\input{lemma_proof_3.tex}}
\subfigure[Curved triangle $\breve{T}_{\bar{e}}$ corresponding to triangle $T_{\bar{e}}$.]{\label{fig:proof_step_4}\input{lemma_proof_4.tex}}
\end{center}
\caption{Steps of the construction of $G$ as described in the proof of Lemma \ref{lmm:homeomorphism}.}
\end{figure}

\begin{lemma}[Parameterisation of the exact geometry]
\label{lmm:homeomorphism}
{Let $h$ be sufficiently small (depending on the curvature of $\Gamma$ and the mesh regularity parameter $\gamma_1$) and} let $\mathcal{E}_h$ fulfil assumptions \ref{B1}-\ref{A2}. There exists a homeomorphism $G:\Omega_h \rightarrow \Omega$ such that $G \in W^{1,\infty}(\Omega_h)$ and
\begin{align}
\label{estimate_1}
&G|_{\Gamma_h} = \bolda|_{\Gamma_h};\\
\label{estimate_2}
&G|_{\Omega_h \setminus \Omega_B} = Id;\\
\label{estimate_3}
&\|JG - Id\|_{L^\infty(\Omega_B)} \leq Ch;\\
\label{estimate_4}
&\|\det(JG) - 1\|_{L^\infty(\Omega_B)} \leq Ch;\\
\label{estimate_5}
&\|G - Id\|_{L^\infty(\Omega_B)} \leq Ch^2,
\end{align}
where $\bolda$ is the normal projection defined in Lemma \ref{lmm:fermi}, $JG$ is the Jacobian of $G$ and $C$ is a constant that depends on $\Gamma$ and the constants $\gamma_1$, $\gamma_2$ are those considered in Assumptions \ref{A1}-\ref{A2}.
\begin{proof}
Let $E\in \mathcal{E}_h$. By Assumption \ref{A2}, $E$ is star-shaped with respect to a ball $\mathcal{B}_E$ of diameter $R_E \geq \gamma_2 h_E$, let $\boldx_E$ be the center of such ball {(Fig. \ref{fig:proof_step_1})}. For $\bar{e} \in\mathcal{BE}(E)$ let $T_{\bar{e}}$ be the triangle spanned by $\bar{e}$ and $\boldx_E$ {(Fig. \ref{fig:proof_step_2})}. Let us now consider the collection of all the $T_{\bar{e}}$'s defined as follows:
$\mathcal{BT}_h := \{T_{\bar{e}} | \bar{e}\in E, \ E \in \mathcal{E}_h\}$.
We need to prove that $\mathcal{BT}_h$ is quasi-uniform, i.e we need to prove that for all $E\in\mathcal{E}_h$ and $\bar{e}\in\mathcal{BE}(E)$, the triangle $T_{\bar{e}}$ has an inscribed ball of diameter greater or equal to $\bar{\gamma} h_E$, where the constant $\bar{\gamma} > 0$ depends on $\gamma_1$ and $\gamma_2$, only. To this end, if $h_{\bar{e}}$ is the height of $T_{\bar{e}}$ relative to the basis $\bar{e}$ {(Fig. \ref{fig:proof_step_3})}, then we have
\begin{align}
|h_{\bar{e}}| \geq R_E \geq \gamma_1 h_E, \qquad \text{from Assumption } \ref{A1};\\
|\bar{e}| \geq \gamma_2 h_E, \qquad \text{from Assumption } \ref{A2}.
\end{align}
In addition, since no edge of $T_{\bar{e}}$ is longer than $h_E$, our claim follows.\\
Let now $\breve{T}_{\bar{e}}$ be the curved triangle corresponding to $T_{\bar{e}}$ {(Fig. \ref{fig:proof_step_4})}, {which again is well-defined by assuming $h$ sufficiently small depending on $\Gamma$ and $\bar{\gamma}$, which in turn depends on $\gamma_1$ and $\gamma_2$.} Since the triangulation $\mathcal{BT}_h$ is quasi-uniform, then, from \cite[Proposition 4.7 and its proof]{elliottranner2013finite} there exists a diffeomorphism $G_{\bar{e}}: T_{\bar{e}} \rightarrow \breve{T}_{\bar{e}}$ such that
\begin{align}
\label{extension_of_normal_projection}
&G_{\bar{e}}(\boldx) = \bolda(\boldx), \qquad \forall \boldx \in \bar{e};\\
\label{continuity_across_edges}
&G_{\bar{e}}(\boldx) = \boldx, \qquad \forall \boldx \in \partial T_{\bar{e}} \setminus \bar{e};\\
\label{homeomorphism_regularity_1}
&\|JG_{\bar{e}} - Id\|_{L^\infty(T_{\bar{e}})} \leq Ch;\\
\label{homeomorphism_regularity_2}
&\|\det(JG)_{\bar{e}} - 1\|_{L^\infty(T_{\bar{e}})} \leq Ch;\\
\label{homeomorphism_regularity_3}
&\|G_{\bar{e}} - Id\|_{L^\infty(T_{\bar{e}})} \leq Ch^2,
\end{align} 
where $C$ depends only on $\bar{\gamma}$, which in turn depends only on $\gamma_1$ and $\gamma_2$. We are ready to construct the mapping $G:\Omega_h \rightarrow\Omega$ as follows:
\begin{equation}
\label{global_homeomorphism}
G(\boldx) =
\begin{cases}
G_{\bar{e}}(\boldx), \qquad \text{if} \ \boldx \in T_{\bar{e}} \ \text{for some} \ \bar{e} \in \mathcal{BE}(E), \ E\in\mathcal{E}_h;\\
\boldx, \qquad \text{otherwise}.
\end{cases}
\end{equation}
Property \eqref{continuity_across_edges} can now be rephrased by saying that $G$, restricted to any inner edge of $\Omega_h$, is the identity. In addition, since all the $G_{\bar{e}}$'s are homeomorphisms, we obtain that $G$ is a homeomorphism between $\Omega_h$ and $\Omega$. Finally, from \eqref{extension_of_normal_projection}, \eqref{homeomorphism_regularity_1}, \eqref{homeomorphism_regularity_2} and \eqref{global_homeomorphism} we obtain the desired estimates \eqref{estimate_1}-\eqref{estimate_4}.
\end{proof}
\end{lemma}

\begin{remark}[Virtual elements for bulk-only PDEs]
Lemma \ref{lmm:homeomorphism} has an important consequence in the analysis of boundary approximation for VEMs for the case of bulk PDEs. For triangular finite elements, boundary approximation is a well-understood topic, see for instance \cite{ciarlet2002}. For VEMs, the first work in this direction is \cite{beirao2019curved}, in which a VEM on exact curved polygons is considered, in order to take out the geometric error. However, {in the lowest-order}  VEM on {polygonal domains}, it is {empirically} known that the geometric error does not prevent optimality, {as discussed in} \cite{bertoluzza2019}. To the best of the authors' knowledge, the present study provides, as a by-product, the first rigorous proof of this fact.
\end{remark}

Thanks to Lemma \ref{lmm:homeomorphism} it is possible to define bulk- and surface-lifting operators.
\begin{definition}[Bulk- and surface-lifting operators]
Given $V:\Omega_h\rightarrow \mathbb{R}$ and $W:\Gamma_h \rightarrow\mathbb{R}$, their \emph{lifts} are defined by $V^\ell := V \circ G^{-1}$ and $W^\ell := W \circ G^{-1}$,
respectively. Conversely, given $v:\Omega \rightarrow\mathbb{R}$ and $w:\Gamma \rightarrow\mathbb{R}$, their \emph{inverse lifts} are defined by
$v^{-\ell} := v \circ G$ and $w^{-\ell} := w \circ G$, respectively, with $G:\Omega_h \rightarrow\Omega$ being the mapping defined in Lemma \ref{lmm:homeomorphism}.
\end{definition}
Lemma \ref{lmm:homeomorphism} also enables us to show the equivalence of Sobolev norms under lifting as illustrated next.
\begin{lemma}[Equivalence of norms under lifting]
There exists two constants $c_2 > c_1 > 0$ depending on $\Gamma$ and $\gamma_2$ such that, for all $V: \Omega_h \rightarrow \mathbb{R}$ and for all $W: \Gamma_h \rightarrow \mathbb{R}$,
\begin{align}
\label{equivalence_omega_l2}
&c_1 \|V^\ell\|_{L^2(\Omega_h^\ell)} \hspace*{3mm} \leq \|V\|_{L^2(\Omega_h)}  \hspace*{3mm} \leq c_2\|V^\ell\|_{L^2(\Omega_h^\ell)};\\
\label{equivalence_omega_h1}
&c_1 |V^\ell|_{H^1(\Omega_h^\ell)} \leq |V|_{H^1(\Omega_h)} \leq c_2|V^\ell|_{H^1(\Omega_h^\ell)};\\
\label{equivalence_gamma_l2}
&c_1 \|W^\ell\|_{L^2(\Gamma)}  \hspace*{2mm} \leq \|W\|_{L^2(\Gamma_h)}  \hspace*{2mm} \leq c_2\|W^\ell\|_{L^2(\Gamma)};\\
\label{equivalence_gamma_h1}
&c_1 |W^\ell|_{H^1(\Gamma)} \leq |W|_{H^1(\Gamma_h)} \leq c_2|W^\ell|_{H^1(\Gamma)}.
\end{align}
\begin{proof}
Estimates \eqref{equivalence_omega_l2}-\eqref{equivalence_omega_h1} are found by using the map $G$ introduced in Lemma \ref{lmm:homeomorphism} in the proof of \cite[Proposition 4.9]{elliottranner2013finite}. A proof of \eqref{equivalence_gamma_l2}-\eqref{equivalence_gamma_h1} can be found in \cite[Lemma 4.2]{dziuk2013finite}.
\end{proof}
\end{lemma}

We are ready to estimate the effect of lifting on bulk and surface integrals.
\begin{lemma}[Geometric error of lifting]
\label{lmm:geometric_errors_bilinear forms}
If $u, \varphi\in H^{1}(\Omega)$, then
\begin{align}
\label{lifting_error_bulk_nabla}
&\left |\intomega \nabla u \cdot \nabla \varphi - \intomegah \nabla u^{-\ell} \cdot \nabla \varphi^{-\ell} \right| \leq Ch|u|_{H^1(\Omega_B^\ell)}|\varphi|_{H^1(\Omega_B^\ell)},\\
\label{lifting_error_bulk_mass}
&\left |\intomega u \varphi - \intomegah u^{-\ell} \varphi^{-\ell} \right| \leq Ch \|u\|_{L^{2}(\Omega_B^\ell)}\|\varphi\|_{L^{2}(\Omega_B^\ell)},
\end{align} 
where $C$ depends on $\Gamma$, $\gamma_1$ and $\gamma_2$. If $v,\psi\in H^1(\Gamma)$, then
\begin{align}
\label{lifting_error_surf_nabla}
&\left|\intgamma \nablagamma v \cdot \nablagamma \psi - \intgammah \nablagammah v^{-\ell} \cdot \nablagammah \psi^{-\ell} \right| \leq Ch^2|v|_{H^1(\Gamma)} |\psi|_{H^1(\Gamma)};\\
\label{lifting_error_surf_mass}
&\left|\intgamma v \psi - \intgammah v^{-\ell} \psi^{-\ell} \right| \leq Ch^2\|v\|_{L^2(\Gamma)} \|\psi\|_{L^2(\Gamma)},
\end{align}
where $C$ depends on $\Gamma$, $\gamma_1$ and $\gamma_2$.
\begin{proof}
To prove \eqref{lifting_error_bulk_nabla}-\eqref{lifting_error_bulk_mass} it is sufficient to use the bulk geometric estimates \eqref{estimate_2}-\eqref{estimate_4} in the proof of \cite[Lemma 6.2]{elliottranner2013finite}. A proof of \eqref{lifting_error_surf_nabla}-\eqref{lifting_error_surf_mass} can be found in \cite{dziuk2013finite}.
\end{proof}
\end{lemma}

Lemma \ref{lmm:geometric_errors_bilinear forms} allows to rephrase integrals on the exact domain as integrals on the discrete domain, and vice-versa, up to a small error that is $O(h)$ in the bulk and $O(h^2)$ on the surface.

\begin{remark}[Preservation of regularity under lifting]
For $E\in\mathcal{E}_h$ the inverse lift of an $H^2(\breve{E})$ function is not, in general, $H^2(E)$, cf. Remark \ref{rmk:regularity_normal_projection} and Lemma \ref{lmm:homeomorphism}. This problem does not arise in the context of triangular BSFEM in the presence of curved boundaries, see for instance \cite{Kovacs_2016}, or bulk virtual elements in the absence of curved boundaries, see for instance \cite{adak2019convergence}. In the context of bulk- or bulk-surface VEM in the presence of curved boundaries, our analysis requires \emph{full $H^2$-regularity of the exact solution {mapped on} the polygonal domain}. {Because $u^{-\ell}$ does not retain the regularity class of $u$}, we need {an alternative {mapping} of $u$ on $\Gamma_h$. To this end, we recall the following Theorem.}

\begin{theorem}[Sobolev extension theorem]
\label{thm:sobolev_extension_theorem}
Assume that $\Omega$ has a Lipschitz boundary $\Gamma$, let $r\in\mathbb{N}$ and $p\in [1,+\infty]$. Then, for any function $u\in W^{r,p}(\Omega)$, there exists an extension $\tilde{u} \in W^{r,p}(\mathbb{R}^2)$ such that $\tilde{u}_{|\Omega} = u$ and
\begin{equation}
\label{sobolev_extension}
\|\tilde{u}\|_{W^{r,p}(\mathbb{R}^2)} \leq C\|u\|_{W^{r,p}(\Omega)},
\end{equation}
where $C$ depends on $\Omega$ and $r$.
\begin{proof}
See \cite{Stein_1971}.
\end{proof}
\end{theorem}

{We are now able to approximate the exact solution $u$ through $\tilde{u}_{|\Omega_h}$, that is the restriction to $\Omega_h$ of the Sobolev extension $\tilde{u}$ of the exact solution $u$. We present the following variant of Lemma \ref{lmm:geometric_errors_bilinear forms} in order to quantify the geometric error of the Sobolev extension.}
\end{remark}

\begin{lemma}[Geometric error of Sobolev extension]
If $0 < \gamma < 1$, there exist $C>0$ and $C_\gamma > 0$ such that
\begin{align}
\label{error_between_lift_and_extension_l2}
&\|\tilde{u} - u^{-\ell}\|_{L^2(\Omega_h)} \leq Ch^2\|u\|_{H^{1+3/4}(\Omega)}, \quad \forall\ u\in H^{1+3/4}(\Omega);\\
\label{error_between_lift_and_extension}
&|\tilde{u} - u^{-\ell}|_{H^1(\Omega_h)} \leq Ch^{3/2}\|u\|_{H^2(\Omega)} + C_\gamma h^{1/2 + 2\gamma}\|u\|_{H^{2+\gamma}(\Omega)}, \quad \forall\ u\in H^{2 + \gamma}(\Omega).
\end{align}
\begin{proof}
By using \eqref{sobolev_extension}, \eqref{sobolev_inequality_holder} with $\gamma = \frac{3}{4}$, \eqref{estimate_2} and \eqref{estimate_5} we have that
\begin{equation}
\label{my_lemma_proof_0}
\begin{split}
&\|\tilde{u}-u^{-\ell}\|_{L^2(\Omega_h)} = \|\tilde{u} - \tilde{u}\circ G\|_{L^2(\Omega_h)} \leq C\|\tilde{u}\|_{H^{1+3/4}(\Omega_h)}\|(Id-G)^{3/4}\|_{L^2(\Omega_h)}\\
= &C\|u\|_{H^{1+3/4}(\Omega)}\|(Id-G)^{3/4}\|_{L^2(\Omega_B)} \leq C\|u\|_{H^{1+3/4}(\Omega)} |\Omega_B|^{1/2}\|Id-G\|_{L^\infty(\Omega_B)}^{3/4}\\
\leq &Ch^{1/2} h^{3/2}\|u\|_{H^{1+3/4}(\Omega)} = Ch^2\|u\|_{H^{1+3/4}(\Omega)},
\end{split}
\end{equation}
which proves \eqref{error_between_lift_and_extension_l2}. Notice that, in the last line of \eqref{my_lemma_proof_0}, the $h^{1/2}$ term is the effect of the Sobolev extension being exact except on the discrete narrow band $\Omega_B$, while the $h^{3/2}$ terms is the approximation accuracy, which is intuitively justified by the Sobolev index $1+3/4$ being larger than the exponent $3/2$. Using \eqref{narrow_band_inequality}, \eqref{sobolev_extension}, \eqref{estimate_3} and \eqref{estimate_5} we have that
\begin{align}
\label{my_lemma_proof_1}
&|\tilde{u} - u^{-\ell}|_{H^1(\Omega_h)} = \|\nabla \tilde{u} - (JG^T\ \nabla\tilde{u})\circ G\|_{L^2(\Omega_h)}\\
\notag
\leq &\|(Id-JG^T\circ G)\|_{L^\infty(\Omega_h)}\|\nabla \tilde{u}\|_{L^2(\Omega_B)} + \|JG^T\circ G\|_{L^\infty(\Omega_h)}\|\nabla\tilde{u} - \nabla\tilde{u}\circ G\|_{L^2(\Omega_h)}\\
\notag
\leq &Ch\|\nabla\tilde{u}\|_{L^2(\Omega_B)} + C\|\nabla\tilde{u} - \nabla\tilde{u}\circ G\|_{L^2(\Omega_h)} \leq Ch^{3/2}\|u\|_{H^2(\Omega)} + C\|\nabla\tilde{u} - \nabla\tilde{u}\circ G\|_{L^2(\Omega_h)}.
\end{align}
Since $\tilde{u} \in H^{2+\gamma}(\Omega_h)$, then $\nabla\tilde{u} \in H^{1+\gamma}(\Omega_h)$. Hence, by reasoning as in \eqref{my_lemma_proof_0} we have that
\begin{equation}
\label{my_lemma_proof_3}
\|\nabla\tilde{u} - \nabla\tilde{u}\circ G\|_{L^2(\Omega_h)} \leq C_\gamma h^{1/2 + 2\gamma}\|u\|_{H^{2+\gamma}(\Omega)}.
\end{equation}
By substituting \eqref{my_lemma_proof_3} into \eqref{my_lemma_proof_1} we get the desired estimate.
\end{proof}
\end{lemma}

\subsection{Virtual element space {and operators} in the bulk}
\label{sec:vem_space_bulk}
In this section, we define virtual element spaces on polygons and polygonal domains by following \cite{da2017high}. Let $E$ be a polygon in $\mathbb{R}^2$. A preliminary virtual element space on $E$ is given by
\begin{equation}
\tilde{\mathbb{V}}(E) := \Big\{v \in H^1(E) \cap \mathcal{C}^0(E) \Big|\ v_{|e} \in\mathbb{P}_1(e),\ \forall\ e\in\text{edges}(E) \land \Delta v \in \mathbb{P}_1(E) \Big\},
\end{equation}
where $\mathbb{P}_1(E)$ is the space of linear polynomials on the polygon $E$. The functions in $\mathbb{V}(E)$ are not known in closed form, but we are able to use them in a spatially discrete method, hence the name \emph{virtual}. Let us consider the elliptic projection $\Pi^\nabla_E : \tilde{\mathbb{V}}(E) \rightarrow \mathbb{P}_1(E)$ defined by
\begin{equation}
\displaystyle\int_E \nabla (v - \Pi^\nabla_E v) \cdot \nabla p_1 = 0 \quad \forall\ p_1 \in \mathbb{P}_1(E); \quad \text{and} \quad \displaystyle\int_{\partial E} (v - \Pi^\nabla_E v) = 0.
\end{equation}
Using Green's formula, it is easy to see that the operator $\Pi_E^\nabla$ is computable, see \cite{ahmad2013equivalent} for the details. The so-called \emph{enhanced virtual element space} in two dimensions is now defined as follows:
\begin{equation}
\mathbb{V}(E) := \left\{v\in \tilde{\mathbb{V}}(E) \middle|\ \int_E vp_1 = \int_E (\Pi_E^\nabla v)p_1,\ \forall\ p_1\in\mathbb{P}_1(E)\right\}.
\end{equation}
The practical usability of the space $\mathbb{V}(E)$ stems from the following result.
\begin{proposition}[Degrees of freedom]
\label{prop:dofs}
Let $n\in\mathbb{N}$. If $E$ is a polygon with $n$ vertices $\boldx_i$, $i=1,\dots,n$, then $\dim(\mathbb{V}(E)) = n$ and each function $v \in \mathbb{V}(E)$ is uniquely defined by the nodal values $v(\boldx_i)$, $i=1,\dots,n$. Hence, the nodal values constitute a set of degrees of freedom.
\begin{proof}
See \cite{ahmad2013equivalent}.
\end{proof}
\end{proposition}

For $s=1,2$ we define the broken bulk Sobolev seminorms as $|u|_{s,\Omega,h} := \sum_{E\in\mathcal{E}_h} |u_{|E}|_{H^s(E)}$. The approximation properties of the space $\mathbb{V}(E)$ are given by the following result.
\begin{proposition}[Projection error on $\mathbb{P}_1(E)$]
For all $E\in\mathcal{E}_h$, $s\in \{1,2\}$ and $w\in H^s(E)$ there exists $w_\pi \in \mathbb{P}_1(E)$ such that
\begin{equation}
\label{projection_error}
\|w-w_\pi\|_{L^2(E)} + h_E |w-w_\pi|_{H^1(E)} \leq C h_E^s |w|_{H^s(E)},
\end{equation}
where $C$ is a constant that depends only on $\gamma_1$.
\begin{proof}
See \cite{ahmad2013equivalent}.
\end{proof}
\end{proposition}

\begin{remark}[Regularity of $\mathbb{V}(E)$ functions]
\label{rmk:regularity_vem_space}
If $E$ is a convex polygon, from the properties of the Poisson problem on convex Lipschitz domains, it holds that $\mathbb{V}(E) \subset H^2(E)$, see for instance \cite{savare1998regularity}. Otherwise, if $E$ is non-convex, we may only assert that $\mathbb{V}(E) \subset H^{1+\varepsilon}(E)$ for $0\leq \varepsilon < 1/2$, see \cite{savare1998regularity}. In either case, $\mathbb{V}(E) \subset\mathcal{C}^0(E)$, see Theorem \ref{thm:sobolev_embedding}. We will account for this regularity issue in devising a numerical method with optimal convergence.
\end{remark}

The discontinuous and continuous bulk virtual element spaces are defined by pasting local spaces:
\begin{equation}
\mathbb{V}_{\Omega,h} := \{v:\Omega_h \rightarrow\mathbb{R} \ | \ v_{|E} \in \mathbb{V}(E), \ \forall\ E \in \mathcal{E}_h\}; \quad \text{and} \quad \mathbb{V}_\Omega := \mathbb{V}_{\Omega,h} \cap \mathcal{C}^0(\Omega_h).
\end{equation}
Thanks to Remark \ref{rmk:regularity_vem_space}, the only source of discontinuity in $\mathbb{V}_{\Omega,h}$ are jumps across edges. In $\mathbb{V}_\Omega$ we consider the Lagrange basis $\{\varphi_i, \ i=1,\dots,N\}$ where, for each $i=1,\dots,N$, $\varphi_i$ is the unique $\mathbb{V}_\Omega$ function such that $\varphi_i(\boldx_j) = \delta_{ij},$ for all $j = 1,\dots,N$, with $\delta_{ij}$ being the Kronecker symbol. The set $\{\varphi_i, \ i=1,\dots,N\}$ is a basis of $\mathbb{V}_\Omega$ thanks to Proposition \ref{prop:dofs}.

In the remainder of this section, let $E$ be an element of $\Omega_h$. The \emph{stabilizing form} $S_E: \mathbb{V}(E) \times \mathbb{V}(E) \rightarrow \mathbb{R}$, is defined by
\begin{align}
S_E(v,w) := \sum_{P \in \text{ vertices }(E)} v(P)w(P), \qquad \forall \ v,\, w\in \mathbb{V}(E).
\end{align}
The $L^2$ projector $\Pi^0_E: \mathbb{V}(E) \rightarrow\mathbb{P}_1(E)$ is defined as follows: for $w \in \mathbb{V}(E)$:
\begin{align}
\label{L2_projector}
\int_E (w - \Pi^0_E w)p_1 = 0, \qquad \forall\ \ p_1\in\mathbb{P}_1(E).
\end{align}
As shown in \cite{ahmad2013equivalent}, $\Pi_E^0$ is computable because $\Pi_E^0 = \Pi_E^\nabla$. Even if $\Pi_E^0$ is not a new projector, the presentation and the analysis of the method benefit from the usage of the equivalent definition \eqref{L2_projector}. Moreover, since $\Pi_E^0 = \Pi_E^\nabla$, the boundedness property of projection operators in Hilbert spaces translates to
\begin{equation}
\label{L2_projector_boundedness}
\|\Pi^0_E w\|_{L^2(E)} \leq \|w\|_{L^2(E)} \quad \text{and} \qquad |\Pi^0_E w|_{H^1(E)} \leq |w|_{H^1(E)}.
\end{equation}
We are now ready to introduce the approximate $L^2$ bilinear form $m_E: \mathbb{V}(E) \times \mathbb{V}(E) \rightarrow\mathbb{R}$ and the approximate gradient-gradient bilinear form $a_E: \mathbb{V}(E) \times  \mathbb{V}(E) \rightarrow\mathbb{R}$, defined as follows:
\begin{align}
&m_E(v,w) := \int_E (\Pi^0_E v)(\Pi^0_E w) + \area(E) S_E(v-\Pi^0_E v, w - \Pi^0_E w);\\
&a_E(v,w) := \int_E (\nabla \Pi^\nabla_E v) \cdot (\nabla \Pi^\nabla_E w) + h_E S_E(v-\Pi^\nabla_E v, w - \Pi^\nabla_E w),
\end{align}
for all $v,w\in \mathbb{V}(E)$, where $h_E$ is the diameter of $E$. The following result easily follows from the construction of the bilinear forms of $a_E$ and $m_E$.
\begin{proposition}[Stability and consistency]
\label{prop:stab_cons}
The bilinear forms $a_h$ and $m_h$ are consistent, i.e. for all $v \in \mathbb{V}(E)$ and $p \in \mathbb{P}_1(E)$
\begin{align}
\label{consistency}
a_E(v,p) = \int_E \nabla v \cdot \nabla p; \quad \text{and} \quad m_E(v,p) = \int_E vp.
\end{align}
The bilinear forms $a_h$ and $m_h$ are stable, meaning that there exists two constants $0 < \alpha_* < \alpha^*$ depending on $\gamma_2$ such that, for all $v\in\mathbb{V}(E)$
\begin{align}
\label{stability_stifness}
&\alpha_* \int_E \nabla v \cdot \nabla v \leq a_E(v,v) \leq \alpha^* \int_E \nabla v \cdot \nabla v;\\
\label{stability_mass}
&\alpha_* \int_E v^2 \leq m_E(v,v) \leq \alpha^* \int_E v^2.
\end{align}
\begin{proof}
See \cite{beirao2013basic}.
\end{proof}
\end{proposition}
We observe from \eqref{stability_stifness} and \eqref{stability_mass} that the error in the approximate bilinear forms $a_E$ and $m_E$ is not a function of the meshsize $h$, see also \cite{beirao2013basic}. Nevertheless, we will show that the method retains optimal convergence thanks to the consistency properties \eqref{consistency}.
The global bilinear forms are defined by pasting the corresponding local bilinear forms as follows:
\begin{align}
&m_h(v,w) := \sum_{E\in\mathcal{E}_h}  m_E(v_{|E}, w_{|E}); \quad \text{and} \quad a_h(v,w) := \sum_{E\in\mathcal{E}_h}  a_E(v_{|E}, w_{|E}),
\end{align}
for all $v,w\in \mathbb{V}_\Omega$. A consequence of Proposition \ref{prop:stab_cons} is that $m_h: \mathbb{V}_\Omega \times \mathbb{V}_\Omega \rightarrow\mathbb{R}$ is positive definite, while $a_h: \mathbb{V}_\Omega \times \mathbb{V}_\Omega \rightarrow\mathbb{R}$ is positive semi-definite.

{The bilinear form $m_h$} is not sufficient to discretise load terms like $\intomega f\varphi$, because $f$ is not in the space $\mathbb{V}_\Omega$. We resolve this issue by combining the approaches in \cite{adak2019convergence} and \cite{frittelli2017preserving, frittelli2018virtual}. From \cite{adak2019convergence} we take the usage of the projection operator $\Pi^\nabla$, while from \cite{frittelli2017preserving, frittelli2018virtual} we take the usage of the Lagrange interpolant. In the context of virtual elements, the Lagrange interpolant is defined as follows.
\begin{definition}[Virtual Lagrange Interpolant]
Given an element-wise continuous function $f : \Omega_h \rightarrow \mathbb{R}$, $f_{|E} \in\mathcal{C}(E)$ for all $E\in\mathcal{E}_h$, the \emph{virtual Lagrange interpolant} $I_\Omega f$ of $f$ is the unique $\mathbb{V}_{\Omega,h}$ function such that $I_\Omega f_{|E}(\boldx_i) = f(\boldx_i)$ for all $i$, with $\boldx_i\in\ \text{nodes}(E)$ and $E\in\mathcal{E}_h$.
In particular, if $f \in \mathcal{C}(\Omega_h)$, then $I_\Omega f$ is the unique $\mathbb{V}_\Omega$ function such that $I_\Omega f(\boldx_i) = f(\boldx_i)$ for all $i=1,\dots,N$.
\end{definition}
The following result provides an estimate of the interpolation error in the bulk.
\begin{proposition}[Interpolation error in the bulk]
If $w:\Omega_h \rightarrow\mathbb{R}$ is such that $w_E \in H^2(E)$ for all $E\in\mathcal{E}_h$, then the interpolant $I_\Omega(w)$ fulfils
\begin{equation}
\label{interpolation_error_bulk}
\|w-I_\Omega(w)\|_{L^2(\Omega_h)} + h|w-I_\Omega(w)|_{1,\Omega,h} \leq Ch^2|w|_{2,\Omega,h},
\end{equation}
where $C>0$ depends only on $\gamma_1$.
\begin{proof}
See \cite{ahmad2013equivalent}.
\end{proof}
\end{proposition}
Unlike projection operators, we may not assert that the interpolant $I_\Omega$ is bounded in $L^2(\Omega_h)$. However from \eqref{interpolation_error_bulk} we have the \emph{quasi-boundedness} of $I_\Omega$ in $L^2(\Omega_h)$:
\begin{equation}
\label{interpolant_bulk_quasi_bounded}
\|I_\Omega(w)\|_{L^2(\Omega_h)} \leq \|w\|_{L^2(\Omega_h)} + Ch^2|w|_{2,\Omega,h},
\end{equation}
where $C>0$ depends only on $\gamma_1$.

\subsection{Finite element space and {operators on the surface}}
\label{sec:ingredients_surface}
Let $F\in\mathcal{F}_h$ be an edge of the approximated curve $\Gamma_h$. The local finite element space on $F$ is the space $\mathbb{P}_1(F)$ of linear polynomials on $F$. The global finite element space on $\Gamma_h$ is defined by $\mathbb{V}_\Gamma := \Big\{v \in \mathcal{C}^0(\Gamma_h) \Big|\ v_{|F} \in \mathbb{P}_1(F),\ \text{for all}\ F \in \mathcal{F}_h \Big\}$,
i.e. the space of piecewise linear functions on the approximated curve $\Gamma_h$. On $\mathbb{V}_\Gamma$ we consider the Lagrange basis $\{\psi_k, \ k=1,\dots, M\}$
where, for each $k=1,\dots, M$, $\psi_k$ is the unique $\mathbb{V}_\Gamma$ function such that
$\psi_k(\boldx_l) = \delta_{kl}$ for all $l=1,\dots,M$. It is easy to see that the Lagrange basis functions of $\mathbb{V}_\Gamma$ are the restrictions to $\Gamma_h$ of the first $M$ Lagrange basis functions of $\mathbb{V}_\Omega$, that is:
\begin{equation}
\label{basis_functions_relation}
\varphi_{k | \Gamma_h} = \psi_k, \qquad \forall\ k=1,\dots,M.
\end{equation}
Before introducing the spatially discrete formulations, we are left to treat terms like $\intgamma g\varphi$, since $g$ is not in the boundary finite element space $\mathbb{V}_\Gamma$.
\begin{definition}[Surface Lagrange interpolant]
If $g: \Gamma_h \rightarrow\mathbb{R}$ is a continuous function, the \emph{Lagrange interpolant} $I_\Gamma(g)$ of $g$ is the unique $\mathbb{V}_\Gamma$ function such that $I_\Gamma(g)(\boldx_i) = g(\boldx_i)$ for all $i=1,\dots,M$.
\end{definition}
We consider the broken surface Sobolev norm $|v|_{2,\Gamma,h} := \sum_{F\in\mathcal{F}_h} |V_{|F}|_{H^2(F)}$.
The following are basic properties of Lagrange interpolation.
\begin{lemma}[Properties of Lagrange interpolation on the surface]
Let $v\in \mathcal{C}^0(\Gamma_h)$ such that, for every $F\in\mathcal{F}_h$, $v_{|F}\in H^2(F)$. Then
\begin{equation}
\label{interpolation_error_surface}
\|v-I_\Gamma(v)\|_{L^2(\Gamma_h)} + h|v-I_\Gamma(v)|_{H^1(\Gamma_h)} \leq Ch^2 |v|_{2,\Gamma,h}.
\end{equation}
For any $w\in \mathcal{C}^0(\Gamma_h)$, the surface Lagrange interpolant fulfils
\begin{equation}
\label{interpolation_surface_monotonic}
{|I_\Gamma(w)| \leq I_\Gamma(|w|).}
\end{equation}
\begin{proof}
See \cite{ciarlet2002} for \eqref{interpolation_error_surface}, while \eqref{interpolation_surface_monotonic} {is a consenquence of $I_\Gamma$ being monotonic, i.e. $I_\Gamma (v) \leq I_\Gamma (w)$ for any $v,w \in \mathcal{C}^0(\Gamma_h)$ such that $v\leq w$.}
\end{proof}
\end{lemma}

\subsection{The spatially discrete formulations}
\label{sec:spatially_discrete_formulations}
We are now ready to introduce the \emph{BSVEM discretisation} of the weak problems \eqref{elliptic_problem_weak_form} and \eqref{weak_formulation}. The discrete counterpart of the elliptic problem \eqref{elliptic_problem_weak_form} is: find $U \in \mathbb{V}_\Omega$ and $V \in \mathbb{V}_\Gamma$ such that
\begin{equation}
\label{elliptic_problem_BSVEM}
\begin{cases}
\vspace{2mm}
\displaystyle a_h(U, \varphi) + m_h(U, \varphi) + \intgammah (\alpha U -\beta V)\varphi = m_h(I_\Omega(f), \varphi);\\
\displaystyle \intgammah \nablagammah V\cdot\nablagammah \psi + \intgammah (-\alpha U + (\beta+1)V) \psi = \intgammah I_\Gamma(g)\psi,
\end{cases}
\end{equation}
for all $\varphi \in \mathbb{V}_\Omega$ and $\psi \in \mathbb{V}_{\Gamma}$. We express the spatially discrete solution $(U,V)$ in the Lagrange bases as follows:
\begin{align}
U(\boldx) = \sum_{i=1}^N \xi_i \varphi_i(\boldx), \quad \boldx \in \Omega_h;\quad \text{and} \quad 
V(\boldx) = \sum_{k=1}^M \eta_k \psi_k(\boldx), \quad \boldx \in \Gamma_h.
\end{align}
Hence, problem \eqref{elliptic_problem_BSVEM} is equivalent to: find $\boldxi := (\xi_i, \dots, \xi_N)^{T} \in \mathbb{R}^{N}$ and $\boldeta := (\eta_1,\dots, \eta_M)^{T} \in\mathbb{R}^{M}$ such that 
\begin{equation}
\label{elliptic_problem_BSVEM_basis}
\begin{cases}
\vspace{2mm}
\begin{aligned}
\displaystyle\sum_{i=1}^N\displaystyle \xi_i a_h(\varphi_i, \varphi_j) + \sum_{i=1}^N \xi_i m_h(\varphi_i, \varphi_j) + \alpha & \sum_{k=1}^M \xi_k \intgammah \varphi_k\varphi_l\\
 - \beta  \sum_{k=1}^M \eta_k \intgammah \psi_k\varphi_l = &\sum_{i=1}^N f(\boldx_i) m_h(\varphi_i, \varphi_j);
\end{aligned}\\
\begin{aligned}
\displaystyle \sum_{k=1}^M \eta_k \intgammah \nablagammah \psi_k \cdot\nablagammah \psi_l -\alpha & \sum_{k=1} ^ M \xi_k \intgammah \varphi_k\psi_l\\
+ (\beta+1)  \sum_{k=1} ^ M \eta_k \intgammah \psi_k \psi_l = & \sum_{k=1}^M g(\boldx_k)\intgammah \psi_k\psi_l,
\end{aligned}
\end{cases}
\end{equation}
for all $j=1,\dots,N$ and $l=1,\dots,M$. We consider the matrices $A_\Omega = (a_{i,j}^{\Omega}) \in\mathbb{R}^{N\times N}$,  $M_\Omega = (m_{i,j}^{\Omega}) \in\mathbb{R}^{N \times N}$, $A_\Gamma = (a_{k,l}^{\Gamma}) \in \mathbb{R}^{M \times M}$, $M_{\Gamma} = (m_{k,l}^{\Gamma}) \in\mathbb{R}^{M \times M}$ defined as follows:
\begin{align}
\label{stiffness_matrix_bulk}
&a_{i,j}^{\Omega} := a_h(\varphi_i, \varphi_j); \quad \text{and}\quad m_{i,j}^{\Omega} := m_h(\varphi_i, \varphi_j), \qquad i,j=1,\dots,N;\\
\label{stiffness_matrix_surface}
&a_{k,l}^{\Gamma} := \intgammah \nablagammah \psi_k\cdot \nablagammah \psi_l; \quad \text{and} \quad m_{k,l}^{\Gamma} := \intgammah \psi_k \psi_l, \qquad k,l=1,\dots,M.
\end{align}
Moreover, we consider the column vectors $\boldf \in\mathbb{R}^{N}$ and $\boldg\in \mathbb{R}^{M}$ defined by
\begin{equation}
\boldf :=
\left[\begin{matrix}
f(\boldx_1)\\
\vdots\\
f(\boldx_N)
\end{matrix}\right];
\qquad \text{and}
\qquad
\boldg :=
\left[\begin{matrix}
g(\boldx_1)\\
\vdots\\
g(\boldx_M)
\end{matrix}\right].
\end{equation}
By using \eqref{basis_functions_relation} we can now rewrite the discrete formulation \eqref{elliptic_problem_BSVEM_basis} in matrix-vector form as a block $(N+M)\times(N+M)$ linear algebraic system:
\begin{equation}
\label{elliptic_problem_BSVEM_linear_system}
\begin{cases}
A_\Omega \boldxi + M_\Omega \boldxi + \alpha RM_\Gamma R^T\boldxi - \beta R M_\Gamma \boldeta = M_\Omega \boldf;\\
A_\Gamma \boldeta - \alpha M_\Gamma R^T\boldxi + (\beta + 1) M_\Gamma \boldeta = M_\Gamma \boldg.
\end{cases}
\end{equation}
In compact form, \eqref{elliptic_problem_BSVEM_linear_system} reads
\begin{equation}
 \label{elliptic_problem_BSVEM_linear_system_compact}
\left[\begin{matrix}
A_\Omega + M_\Omega + \alpha RM_\Gamma R^T     & - \beta R M_\Gamma\\
- \alpha M_\Gamma R^T                                            & A_\Gamma + (\beta + 1) M_\Gamma
\end{matrix}\right]
\left[\begin{matrix}
\boldxi\\
\boldeta
\end{matrix}\right]
=
\left[\begin{matrix}
M_\Omega\boldf\\
M_\Gamma\boldg
\end{matrix}\right],
\end{equation}
\noindent
which is uniquely solvable, {as a consequence of the positive semi-definiteness of $a_h$ and the positive definiteness of $m_h$.} 

The spatial discretisation of the parabolic problem \eqref{weak_formulation} is: find $U\in L^2([0,T]; \mathbb{V}_\Omega)$ and $V \in L^2([0,T]; \mathbb{V}_\Gamma)$ such that
\begin{equation}
\label{parabolic_problem_BSVEM_projection}
\hspace*{-2mm}
\begin{cases}
\vspace{2mm}
\displaystyle m_h \left(\Ut, \varphi\right) +  d_u a_h(U, \varphi) = m_h(\Pi^0 I_\Omega(q(\Pi^0 U)), \varphi) + \intgammah I_\Gamma(s(U,V))\varphi;\\
\displaystyle\intgammah \Vt \psi +  d_v \intgamma \nablagammah V \cdot \nablagammah \psi + \intgammah I_\Gamma(s(U,V))\psi = \intgammah I_\Gamma(r(U,V))\psi,
\end{cases}
\end{equation} 
for all $U \in L^2([0,T]; \mathbb{V}_{\Omega})$ and $V \in L^2([0,T]; \mathbb{V}_\Gamma)$. The discrete initial conditions are prescribed as follows
\begin{equation}
\label{parabolic_problem_BSVEM_IC}
U_0 = I_\Omega(u_0); \quad \text{and} \quad V_0 = I_\Gamma(v_0).
\end{equation}

\begin{remark}[Special cases]
If every element $E\in\mathcal{E}_h$ is convex or $f$ is linear, optimal convergence is retained by replacing the term $m_h(\Pi^0 I_\Omega(q(\Pi^0 U)), \varphi)$ with $m_h(I_\Omega(q(U)), \varphi)$, i.e by removing the projection operator $\Pi^0$. By expressing the time-dependent semi-discrete solution $(U,V)$ in the Lagrange bases as follows
\begin{align}
U(\boldx, t) = \sum_{i=1}^N \xi_i(t) \varphi_i(\boldx), \quad \boldx \in \Omega_h; \quad
 \text{and} \quad V(\boldx, t) = \sum_{k=1}^M \eta_k(t) \psi_k(\boldx), \quad \boldx \in \Gamma_h,
\end{align}
the fully discrete problem can be written in matrix-vector form as an $(M+N)\times (M+N)$ nonlinear ODE system of the form:
\begin{equation}
\label{parabolic_problem_BSVEM_ODE}
\begin{cases}
&M_\Omega \dot{\boldxi}(t) + {d_u}K_\Omega \boldxi(t)= M_\Omega q(\boldxi(t)) + R M_\Gamma h(R^{T}\boldxi(t), \boldeta(t));\\
&M_\Gamma \dot{\boldeta}(t) + {d_v}K_\Gamma \boldeta(t) = - M_\Gamma h(R^{T}\boldxi(t), \boldeta(t)) + M_\Gamma r(R^{T}\boldxi(t), \boldeta(t)).
\end{cases}
\end{equation}
\end{remark}

\begin{remark}[Implementation]
Thanks to the reduction matrix $R$, we are able to implement the spatially discrete problems \eqref{elliptic_problem_BSVEM_linear_system} and \eqref{parabolic_problem_BSVEM_ODE} by using only two kinds of mass matrix ($M_\Omega$ and $M_\Gamma$), two kinds of stiffness matrix ($A_\Omega$ and $A_\Gamma$) and $R$ itself. In previous works on BSFEM as illustrated in \cite{Madzvamuse_2016}, five kinds of mass matrix were used to evaluate the integrals appearing in the spatially discrete formulation. We stress once again that, since the pre-existing BSFEM is a special case of the proposed BSVEM, this work provides, as a by-product, an optimised matrix implementation of the BSFEM.
\end{remark}

\section{Stability and convergence analysis}
\label{sec:convergence_analysis}

The spatially discrete parabolic problem \eqref{parabolic_problem_BSVEM_projection} fulfils the following stability estimates.

\begin{lemma}[Stability estimates for the spatially discrete parabolic problem \eqref{parabolic_problem_BSVEM_projection}]
There exists $C>0$ depending on $q$, $r$, $s$ and $\Omega$ such that
\begin{align}
\label{stability_parabolic_semidiscrete_1}
\sup_{t\in[0,T]}\Big(\|U\|_{L^2(\Omega_h)}^2 + \|V\|_{L^2(\Gamma_h)}^2\Big) + &\int_0^T \Big(|U|_{H^1(\Omega_h)}^2 + |V|_{H^1(\Gamma_h)}^2\Big)\\
\notag
\leq &C\Big(1+ \|U_0\|_{L^2(\Omega_h)}^2 + \|V_0\|_{L^2(\Gamma_h)}^2\Big)\exp(CT);\\
\label{stability_parabolic_semidiscrete_2}
\sup_{t\in[0,T]}\Big(|U|_{H^1(\Omega_h)}^2 + |V|_{H^1(\Gamma_h)}^2\Big) + &\int_0^T \Big(\|\dot{U}\|_{L^2(\Omega_h)}^2 + \|\dot{V}\|_{L^2(\Gamma_h)}^2\Big)\\
\notag
\leq &C\Big(1 + \|U_0\|_{H^1(\Omega_h)}^2 + \|V_0\|_{H^1(\Gamma_h)}^2\Big)\exp(CT).
\end{align}
\begin{proof}
The proof relies on standard energy techniques. By choosing $\varphi = U$ and $\psi = V$ in \eqref{parabolic_problem_BSVEM_projection}, using the Lipschitz continuity of $q$, $r$, $s$, the Young's inequality and summing over the equations we have
\begin{equation}
\label{stability_proof_1}
\begin{split}
\frac{1}{2}\frac{\mathrm{d}}{\mathrm{d}t} \left(m_h(U, U) + \|V\|_{L^2(\Gamma_h)}^2\right) + d_u a_h(U,U) + d_v |V|_{H^1(\Gamma_h)}^2\\
\leq C\left(1 + \|U\|_{L^2(\Omega_h)}^2 + \|V\|_{L^2(\Gamma_h)}^2\right) + c \|U_{|\Gamma_h}\|_{L^2(\Gamma_h)}^2,
\end{split}
\end{equation}
where $c>0$ is arbitrarily small, thanks to the Young's inequality. By applying \eqref{equivalence_omega_h1}, \eqref{trace_inequality} and \eqref{stability_stifness} to the last term in \eqref{stability_proof_1} we can choose $c$ such that we have
\begin{equation}
\label{stability_proof_1_5}
\begin{split}
\frac{1}{2}\frac{\mathrm{d}}{\mathrm{d}t} \left(m_h(U, U) + \|V\|_{L^2(\Gamma_h)}^2\right) + d_u a_h(U,U) + d_v |V|_{H^1(\Gamma_h)}^2\\
\leq C\left(1 + \|U\|_{L^2(\Omega_h)}^2 + \|V\|_{L^2(\Gamma_h)}^2\right) +\frac{d_u}{2} a_h(U,U).
\end{split}
\end{equation}
By using \eqref{stability_mass} into \eqref{stability_proof_1_5} we have
\begin{equation}
\label{stability_proof_2}
\begin{split}
\frac{1}{2}\frac{\mathrm{d}}{\mathrm{d}t} \left(m_h(U, U) + \|V\|_{L^2(\Gamma_h)}^2\right) \leq &C\left(m_h(U,U) + \|V\|_{L^2(\Omega_h)}^2\right)\\
+ &C - \frac{d_u}{2} a_h(U,U) - d_v|V|_{H^1(\Gamma_h)}^2.
\end{split}
\end{equation}
An application of Gr\"{o}nwall's lemma to \eqref{stability_proof_2} yields 
\begin{equation}
\begin{split}
\sup_{t\in[0,T]}\Big(m_h(U,U) + &\|V\|_{L^2(\Gamma_h)}^2\Big) \leq \Big(1+ m_h(U_0, U_0) + \|V_0\|_{L^2(\Gamma_h)}^2\Big)\exp(CT) - 1\\
- &\int_0^T \Big(\frac{d_u}{2} a_h(U,U) + d_v |V|_{H^1(\Gamma_h)}^2\Big)\exp\{C(T-t)\}\mathrm{d}t,
\end{split}
\end{equation}
which yields \eqref{stability_parabolic_semidiscrete_1} after an application of \eqref{equivalence_omega_h1}. Similarly, by choosing $\varphi = \dot{U}$ and $\psi = \dot{V}$ in \eqref{parabolic_problem_BSVEM_projection} and summing over the equations we have
\begin{equation}
\label{stability_proof_3}
\begin{split}
&m_h(\dot{U}, \dot{U}) + \|\dot{V}\|_{L^2(\Gamma_h)}^2 + \frac{1}{2}\frac{\mathrm{d}}{\mathrm{d}t}\left(d_u a_h(U,U) + d_v |V|_{H^1(\Gamma_h)}^2\right)\\
\leq &C\left(1 + \|U\|_{L^2(\Omega_h)}^2 + \|V\|_{L^2(\Gamma_h)}^2\right) + C a_h(U,U) + \frac{1}{2}\Big(m_h(\dot{U}, \dot{U}) + \|\dot{V}\|_{L^2(\Gamma_h)}^2\Big),
\end{split}
\end{equation}
where we have exploited the Lipschitz continuity of $q$, $r$, $s$ and  the Young's inequality and \eqref{stability_mass}. From \eqref{stability_proof_3} we immediately get
\begin{equation}
\label{stability_proof_4}
\begin{split}
&\frac{1}{2}\Big(m_h(\dot{U}, \dot{U}) + \|\dot{V}\|_{L^2(\Gamma_h)}^2\Big) + \frac{1}{2}\frac{\mathrm{d}}{\mathrm{d}t}\Big(d_u a_h(U,U) + d_v |V|_{H^1(\Gamma_h)}^2\Big)\\
\leq &C\Big(a_h(U,U) + |V|_{H^1(\Gamma_h)}^2\Big) + C\left(1 + \|U\|_{L^2(\Omega_h)}^2 + \|V\|_{L^2(\Gamma_h)}^2\right).
\end{split}
\end{equation}
By applying Gr\"{o}nwall's lemma and then using \eqref{stability_stifness}-\eqref{stability_mass} we obtain \eqref{stability_parabolic_semidiscrete_2}.
\end{proof}
\end{lemma}

To derive error estimates for the spatially discrete solution we need suitable bulk and surface Ritz projections. The surface Ritz projection is taken from \cite{Elliott_ranner_cahn_hilliard}, while the bulk Ritz projection is tailor-made.
\begin{definition}[Surface Ritz projection]
The \emph{surface Ritz projection} of a function $v\in H^1(\Gamma)$ is the unique function $\mathcal{R}v \in \mathbb{V}_\Gamma$ such that
\begin{align}
\label{surface_ritz}
\intgammah \nablagammah\mathcal{R}v\cdot\nablagammah\psi = \intgamma \nablagamma v\cdot\nablagamma\psi \quad \forall \psi\in \mathbb{V}_\Gamma; \quad \text{and} \quad \intgammah \mathcal{R}v = \intgamma v.
\end{align}
\end{definition}

\begin{theorem}[Error bounds for the surface Ritz projection]
The surface Ritz projection fulfils the optimal a priori error bound
\begin{equation}
\label{surface_ritz_estimate_L2_and_H1}
\|v - (\mathcal{R}v)^\ell\|_{L^2(\Gamma)} + h\|v - (\mathcal{R}v)^\ell\|_{H^1(\Gamma)} \leq Ch^2\|v\|_{H^2(\Gamma)},
\end{equation}
where $C>0$ depends only on $\Gamma$ and $\ell$ is the lifting.
\begin{proof}
See \cite{Elliott_ranner_cahn_hilliard}.
\end{proof}
\end{theorem}

We now define a suitable Ritz projection in the bulk.
\begin{definition}[Bulk Ritz projection]
The \emph{bulk Ritz projection} of a function $u\in H^1(\Omega)$ is the unique function $\mathcal{R}u \in \mathbb{V}_\Omega$  such that
\begin{align}
\label{ritz_projection_definition}
a_h(\mathcal{R}u, \psi) = \intomega \nabla u\cdot\nabla\psi^\ell \quad \forall\ \psi\in \mathbb{V}_\Omega; \quad \text{and} \quad \mathcal{R}u_{|\Gamma_h} = I_\Gamma (u^{-\ell}),
\end{align}
where $-\ell$ is the inverse lifting.
\end{definition}

In the following theorems we show that the bulk Ritz projection of a sufficiently regular function fulfils optimal a priori error bounds in $H^1(\Omega)$, $H^1(\Gamma)$, $L^2(\Gamma)$ and $L^2(\Omega)$ norms.
\begin{theorem}[$H^1(\Omega)$ a priori  error bound for the bulk Ritz projection]
\label{thm:ritz_proj_H1}
For any $u\in H^{2+1/4}(\Omega)$ it holds that
\begin{equation}
\label{bulk_ritz_estimate_bulk_H1}
|u - (\mathcal{R}u)^\ell|_{H^1(\Omega)} \leq Ch\|u\|_{H^2(\Omega)} + Ch\|u\|_{H^{2+1/4}(\Omega)},
\end{equation}
where $\ell$ is the lifting {and $C$ depends on $\Omega$ and the constants $\gamma_1$ and $\gamma_2$ in \ref{A1}-\ref{A2}.}
\begin{proof}
We set $e_h := \mathcal{R}u - \tilde{u}$. From \eqref{lifting_error_bulk_nabla}, \eqref{error_between_lift_and_extension}, \eqref{projection_error}, \eqref{consistency}, \eqref{stability_stifness} and \eqref{ritz_projection_definition} we have
\begin{align*}
&\alpha_*|e_h|_{H^1(\Omega_h)}^2 \leq a_h(e_h, e_h) = a_h(\mathcal{R}u, e_h) - a_h(\tilde{u}, e_h)\\
= &\intomega \nabla u\cdot\nabla e_h^{\ell} - \sum_{E\in\mathcal{E}_h} (a_E(\tilde{u}, e_h)) = \intomega \nabla u\cdot\nabla e_h^{\ell} - \sum_{E\in\mathcal{E}_h} (a_E(\tilde{u} - \tilde{u}_\pi, e_h) + a_E(\tilde{u}_\pi,e_h))\\
= &\intomega \nabla u\cdot\nabla e_h^{\ell} - \intomegah \nabla \tilde{u}\cdot\nabla e_h - \sum_{E\in\mathcal{E}_h} a_E(\tilde{u} - \tilde{u}_\pi, e_h)\\
= &\intomega \nabla u \cdot \nabla e_h^{\ell} - \intomegah \nabla u^{-\ell}\cdot \nabla e_h + \intomegah \nabla u^{-\ell}\cdot \nabla e_h  - \intomegah \nabla \tilde{u}\cdot\nabla e_h - \sum_{E\in\mathcal{E}_h} a_E(\tilde{u} - \tilde{u}_\pi, e_h)\\
\leq & C\left(h\|u\|_{H^2(\Omega)} + h^{3/2}\|u\|_{H^2(\Omega)} + Ch\|u\|_{H^{2+1/4}(\Omega)}\right) |e_h|_{H^1(\Omega_h)},
\end{align*}
which yields, for $h \leq h_0$
\begin{equation}
\label{ritz_projection_energy_1}
|e_h|_{H^1(\Omega_h)} \leq \left(Ch\|u\|_{H^2(\Omega)} + Ch\|u\|_{H^{2+1/4}(\Omega)}\right).
\end{equation}
By using \eqref{equivalence_omega_h1}, \eqref{error_between_lift_and_extension} and \eqref{ritz_projection_energy_1} we get
\begin{equation*}
\begin{split}
|u - (\mathcal{R}u)^\ell|_{H^1(\Omega)} \leq C|u^{-\ell} - \mathcal{R}u|_{H^1(\Omega_h)} \leq &C(|u^{-\ell} - \tilde{u}|_{H^1(\Omega_h)} + |e_h|_{H^1(\Omega_h)})\\
\leq &Ch\|u\|_{H^2(\Omega)} + Ch\|u\|_{H^{2+1/4}(\Omega)}.
\end{split}
\end{equation*}
\end{proof}
\end{theorem}

In order to prove the $L^2$ convergence, care must be taken about inverse trace operators. A consequence of Theorem \ref{thm:trace_theorem} in the Appendix A is that, given $v\in H^1(\Gamma)$, there exists $v_B\in H^1(\Omega)$ such that $\Tr(v_B) = v$ and $\|v_B\|_{H^1(\Omega)} \leq C\|v\|_{H^1(\Gamma)}$. However, for our purposes, we need the existence of a constant $C>0$ such that the bounds  $\|v_B\|_{L^2(\Omega)} \leq C\|v\|_{L^2(\Gamma)}$ and $|v_B|_{H^1(\Omega)} \leq C\|v\|_{H^1(\Gamma)}$ are simultaneously fulfilled, namely we need a $L^2$-\emph{preserving} inverse trace operator. In the following result, we prove the existence of such an operator under special assumptions on the regularity of $\Gamma$. 

\begin{lemma}[$L^2$-preserving inverse trace]
\label{lmm:L2_preserving_inverse_trace}
Assume the boundary $\Gamma$ is $\mathcal{C}^3$. Then, for any $v\in H^1(\Gamma)$ such that $\|v\|_{L^2(\Gamma)}$ is sufficiently small, there exists $v_B\in H^1(\Omega)$ such that
\begin{equation}
\label{L2_preserving_inverse_trace}
\qquad \|v_B\|_{L^2(\Omega)} \leq C\|v\|_{L^2(\Gamma)}; \quad \text{and} \quad |v_B|_{H^1(\Omega)} \leq C\|v\|_{H^1(\Gamma)}.
\end{equation}
\begin{proof}
With $\delta$ as defined in Theorem \ref{thm:co_area_formula}, we take $0 \leq \delta_0 \leq \delta$ and we define
\begin{equation}
\label{smart_inverse_trace}
v_B(\boldx) :=
\begin{cases}
v(\bolda(\boldx))\left(1+\frac{d(\boldx)}{\delta_0}\right), \qquad \text{if}\ \boldx\in U_{\delta_0};\\
0, \qquad \text{if}\ \boldx\in \Omega\setminus U_{\delta_0}. 
\end{cases}
\end{equation}
By using \eqref{co_area_formula} we have
\begin{equation}
\label{inverse_trace_1}
\begin{split}
&\intomega v_B^2(\boldx)\mathrm{d}x = \int_{U_h} v_B^2(\boldx)\mathrm{d}x = \int_{-\delta_0}^0\mathrm{d}s \int_{\Gamma_s} v^2(\bolda(\boldx))\left(1+\frac{d(\boldx)}{\delta_0}\right)^2\mathrm{d}x \\
= &\int_{-\delta_0}^0\left(1+\frac{s}{\delta_0}\right)^2\mathrm{d}s \int_{\Gamma_s} v^2(\bolda(\boldx))\mathrm{d}x = \frac{\delta_0}{3}\int_{\Gamma_s} v^2(\bolda(\boldx))\mathrm{d}x.
\end{split}
\end{equation}
Since the decomposition $(d(\boldx), \bolda(\boldx))$ is unique (see Lemma \ref{lmm:fermi}) and all the points $\boldx\in \Gamma_s$ share the same distance $d(\boldx)$, then the mapping $\bolda_s := \bolda_{|\Gamma_s} : \Gamma_s \rightarrow \Gamma$ is invertible. Moreover, since $\bolda_0 = Id_{|\Gamma}$ (which implies $\|\nabla_\Gamma \bolda\| = 1$) and $\bolda \in\mathcal{C}^2(U)$ (see Remark \ref{rmk:regularity_normal_projection}), we can choose $\delta_0$ small enough such that $0 < c \leq \|\nabla_{\Gamma_s} \bolda_s\| \leq C$. Hence, $\bolda_s^{-1}$ is $\mathcal{C}^2$ as well, which implies that $\|\nabla_{\Gamma} \bolda_s^{-1}\| \leq C$. Hence by setting $\boldy = \bolda(\boldx)$ we have
\begin{equation}
\label{inverse_trace_2}
\int_{\Gamma_s} v^2(\bolda(\boldx))\mathrm{d}x = \intgamma v^2(\boldy)\|\nabla_{\Gamma} \bolda_s^{-1}(\boldy)\|\mathrm{d}\boldy \leq C\|v\|_{L^2(\Gamma)}^2,
\end{equation}
where $\|\cdot\|$ is the Euclidean norm. By combining \eqref{inverse_trace_1} and \eqref{inverse_trace_2} we obtain the first inequality in \eqref{L2_preserving_inverse_trace}, since $\delta_0$ depends only on $\Gamma$. An application of the chain rule and Leibniz's rule yields
\begin{equation*}
\nabla v_B(\boldx) = \nablagamma v(\bolda(\boldx))J\bolda(\boldx)\left(1+\frac{d(\boldx)}{\delta_0}\right) + v(\bolda(\boldx))\frac{\nabla d(\boldx)}{\delta_0}.
\end{equation*}
Since $d(\boldx)$ and $\bolda(\boldx)$ are both $\mathcal{C}^2$ on the compact set $U_\delta \supset U_{\delta_0}$ and $-\delta_0 \leq d(\boldx) \leq 0$, we obtain
\begin{equation}
\label{inverse_trace_3}
\|\nabla v_B(\boldx)\| \leq C\|\nablagamma v(\bolda(\boldx))\| + \frac{C}{\delta_0}|v(\bolda(\boldx))|, \qquad  \forall\ \boldx\in U_{\delta_0},
\end{equation}
where $|\cdot|$ is the absolute value. Thanks to the continuous pasting in \eqref{smart_inverse_trace}, it is easy to show that the bound \eqref{inverse_trace_3} for the distributional gradient $\nabla v_B$ still holds on the junction $\Gamma_{\delta_0}$ between $U_{\delta_0}$ and $\Omega\setminus U_{\delta_0}$. From \eqref{co_area_formula}, \eqref{inverse_trace_3} and the Young's inequality we have
\begin{equation}
\begin{split}
&\intomega \|\nabla v_B(\boldx)\|^2\mathrm{d}\boldx \leq C \int_{-\delta_0}^0\mathrm{d}s\int_{\Gamma_s} \left(\|\nablagamma v(\bolda(\boldx))\|^2 + \frac{|v(\bolda(\boldx))|^2}{\delta_0^2}\right)\mathrm{d}\boldx\\
= &\int_{-\delta_0}^0\mathrm{d}s\int_{\Gamma} \left(\|\nablagamma v(\boldy)\|^2 + \frac{|v(\boldy)|^2}{\delta_0^2}\right)\|\nabla_{\Gamma} \bolda_s^{-1}(\boldy)\|\mathrm{d}\boldy\\
\leq &C \delta_0 \int_{\Gamma} \left(\|\nablagamma v(\boldy)\|^2 + \frac{|v(\boldy)|^2}{\delta_0^2}\right)\mathrm{d}\boldy \leq C\left(\|v\|_{L^2(\Gamma)}^2 + |v|_{H^1(\Gamma)}^2\right),
\end{split}
\end{equation}
which proves the second inequality in \eqref{L2_preserving_inverse_trace}.
\end{proof}
\end{lemma}

\begin{theorem}[$H^1(\Gamma)$, $L^2(\Gamma)$ and $L^2(\Omega)$ error bounds for the bulk Ritz projection]
\label{thm:ritz_proj_L2}
Let $\Omega$ have a $\mathcal{C}^3$ boundary. Then, for any $u\in H^{2+1/4}(\Omega)$ such that $\Tr(u)\in H^2(\Gamma)$ and for $h$ sufficiently small, it holds that
\begin{align}
\label{bulk_ritz_estimate_bulk_L2}
&\|u - (\mathcal{R}u)^{\ell}\|_{L^2(\Omega)} \leq Ch^{2}\left(\|u\|_{H^2(\Omega)} + \|u\|_{H^{2+1/4}(\Omega)}\right),\\
\label{bulk_ritz_estimate_boundary_L2_and_H1}
&\|u - (\mathcal{R}u)^{\ell}\|_{L^2(\Gamma)} + h|u - (\mathcal{R}u)^{\ell}|_{H^1(\Gamma)} \leq Ch^2\|u\|_{H^2(\Gamma)},
\end{align}
with $C$ depending on $\Omega$, $\gamma_1$ and $\gamma_2$. In \eqref{bulk_ritz_estimate_bulk_L2}, the term in $H^{2+1/4}(\Omega)$ norm arises only in the simultaneous presence of curved boundaries and non-triangular boundary elements.
\begin{proof}
See Appendix B.
\end{proof}
\end{theorem}

In the following theorem we prove optimal convergence in $L^\infty([0,T], H^1(\Omega)\times H^1(\Gamma))$ norm for the spatially discrete parabolic problem  \eqref{parabolic_problem_BSVEM_projection}, by harnessing the techniques used in \cite{adak2019convergence}, \cite{dziuk2013finite} and \cite{frittelli2017preserving}.
\begin{theorem}[Convergence of the BSVEM for the parabolic case]
\label{thm:convergence}
Assume that the kinetics $q$, $r$, $s$ are $\mathcal{C}^2$ and globally Lipschitz continuous (or at least Lipschitz on the range of the discrete solution). Assume that the exact solution $(u,v)$ of the parabolic problem \eqref{parabolic_problem} is such that $u$, $u_t\in L^\infty([0,T]; H^{2+1/4}(\Omega))$ and $v$, $v_t$, $\Tr(u)$, $\Tr(u_t)\in L^\infty([0,T];H^2(\Gamma))$. Let $(U,V)$ be the solution of \eqref{parabolic_problem_BSVEM_projection}. Then it holds that
\begin{equation}
\|u - U^\ell\|_{L^\infty([0,T]; L^2(\Omega))} + \|v-V^\ell\|_{L^\infty([0,T]; L^2(\Gamma))} \leq Ch^2,
\end{equation}
where the constant $C$ depends on the diffusion coefficient $d_u$, the final time $T$ and on the following norms:
\begin{itemize}
\item  $\|(u, u_t)\|_{L^\infty([0,T]; H^{2}(\Omega))}$, $\|(v,v_t,\Tr(u),\Tr(u_t))\|_{L^\infty([0,T]; H^2(\Gamma))}$ in any case;
\item $\|(u, u_t)\|_{L^\infty([0,T]; H^{2+1/4}(\Omega))}$ only in the simultaneous presence of curved boundaries and non-triangular boundary elements.
\end{itemize}
\begin{proof}
See Appendix B.
\end{proof}
\end{theorem}

\section{Time discretisation of the parabolic problem}
\label{sec:time_disc}

For the time discretisation of the semi-discrete formulation \eqref{parabolic_problem_BSVEM_ODE} of the parabolic problem \eqref{parabolic_problem} we use the IMEX (IMplicit-EXplicit) Euler method, which approximates diffusion terms implicitly and reaction- and boundary terms explicitly, see for instance \cite{frittelli2017preserving}. This choice is meant to make the implementation as simple as possible. However, if the boundary conditions and/or the reaction terms are linear, the aforementioned terms can be approximated implicitly without involving any additional nonlinear solver. We choose a timestep $\tau > 0$ and we consider the equally spaced discrete times $t_n := n\tau$, for $n=0,\dots, N_T$, with $N_T := \left\lceil \frac{T}{\tau}\right\rceil$. For $n=0,\dots, N_T$ we denote by $\boldxi^n$ and $\boldeta^n$ the numerical solution at time $t_n$. The IMEX Euler time discretisation of \eqref{parabolic_problem_BSVEM_ODE} reads
\begin{equation}
\begin{cases}
\vspace{2mm}
M_\Omega \dfrac{\boldxi^{n+1} - \boldxi^n}{\tau} + {d_u}K_\Omega \boldxi^{n+1} = M_\Omega q(\boldxi^n) + RM_\Gamma s(R^T\boldxi^n, \boldeta^n);\\
M_\Gamma \dfrac{\boldeta^{n+1} - \boldeta^n}{\tau} + {d_v}K_\Gamma \boldeta^{n+1} = M_\Gamma (-s(R^T\boldxi^n, \boldeta^n) + r(R^T\boldxi^n, \boldeta^n)),
\end{cases}
\end{equation}
for $n=0,\dots, N_T-1$, where $R$ is the reduction matrix introduced in \eqref{reduction_matrix}. By solving for $\boldxi^{n+1}$ and $\boldeta^{n+1}$ we obtain the following time stepping scheme:
\begin{equation}
\begin{cases}
(M_\Omega + \tau {d_u}K_\Omega) \boldxi^{n+1} =  M_{\Omega}(\boldxi^n + \tau q(\boldxi^n)) + \tau RM_\Gamma s(R^T\boldxi^n, \boldeta^n);\\
(M_\Gamma + \tau {d_v}K_\Gamma)\boldeta^{n+1} =  M_\Gamma (\boldeta^n + \tau  r(R^T\boldxi^n, \boldeta^n) - \tau s(R^T\boldxi^n, \boldeta^n)),
\end{cases}
\end{equation}
for $n=0,\dots, N_T-1$. Hence, two linear systems must be solved at each time step, where the matrix coefficients are the same for all $n$.

\section{Construction of meshes optimised for matrix assembly}
\label{sec:mesh_advantage}
One of the benefits of using suitable polygonal elements is the reduction in computational complexity of the matrix assembly. This is useful especially for (i) time-independent problems and (ii) time-dependent problems on evolving domains, where matrix assembly might take the vast majority of the overall computational time.

\begin{figure}[ht!]
\begin{center}
\subfigure[Step 1. The bulk $\Omega$, enclosed by the red boundary $\Gamma$, is bounded by the green square $Q$, which is subdivided with a Cartesian grid.]{\label{fig:mesh_advantage_1}\includegraphics[scale=.25]{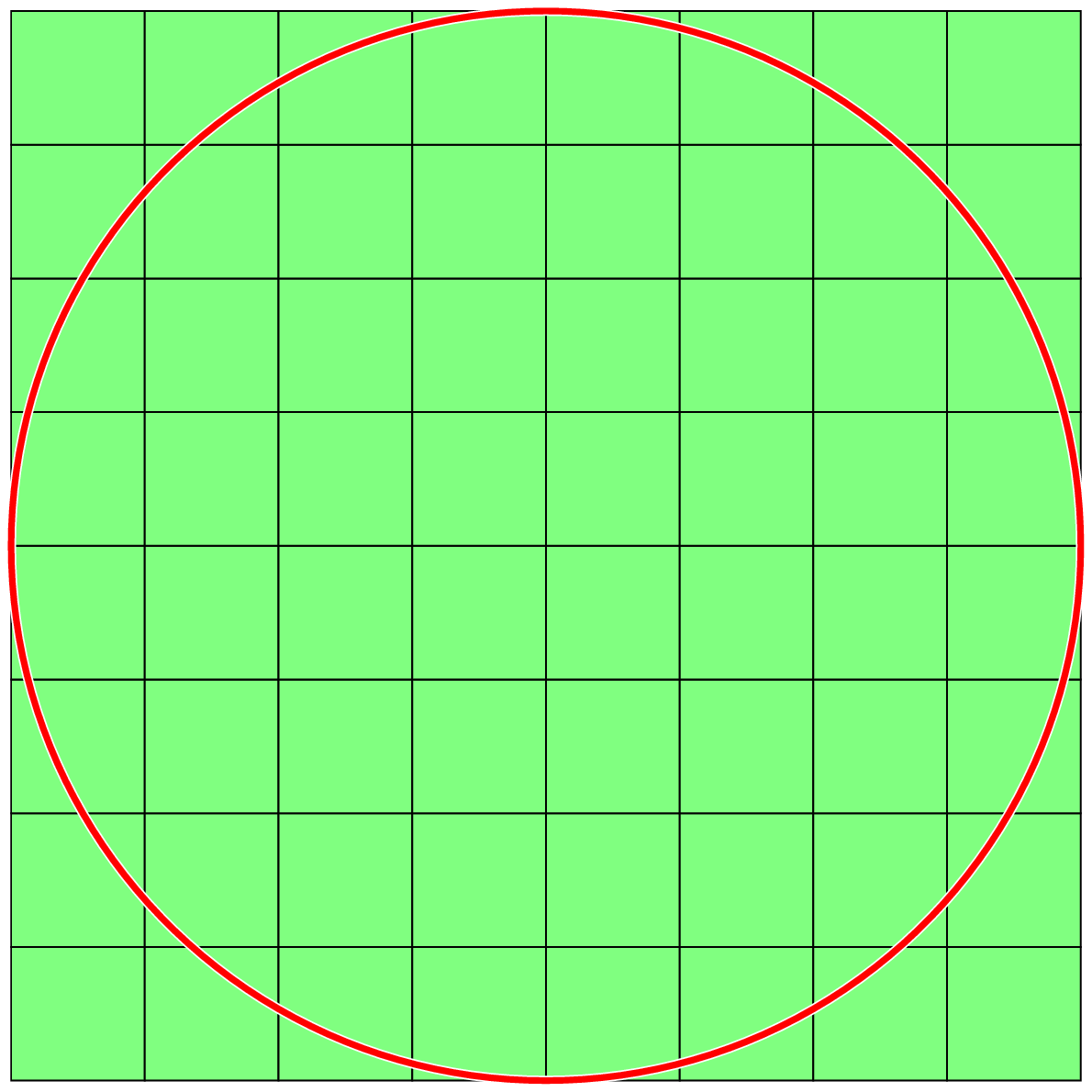}}
\subfigure[Step 2. The mesh elements that are entirely outside $\Gamma$ are discarded.]{\label{fig:mesh_advantage_2}\includegraphics[scale=.25]{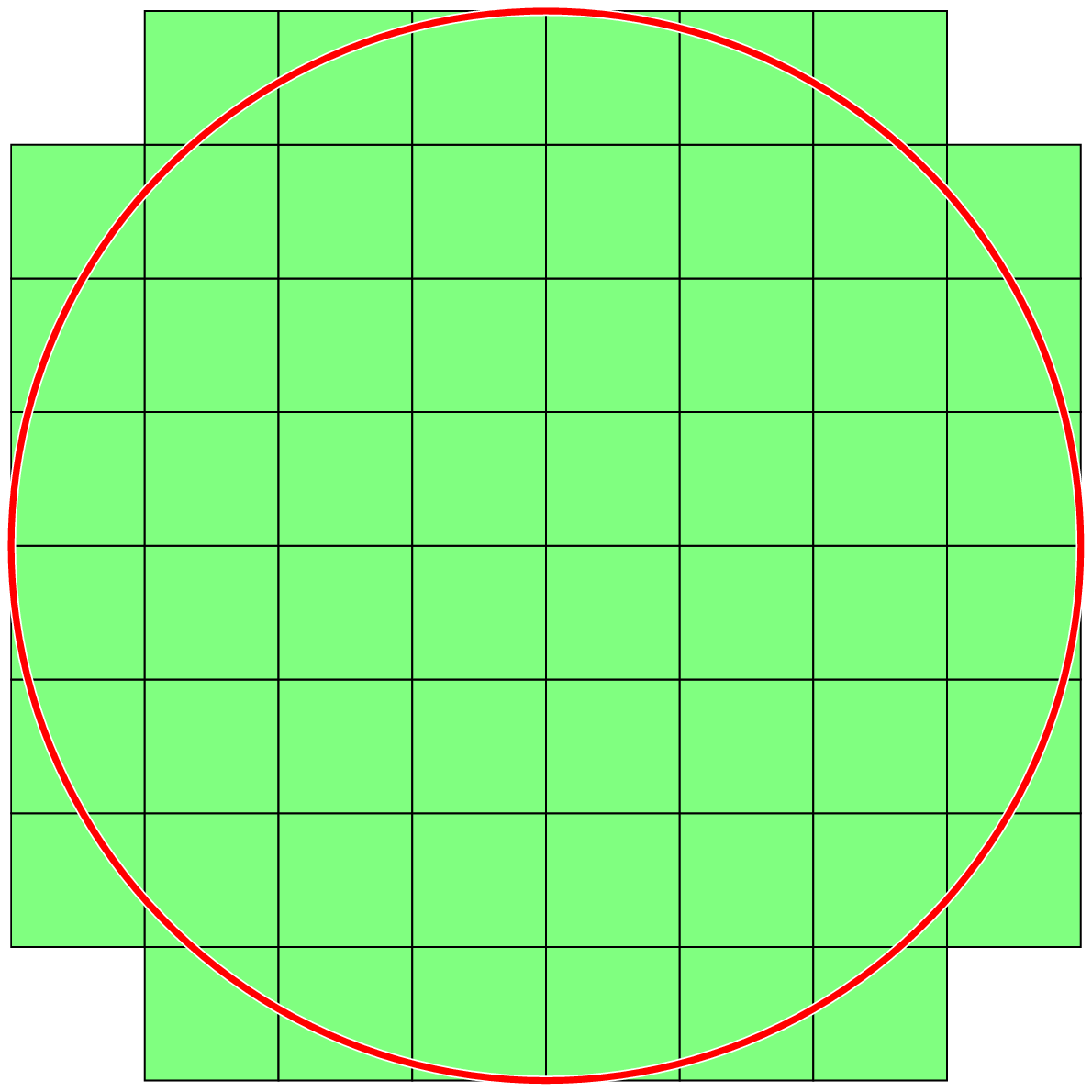}}
\subfigure[Step 3. The elements that intersect the surface $\Gamma$ are cut, thereby producing the purple band of polygonal elements.]{\label{fig:mesh_advantage_3}\includegraphics[scale=.25]{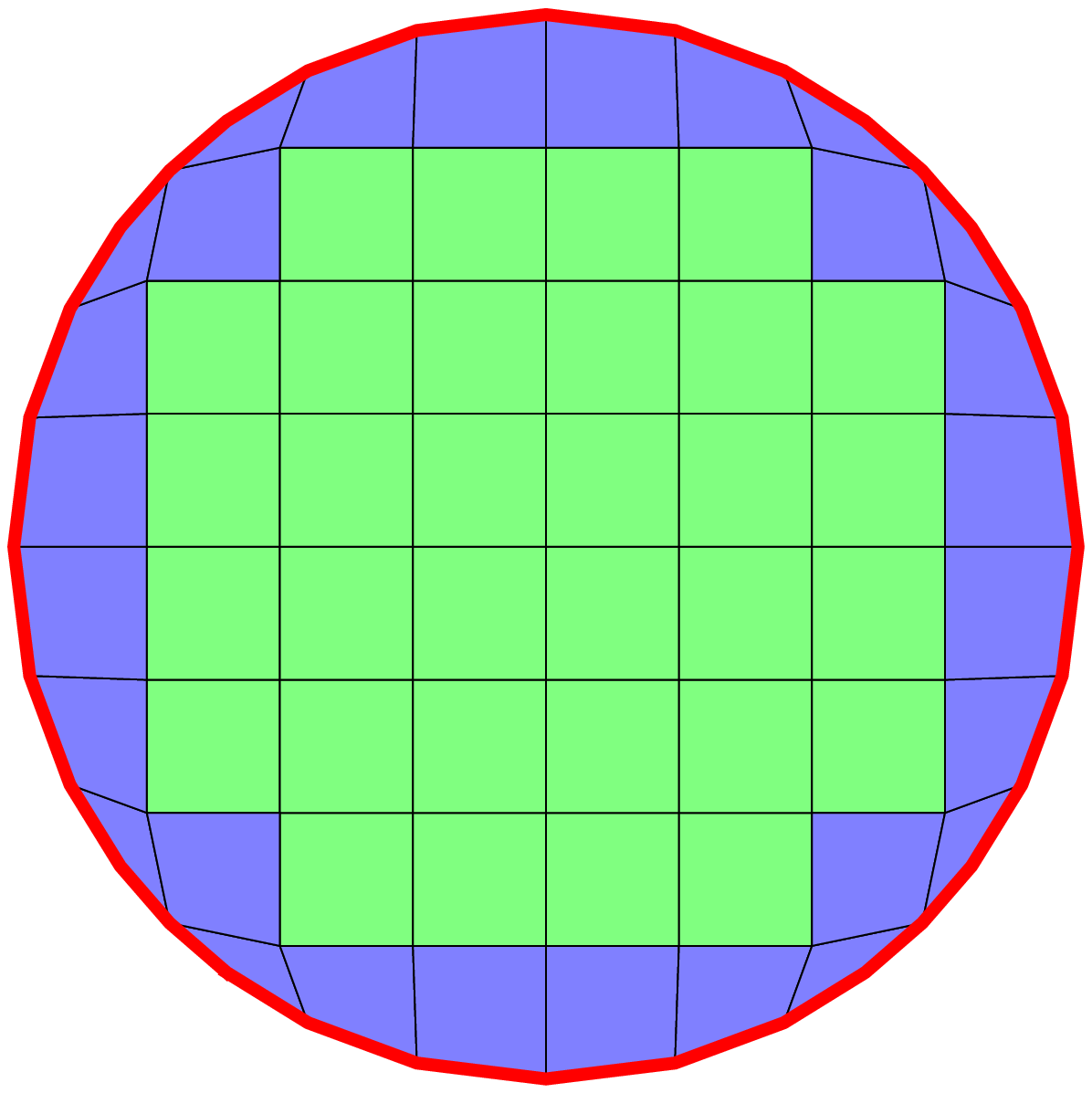}}
\end{center}
\caption{Generation of a polygonal BS mesh that allows for optimised matrix assembly.}
\label{fig:mesh_advantage}
\end{figure}

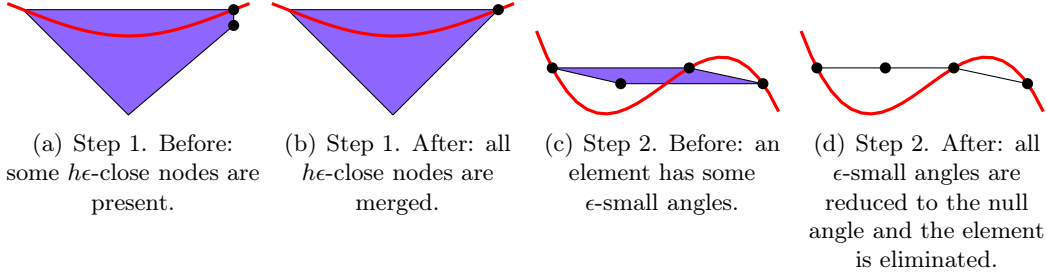
\begin{figure}[ht!]
\begin{center}
\subfigure[Step 1. Before: some $h\epsilon$-close nodes are present.]{\label{fig:post_processing_1_a}\input{post_processing_1_a.tex}}
\subfigure[Step 1. After: all $h\epsilon$-close nodes are merged.]{\label{fig:post_processing_1_b}\input{post_processing_1_b.tex}}
\subfigure[Step 2. Before: an element has some $\epsilon$-small angles.]{\label{fig:post_processing_2_a}\input{post_processing_2_a.tex}}
\subfigure[Step 2. After: all $\epsilon$-small angles are reduced to the null angle and the element is eliminated.]{\label{fig:post_processing_2_b}\input{post_processing_2_b.tex}}
\end{center}
\caption{{Enforcement of conditions \ref{A1}-\ref{A2} on the polygonal mesh for optimised matrix assembly. The neighbouring elements -not shown in the picture- are affected in the process, but properties \ref{A1}-\ref{A2} hold true regardless.}}
\label{fig:post-processing}
\end{figure}

As an example, given an arbitrarily shaped domain $\Omega$ with $\mathcal{C}^1$ boundary $\Gamma$, we construct a polygonal mesh designed for fast matrix assembly, by proceeding as follows. Suppose that the bulk $\Omega$ is contained in a square $Q$. We discretise $Q$ with a Cartesian grid made up of square mesh elements and assume that at least one of such squares is fully contained in the interior of $\Omega$ (Fig. \ref{fig:mesh_advantage_1}). Then we discard the elements that are outside $Q$ (Fig. \ref{fig:mesh_advantage_2}). Finally, we {project on $\Gamma$ the nodes that are outside $\Omega$}, thereby producing a discrete narrow band of irregular {quadrilaterals} (highlighted in purple in Fig. \ref{fig:mesh_advantage_3}). 
{The mesh is then post-processed in two steps as illustrated in Fig. \ref{fig:post-processing} to comply with Assumptions \ref{A1}-\ref{A2}. Let $\varepsilon > 0$ be a threshold. First, all $h\epsilon$-close nodes are merged, see Figs. \ref{fig:post-processing}(a)-(b). After that, the only chance for a polygon with at most four vertices fulfilling \ref{A1} to not fulfil \ref{A2} is the presence of small angles. In this case, all $\epsilon$-small angles are reduced to the null angle, see Figs. \ref{fig:post-processing}(c)-(d). It is worth remarking that Virtual Elements were proven to be robust with respect to distorted elements (see \cite{da2017high}), so $\epsilon$ can be chosen small.}

The resulting mesh $\Omega_h$ is made up of equal square elements, except for the elements that are close to $\Gamma$, as we can see in Fig. \ref{fig:mesh_advantage_3}, which results in faster matrix assembly. Let $h$ be the the meshsize of $\Omega_h$. By construction, $h = h_Q$, where $h_Q$ is the meshsize of the Cartesian grid. Of course, $\Omega_h$ is made up of $\mathcal{O}(1/h_Q^2) = \mathcal{O}(1/h^2)$ elements. However, by definition of Hausdorff dimension, the number of squares of $Q$ that intersect $\Gamma$ is $\mathcal{O}(1/h_Q) = \mathcal{O}(1/h)$, hence the number on non-square elements of $\Omega_h$ is only $\mathcal{O}(1/h)$. This implies that, when assembling the mass- and stiffness- matrices $M_\Omega$ and $A_\Omega$, respectively, only $\mathcal{O}(1/h)$ element-wise local matrices must be actually computed, since the local matrices for a square element are pre-computed.

It is worth remarking that the advantage described in this section becomes even more striking in higher space dimension: if the embedding space has dimension $D\in\mathbb{N}$, $D \geq 2$, then the procedure described in this section reduces the computational complexity of matrix assembly from $\mathcal{O}(1/h^D)$ to $\mathcal{O}(1/h^{D-1})$.  Matrix assembly optimization can be also achieved through alternative approaches, such as cut FEM \cite{burman2016cut} or trace FEM \cite{gross2015trace}. However, in these works, the authors adopt a level set representation of the boundary $\Gamma$, which we do not need in this study, as we exploit the usage of arbitrary polygons to approximate $\Gamma$.

\section{Numerical simulations}
\label{sec:numerical_examples}
{In this section we provide three numerical examples to compare the performances of BSFEM and BSVEM. In the first example, we show (i) optimal convergence in the case of the elliptic problem \eqref{elliptic_model} and (ii) the optimised matrix assembly introduced in Section \ref{sec:mesh_advantage}. In the second example, we show optimal convergence in space and time in the case of the parabolic problem \eqref{parabolic_problem}. In the third example, we compare the BSFEM and BSVEM solutions of the parabolic BS wave pinning model studied in \cite{cusseddu2019}.}

\subsection{Experiment 1: The elliptic problem}
\label{sec:experiment_elliptic}
\begin{figure}
\begin{center}
\subfigure[{Quasi-uniform Delaunay} triangular mesh of the kind used for BSFEM.]{\includegraphics[scale=.3]{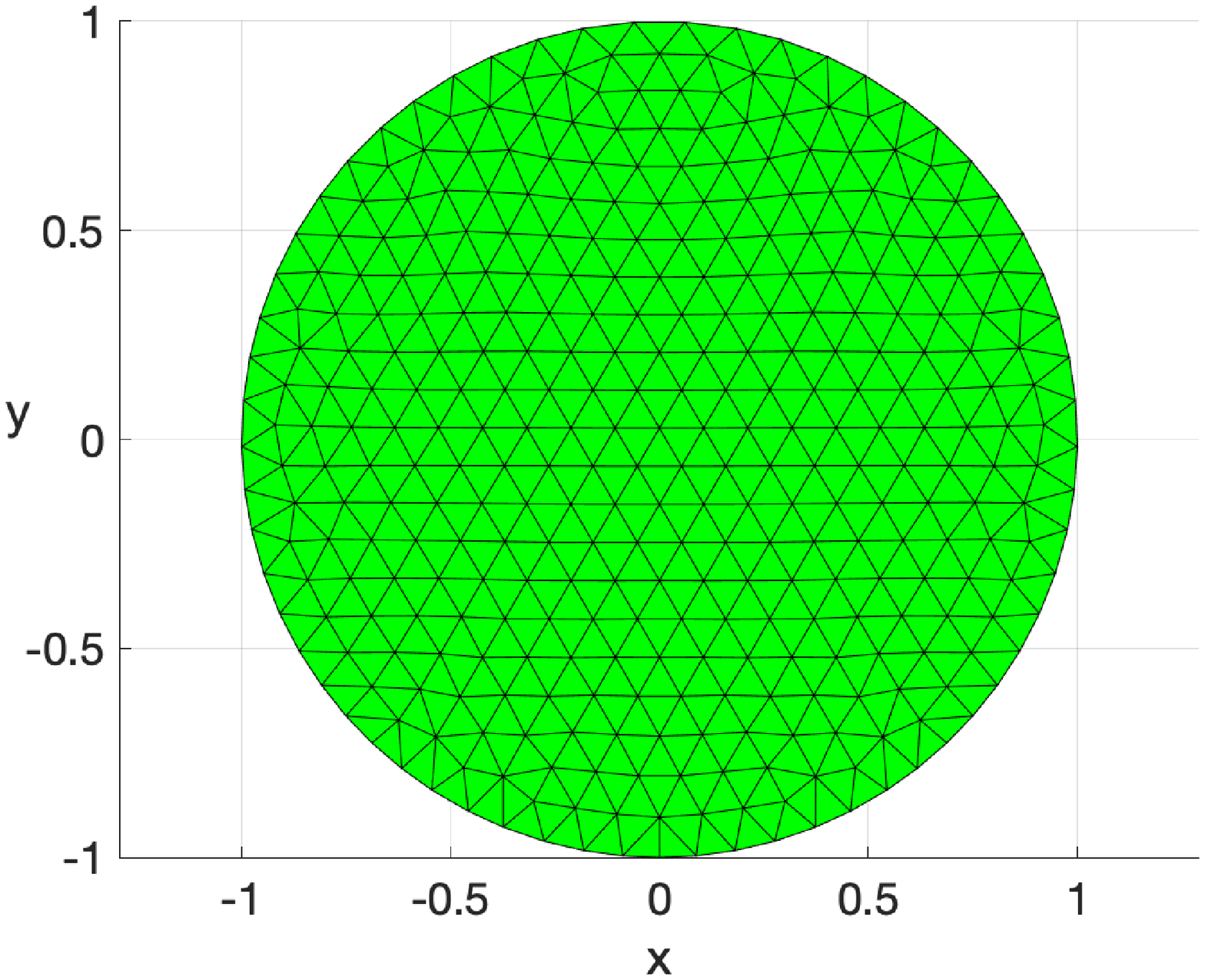}}
\subfigure[{Polygonal} mesh of the kind used for BSVEM.]{\includegraphics[scale=.3]{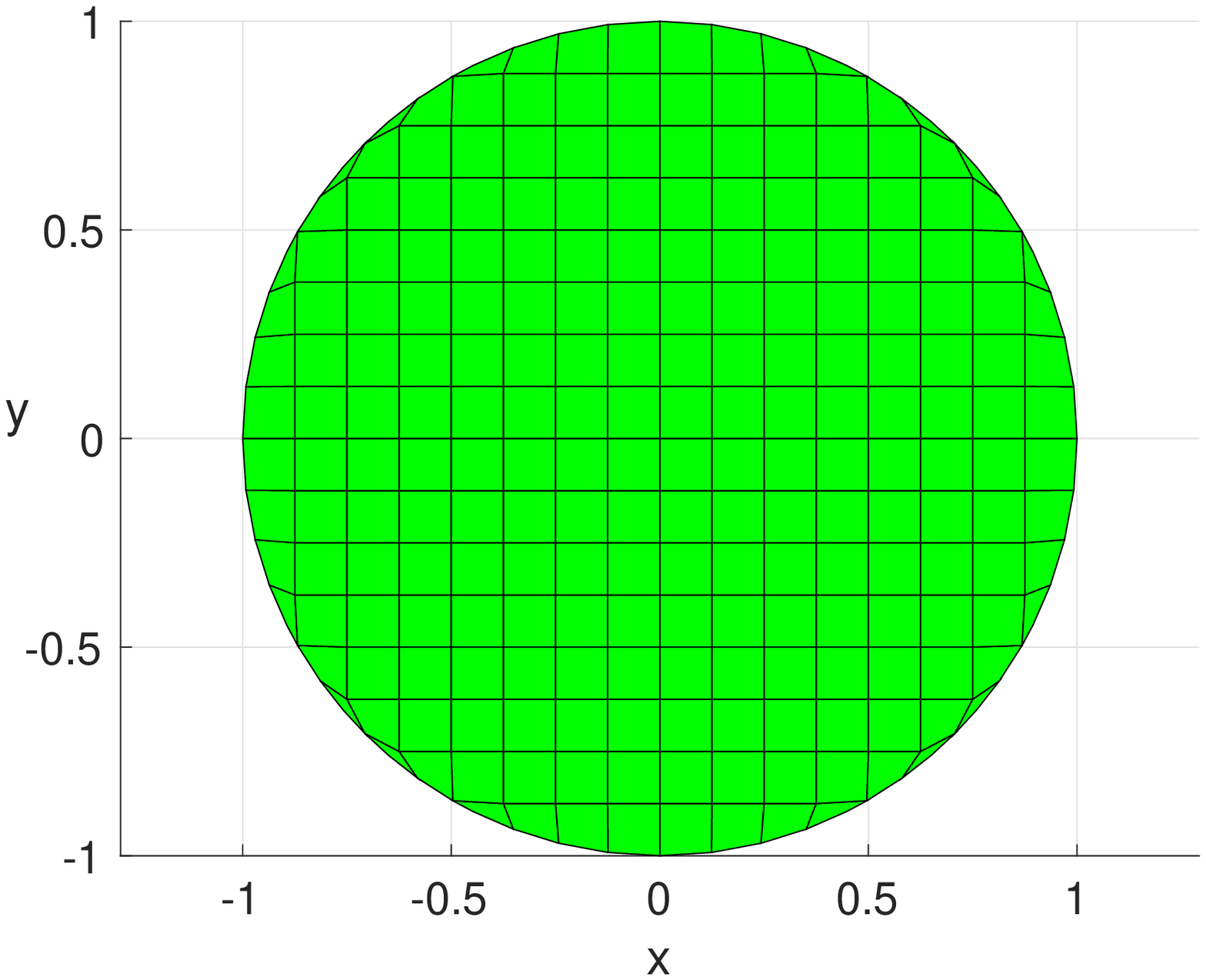}}
\end{center}
\caption{Illustrative representation of the meshes used for the numerical experiments.}
\label{fig:meshes}
\end{figure}

\begin{figure}
\begin{center}
\includegraphics[scale=.3]{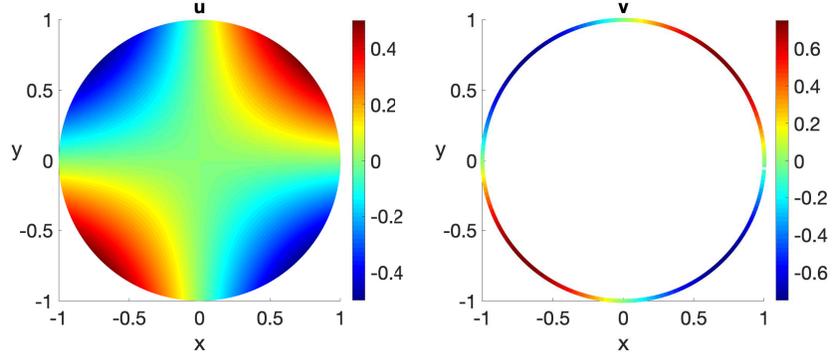}
\end{center}
\caption{Experiment 1 for the elliptic problem \eqref{elliptic_model}. Numerical solution on the finest of the meshes listed in Table \ref{tab:elliptic_table}, with meshsize $h = 4.7611$e-2. Left plot shows the solution $u$ in the bulk and the right plot shows the solution $v$ on the surface (curve) $\Gamma$.}
\label{fig:elliptic_solution}
\end{figure}

We start by constructing an exact solution for the elliptic problem \eqref{elliptic_model} on the unit circle $\Omega := \{(x,y)\in\mathbb{R}^2 | x^2 + y^2 = 1\}$ by using the fact that $-\Delta_\Gamma xy = {4}xy$, i.e. the function $w(x,y) := xy$ is an eigenfunction of the Laplace-Beltrami operator on the unit circumference $\Gamma = \partial \Omega$. Specifically, in \eqref{elliptic_model} we choose the following load terms $f(x,y) := xy$, for $(x,y) \in \Omega$ and $g(x,y) := \left(2+\frac{5(\alpha+2)}{\beta}\right)xy$ for $(x,y) \in \Gamma$. Here $\alpha$ and $\beta$ are the parameters that appear in the model, such that the exact solution is given by $u(x,y) := xy$ for $(x,y) \in \Omega$ and $v(x,y) := \frac{\alpha+2}{\beta}xy$ for $(x,y) \in \Gamma$. For illustrative purposes, we choose $\alpha=1$ and $\beta =2$. We solve the problem {with BSFEM and BSVEM on two respective sequences of five meshes} with decreasing meshsizes. {For BSVEM, we use optimised meshes as in Section \ref{sec:mesh_advantage}, while for BSFEM we use quasi-uniform Delaunay triangular meshes. See Fig. \ref{fig:meshes} for an illustration.} For each mesh, we measure:
\begin{itemize}
\item the {$L^2(\Omega)\times L^2(\Gamma)$ and $L^\infty(\Omega)\times L^\infty(\Gamma)$} errors of the numerical solution and the respective experimental orders of convergence (EOCs), that are quadratic in the meshsize. {In this case BSVEM is almost three times as accurate as BSFEM on similar meshsizes, probably because the mesh reflects the symmetry of the problem};
\item {the condition number $\text{cond}_{ell}$ of the linear system \eqref{elliptic_problem_BSVEM_linear_system_compact}. The BSVEM is approximately four times as ill conditioned as BSFEM on similar meshsizes. An approach to reduce this gap is described in \cite{berrone2017orthogonal}};
\item {only for BSVEM,} the number $N_\Omega$ of all elements and the number $N_\Gamma$  of elements that intersect the boundary. The first is proportional to $\frac{1}{h^2}$, while the latter is proportional to $\frac{1}{h}$, as predicted in Section \ref{sec:mesh_advantage}. This illustrates that only $O(\frac{1}{h})$ local matrices must be computed during matrix assembly, even if the mesh has $O(\frac{1}{h^2})$ elements.
\end{itemize}
In Fig. \ref{fig:elliptic_solution} we show the {BSVEM} numerical solution on the finest mesh. In Tabs. \ref{tab:elliptic_table}-\ref{tab:elliptic_table_bsfem} we show the errors, the respective EOCs, and the numbers $N_\Omega$, $N_\Gamma$ and $\text{cond}_{ell}$.

\begin{table}
\caption{Experiment 1 for the elliptic problem \eqref{elliptic_model}. For BSVEM, EOCs in {both $L^2(\Omega)\times L^2(\Gamma)$ and $L^\infty(\Omega)\times L^\infty(\Gamma)$ norms are} approximately two. Moreover, matrix assembly requires the computation of $N_\Gamma = O(\frac{1}{h})$ local matrices out of $N_\Omega = O(\frac{1}{h^2})$ elements.}
\begin{center}
\input{elliptic_table.tex}
\end{center}
\label{tab:elliptic_table}
\end{table}

\begin{table}
\caption{Experiment 1 for the elliptic problem \eqref{elliptic_model}, solved with BSFEM on a sequence of triangular meshes with meshsizes as similar as possible as those in Tab. \ref{tab:elliptic_table}. The EOCs in {both $L^2(\Omega)\times L^2(\Gamma)$ and $L^\infty(\Omega)\times L^\infty(\Gamma)$ norms are approximately two and the errors almost three times as large w.r.t. the corresponding errors of BSVEM}.}
\begin{center}
\input{elliptic_table_bsfem.tex}
\end{center}
\label{tab:elliptic_table_bsfem}
\end{table}

\subsection{Experiment 2: Linear parabolic problem}
We now want to test convergence for the parabolic problem \eqref{parabolic_problem} in the linear case $q(u) = -u$, $r(u,v) = 2u$, $s(u,v) = \frac{4}{3}v$, $d_u = 1$ and $d_v = \frac{1}{4}$, on the unit circle $\Omega := \{(x,y) \in\mathbb{R}^2 | x^2 +  y^2 \leq 1\}$. We consider the initial conditions $u(x,y,0) = xy$ and $v(x,y,0) = \frac{3}{2}xy$, such that the exact solution $u(x,y,t) = e^{-t}xy$ and $v(x,y,t) = \frac{3}{2}e^{-t}xy$. We solve the problem via BSFEM- and BSVEM-IMEX Euler on the same sequences of meshes of Experiment 1 with coresponding timesteps $\tau_i = 10^{-3} \times 4^{-i}$, $i=0,\dots,4$ and final time $T=1$. The errors in $L^\infty([0,T], L^2(\Omega)\times L^2(\Gamma))$ and $L^\infty([0,T], L^\infty(\Omega)\times L^\infty(\Gamma))$ norms and the respective EOCs are shown in Tabs. \ref{tab:parabolic_table}-\ref{tab:parabolic_table_bsfem}. {For both methods we observe superconvergence, with the BSVEM being more accurate on similar meshsizes.}

\begin{table}
\caption{{Experiment 2 for the parabolic problem \eqref{parabolic_problem}. The EOCs of BSVEM in both $L^\infty([0,T], L^2(\Omega)\times L^2(\Gamma))$ and $L^\infty([0,T], L^\infty(\Omega)\times L^\infty(\Gamma))$ norms appear to be more than two, probably due to the symmetry of the problem. The above norms are abbreviated in the table by $L^\infty(L^2)$ and $L^\infty(L^\infty)$, respectively.}}
\begin{center}
\input{parabolic_table.tex}
\end{center}
\label{tab:parabolic_table}
\end{table}

\begin{table}
\caption{{Experiment 2 for the parabolic problem \eqref{parabolic_problem}. The BSFEM still exhibits superconvergence, but provides less accurate solutions than BSVEM.}}
\begin{center}
\input{parabolic_table_bsfem.tex}
\end{center}
\label{tab:parabolic_table_bsfem}
\end{table}

\subsection{Experiment 3: Nonlinear parabolic problem}
\label{sec:experiment_parabolic}
\begin{figure}
\begin{center}
\subfigure[BSVEM interpolant of the initial condition.]{\includegraphics[scale=.25]{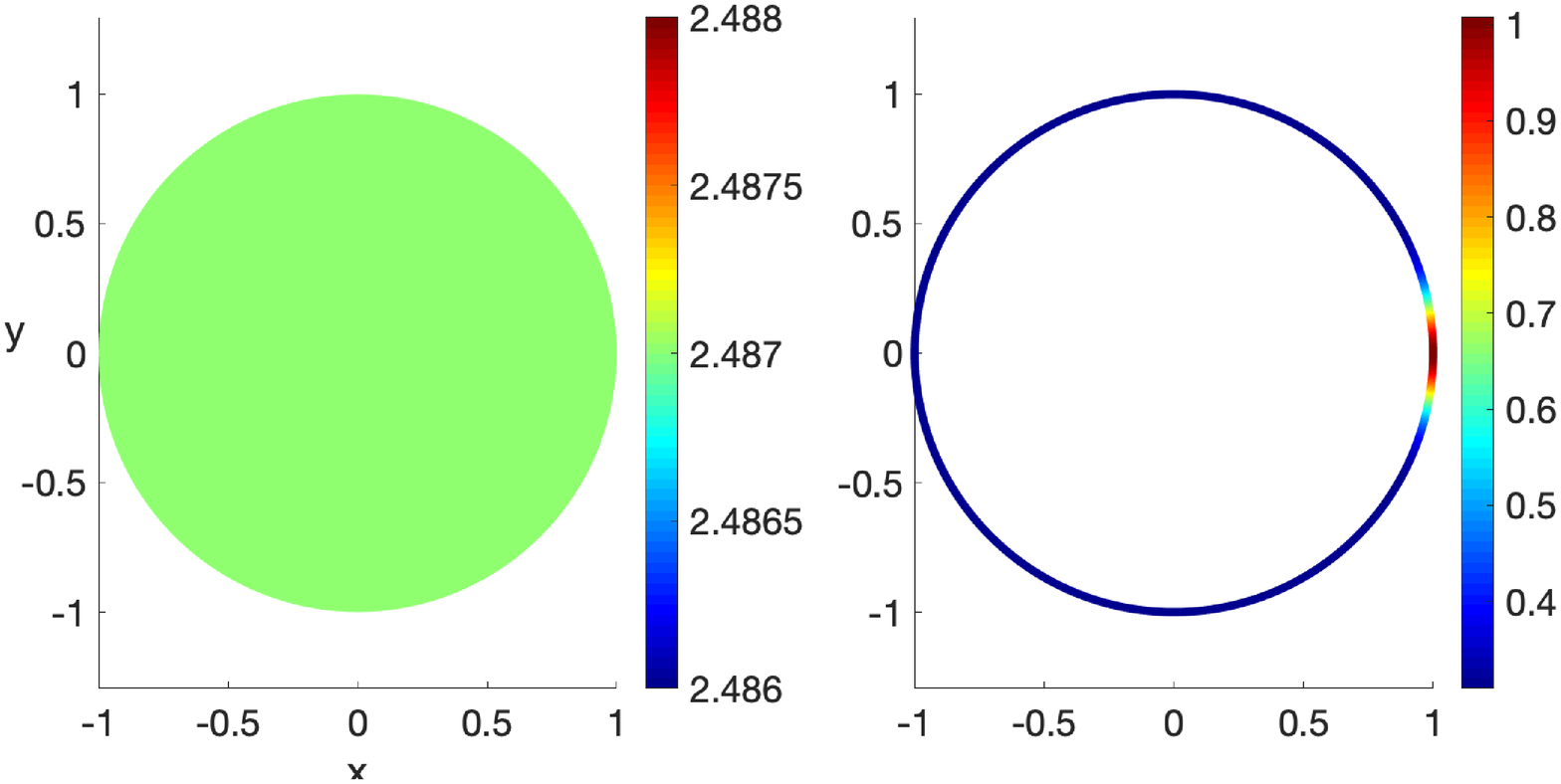}}
\subfigure[BSVEM solution at time $t = 0.1$.]{\includegraphics[scale=.25]{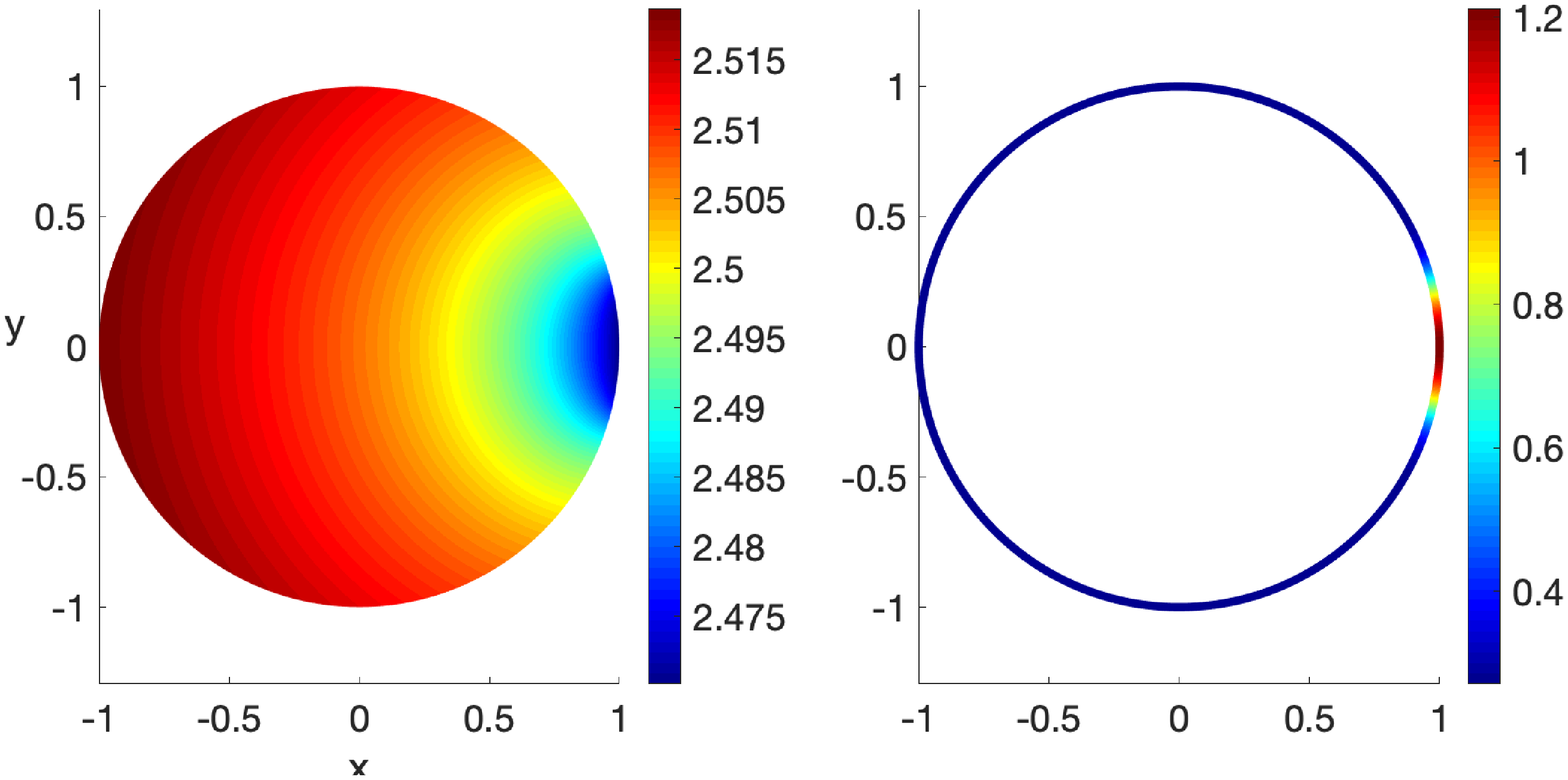}}
\subfigure[BSVEM solution at time $t = 1$.]{\includegraphics[scale=.25]{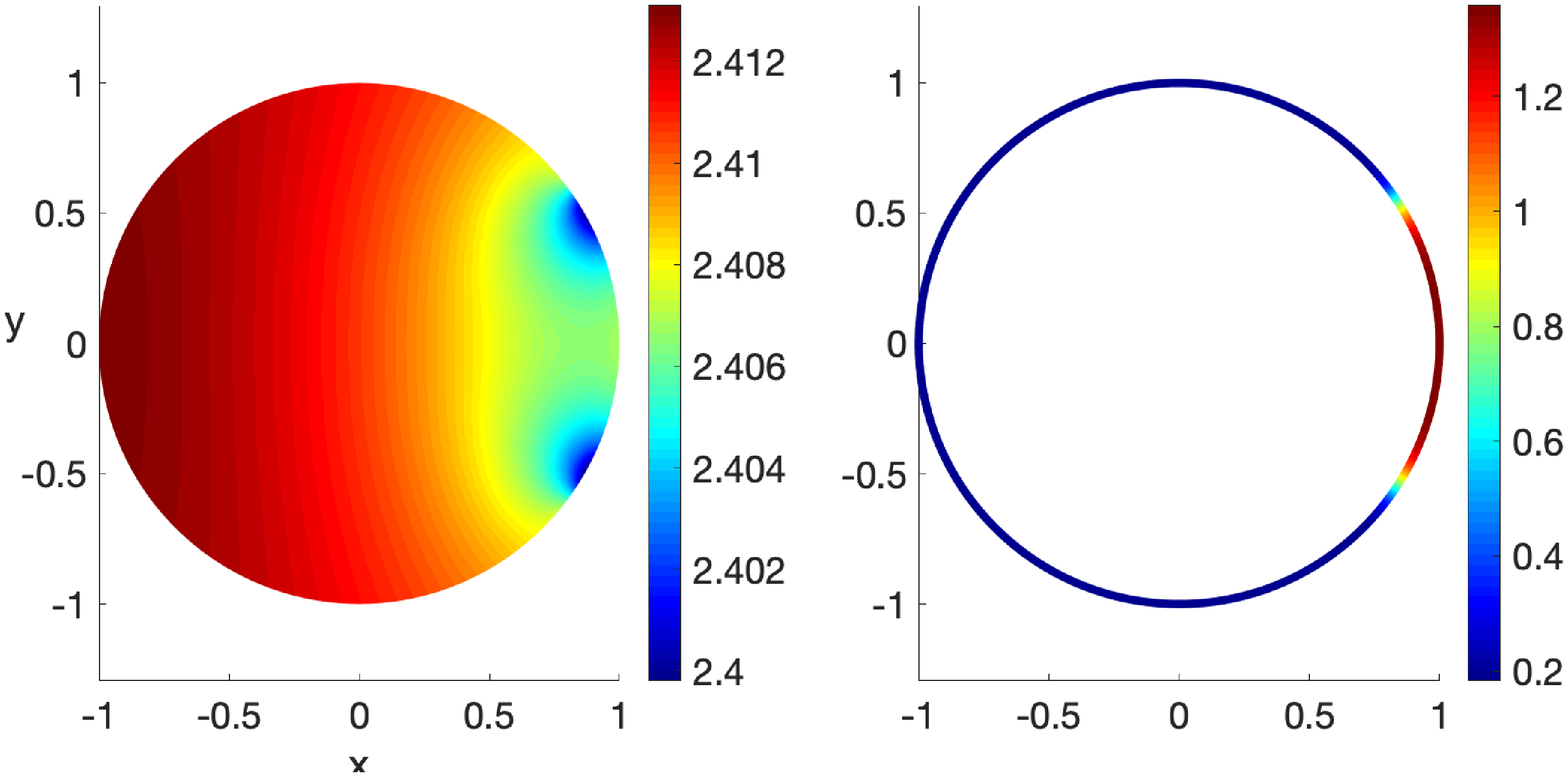}}
\subfigure[BSVEM solution at time $t = 4.5$.]{\includegraphics[scale=.25]{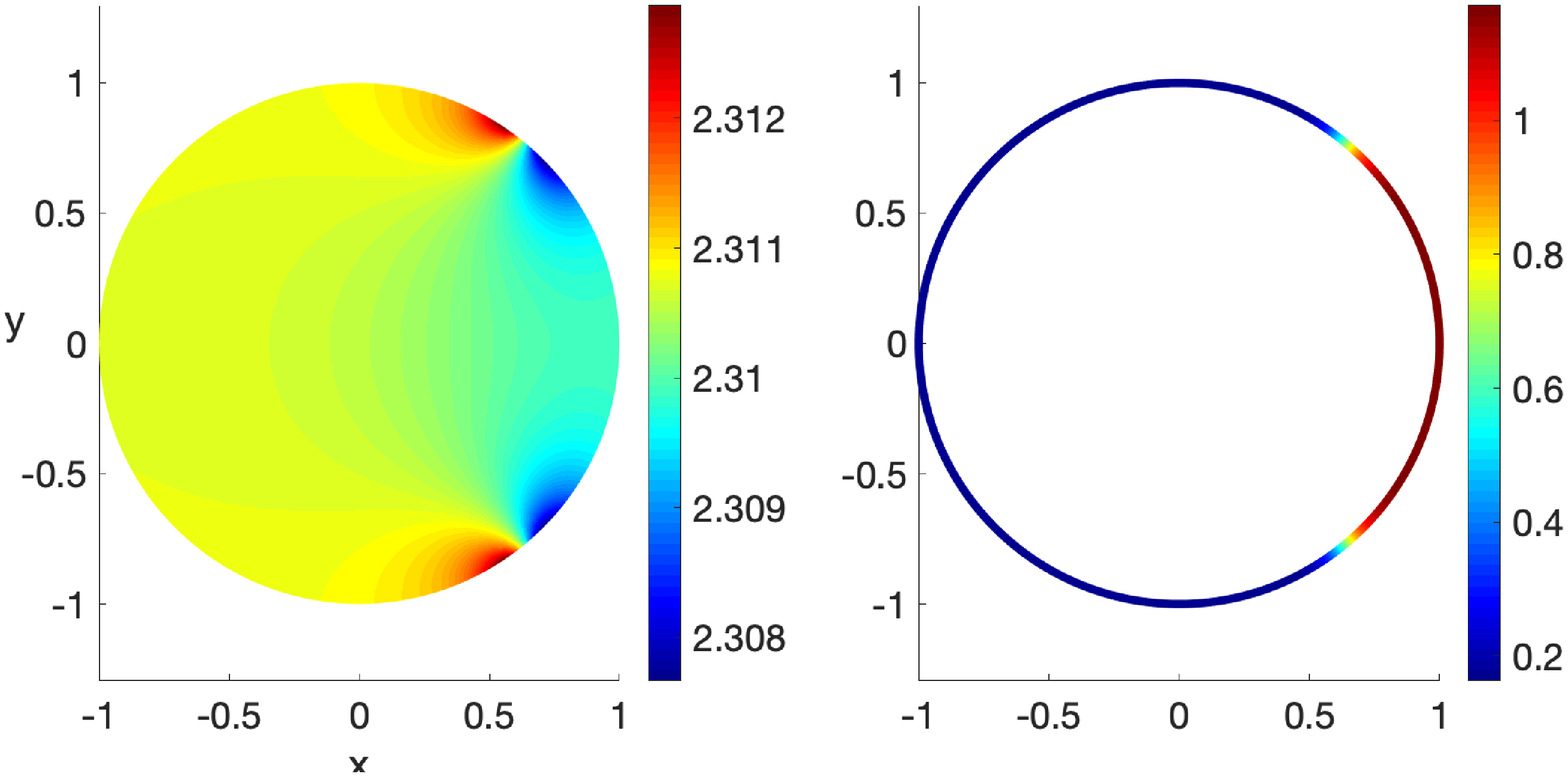}}
\end{center}
\caption{{Experiment 3: BSVEM solutions} of the BSWP model \eqref{BSWP_model} {at different times.}}
\label{fig:simulations}
\end{figure}

In this section we compare the {BSFEM- and BSVEM-IMEX Euler} numerical solutions of the bulk-surface wave pinning (BSWP) model considered in \cite[Fig. 4]{cusseddu2019}. The non-dimensionalised BSWP model seeks to find a bulk concentration $b: \Omega \times [0,T] \rightarrow\mathbb{R}$ and a surface concentration $a:\Gamma\times[0,T] \rightarrow\mathbb{R}$ such that
\begin{equation}
\label{BSWP_model}
\begin{cases}
\vspace*{2mm}
\displaystyle\varepsilon\frac{\partial b}{\partial t} - \Delta b = 0, \qquad \boldx \in \Omega,\\
\vspace*{2mm}
\displaystyle\varepsilon\frac{\partial a}{\partial t} - \varepsilon^2 \Delta_\Gamma a = f(a,b), \qquad  \boldx  \in  \Gamma,\\
-(\boldnu\cdot\nabla b) = f(a,b), \qquad \boldx\in\Gamma,
\end{cases}
\end{equation}
where the kinetic $f$ is of the form $f(a,b) := \left(k_0 + \gamma {(a^2)}/{(1+a^2)} \right)b - a$. We solve the BSWP model \eqref{BSWP_model} on the unit circle $\Omega := \{(x,y) \in\mathbb{R}^2 | x^2 +  y^2 \leq 1\}$. We choose the parameters of the model as $\varepsilon^2 = 0.001$, $k_0 = 0.05$ and $\gamma = 0.79$. The initial condition, plotted in Fig. \ref{fig:simulations}(a), is $b(x,y,0) = 2.487$, for $(x,y)\in\Omega$ and
$a(x,y,0) = 0.309 + 0.35(1 + sign(x))\exp(-20 y^2)$ for $(x,y) \in \Gamma$. The final time is $T=4.5$. 
For the BSFEM we take a {quasi-uniform Delaunay} triangular mesh, while for the BSVEM we take a polygonal mesh designed for optimised matrix assembly, as explained in Section \ref{sec:mesh_advantage}. An illustrative coarser representation of the meshes typically used is shown in Fig. \ref{fig:meshes}. The details of these meshes are reported in Table \ref{tab:mesh_comparison}. As we can see, on almost equal meshsizes, the BSVEM generates (i) a significantly large number of boundary nodes, which translates into better boundary approximation and (ii) less number of elements in the bulk, which implies faster matrix assembly. 
For both spatial discretisations, the time discretisation is computed with timestep $\tau= 2e$-3. {The condition numbers $\text{cond}_\Omega := \text{cond}(M_\Omega + \tau {d_u}K_\Omega)$ and $\text{cond}_\Gamma := \text{cond}(M_\Gamma + \tau {d_v}K_\Gamma)$ of BSVEM on our polygonal mesh are only about $2.5$ times as large as those of BSFEM on a quasi-uniform Delaunay mesh. This result, which could be further improved through an orthogonal polynomial basis approach \cite{berrone2017orthogonal} or through a different time integrator, confirms the well-known robustness of VEM with respect to general polygonal meshes \cite{da2017high}. The BSVEM solution at various times is shown in Fig. \ref{fig:simulations}.} 

\begin{table}
\caption{{Experiment 3 for the wave pinning model \eqref{BSWP_model}:} Details of the meshes used for the BSFEM and BSVEM, respectively. On almost equal meshsizes, the BSVEM mesh has significantly more boundary nodes (providing better boundary approximation) and less bulk elements (which simplifies matrix assembly) at the expense of a slightly higher overall number of bulk nodes {and slightly worse conditioning} compared to the BSFEM.}
\begin{center}
\begin{tabular}{lllllll}
\hline\noalign{\smallskip}
Spatial method & $h$ (meshsize) & \makecell{$N$ (number\\ of nodes)} & \makecell{$M$ (number of\\ boundary nodes)} & \makecell{$|\mathcal{E}_h|$ (number\\ of elements)} & $\text{cond}_\Omega$ & $\text{cond}_\Gamma$\\
\noalign{\smallskip}\hline\noalign{\smallskip}
BSFEM & $3.10e$-2 & $5809$ & $250$ & $11366$ & $8.91e$+2 & $1.55$\\
BSVEM & $3.19e$-2 & $6536$ & $356$ & $6357$ & $2.21e$+3 & $3.21$\\
\noalign{\smallskip}\hline
\end{tabular}
\end{center}
\label{tab:mesh_comparison}
\end{table}

\section{Conclusions}
In this study, we have considered a bulk-surface virtual element method (BSVEM) for the numerical approximation of linear elliptic and semilinear parabolic coupled BSPDE problems on smooth BS domains. The proposed method simultaneously extends the BSFEM for BSRDSs \cite{Madzvamuse_2016} and the VEM for linear elliptic \cite{beirao2013basic} and semilinear parabolic \cite{adak2019convergence} bulk PDEs. {The method has a simplified vector-matrix form that can be exploited also in the special case of the BSFEM considered in \cite{Madzvamuse_2016}.}

{We have introduced polygonal BS meshes in two space dimensions and we have shown that the geometric error arising from domain approximation is $\mathcal{O}(h)$ in the bulk and $\mathcal{O}(h^2)$ on the surface, where $h$ is the meshsize, exactly as in the special case of the triangular BS meshes used in the BSFEM \cite{elliottranner2013finite}. Suitable polygonal meshes reduce the asymptotic computational complexity of matrix assembly from $\mathcal{O}(1/h^2)$ to $\mathcal{O}(1/h)$. In future studies we will show that polygonal meshes can unlock new efficient adaptive refinement strategies.}

We have introduced novel theory to address the challenges of the error analysis. {First, if the exact solution $(u,v)$ is $H^{2+1/4}(\Omega)\times H^2(\Gamma)$, the lifting operator can be replaced with the Sobolev extension operator. Second, if $\Gamma$ is sufficiently smooth, there exists an $L^2$-\emph{preserving} inverse trace operator that preserves the $L^2$ norm and the $H^1$ seminorm up to the same scale factor. Third, we have introduced a tailor-made Ritz projection in the bulk for the VEM that accounts for boundary conditions and that fulfils optimal error bounds. By using our bulk-Ritz projection we have drawn two consequences.} First, the lowest order bulk-VEM \cite{beirao2013basic} retains optimal convergence even in the simultaneous presence of curved boundaries and non-zero boundary conditions, a result that lacked a rigorous proof in the literature. Second, the proposed BSVEM possesses optimal spatial convergence, that is $\mathcal{O}(h^2)$ in $L^\infty([0,T], L^2(\Omega)\times L^2(\Gamma))$ norm, where $T$ is the final time.

{Numerical examples validate our findings in terms of (i) convergence rate in space and time for both the elliptic and the parabolic case and (ii) the computational advantages of polygonal meshes.} {Given the encouraging results, the extension of BSVEM to higher space dimension and different BSPDE problems is under development. Moreover, the generalisation of problem \eqref{parabolic_problem} to \emph{evolving} BS domains, which would comprise additional models addressed in the literature such as \cite{MacDonald_2016}, is left for future work.}

\section*{Acknowledgements and Funding}
MF’s research was partially funded by the Italian National Group of Scientific Computing (GNCS-INdAM), by the University of Sussex and by Regione Puglia (Italy) through the research programme REFIN -Research for Innovation (protocol code 901D2CAA, project number UNISAL026). AM's work was partially funded by the EPSRC grants (EP/J016780/1, EP/T00410X/1), the Leverhulme Trust Research Project Grant (RPG-2014-149), the European Union Horizon 2020 research and innovation programme under the Marie Sk\l{}odowska-Curie grant agreement (No 642866), the Commission for Developing  Countries, and was partially supported by a grant from the Simons Foundation. AM is a Royal Society Wolfson Research Merit Award Holder funded generously by the Wolfson Foundation. AM is a Distinguished Visiting Scholar to the Universit\'{a} degli Studi di Bari Aldo Moro, Bari, Italy and the Department of Mathematics, University of Johannesburg, South Africa.\\
The authors want to thank Professor Endre S\"{u}li (University of Oxford, UK) for the critical review of this work.

\section*{Conflict of interest}
The authors declare that they have no conflict of interest.

\bibliographystyle{plainurl}
\bibliography{bibliography}   

\appendix
\section*{Appendix A: Preliminary definitions and results}
\setcounter{equation}{0}
\renewcommand{\theequation}{A.\arabic{equation}}
\setcounter{theorem}{0}
\renewcommand{\thetheorem}{A\arabic{theorem}}
\setcounter{definition}{0}
\renewcommand{\thedefinition}{A\arabic{definition}}
\setcounter{lemma}{0}
\renewcommand{\thelemma}{A\arabic{lemma}}
\setcounter{remark}{0}
\renewcommand{\theremark}{A\arabic{remark}}
\setcounter{corollary}{0}
\renewcommand{\thecorollary}{A\arabic{corollary}}

In this Appendix we provide preliminary definitions, results and notations adopted throughout the article. Unless explicitly stated, definitions and results are taken from \cite{dziuk2013finite}.

\subsection*{A1. Surfaces and differential operators on surfaces}
Let $\Omega \subset \mathbb{R}^2$ be a compact set such that its boundary $\Gamma := \partial \Omega \subset \mathbb{R}^2$ is a $\mathcal{C}^k$, $k \geq 2$ curve. Since $\Gamma$ can be regarded as the zero level set of the \emph{oriented distance function} $d:\mathbb{R}^2 \rightarrow\mathbb{R}$ defined by
\begin{equation*}
d(\boldx) := 
\begin{cases}
-&\inf \{\|\boldx-\boldy\| : \boldy\in \Gamma\} \qquad \text{if } \boldx \in \Omega;\\
&0 \hspace*{27mm} \text{if } \boldx \in\Gamma;\\
&\inf \{\|\boldx-\boldy\| : \boldy \in \Gamma\} \qquad \text{if } \boldx \in \mathbb{R}^2 \setminus \Omega,
\end{cases}
\end{equation*}
then the outward unit vector field $\boldnu:\Gamma \rightarrow\mathbb{R}^3$ can be defined by
\begin{equation}
\label{normal_vector_field}
\boldnu(\boldx) := \frac{\nabla d(\boldx)}{\|\nabla d(\boldx)\|}, \qquad \boldx\in\Gamma,
\end{equation}
see for instance \cite{elliottranner2013finite}.

\begin{lemma}[Fermi coordinates]
\label{lmm:fermi}
If $\Gamma$ is a $\mathcal{C}^k$, $k\geq 2$ surface, there exists an open neighbourhood $U \subset \mathbb{R}^3$ of $\Gamma$ such that every $\boldx\in U$ admits a unique decomposition of the form
\begin{equation*}
\boldx = \bolda(\boldx) + d(\boldx)\boldnu(\bolda(\boldx)), \qquad  \bolda(\boldx) \in \Gamma.
\end{equation*}
The maximal open set $U$ with this property is called the \emph{Fermi stripe} of $\Gamma$ (see Fig. \ref{fig:exact_domain}), $\bolda(\boldx)$ is called the \emph{normal projection} onto $\Gamma$ and $(\bolda(\boldx), d(\boldx))$ are called the \emph{Fermi coordinates} of $\boldx$. The oriented distance function fulfils $d\in \mathcal{C}^k(U)$.
\begin{proof}
See \cite{dziuk2013finite}.
\end{proof}
\end{lemma}

\begin{definition}[$\mathcal{C}^1(\Gamma)$ functions]
A function $u:\Gamma \rightarrow\mathbb{R}$ is said to be $\mathcal{C}^1(\Gamma)$ if there exist an open neighbourhood $U$ of $\Gamma$ and a $\mathcal{C}^1$ function $\hat{u}:U\rightarrow\mathbb{R}$ such that $\hat{u}_{|\Gamma}  = u$, i.e. $\hat{u}$ is a $\mathcal{C}^1$ extension of $u$ off $\Gamma$.
\end{definition}

\begin{definition}[Tangential gradient and tangential derivatives]
The \emph{tangential gradient} $\nabla_\Gamma u$ of a function $u\in\mathcal{C}^1(\Gamma)$ is defined by
\begin{equation}
\label{tangential_gradient}
\nabla_\Gamma u(\boldx) := \nabla \hat{u}(\boldx) - (\nabla \hat{u}(\boldx) \cdot \boldnu(\boldx)) \boldnu(\boldx), \qquad \boldx \in \Gamma.
\end{equation}
The result of the computation in \eqref{tangential_gradient} is independent of the choice of the extension $\hat{u}$. The components $D_x u$ and $D_y u$ of the tangential gradient $\nabla_\Gamma u$ are called the \emph{tangential derivatives} of $u$.
\end{definition}

\begin{definition}[$\mathcal{C}^k(\Gamma)$ functions]
For $k\in\mathbb{N}$, $k>1$, a function $u:\Gamma \rightarrow\mathbb{R}$ is said to be $\mathcal{C}^k(\Gamma)$ if it is $\mathcal{C}^1(\Gamma)$ and its tangential derivatives are $\mathcal{C}^{k-1}(\Gamma)$.
\end{definition}

\begin{remark}[Regularity of normal projection]
\label{rmk:regularity_normal_projection}
Consider the function
\begin{equation*}
F(\boldx,\boldy) := \boldy - \boldx + d(\boldx)\boldnu(\boldy), \qquad \boldx \in U, \ \boldy \in \Gamma.
\end{equation*}
Since $d \in \mathcal{C}^k(U)$ from Lemma \ref{lmm:fermi} and $\boldnu \in \mathcal{C}^{k-1}(\Gamma)$ from \eqref{normal_vector_field}, then $F \in \mathcal{C}^{k-1}(\Omega \times \Gamma)$. Since the normal projection $\bolda(\boldx)$ of $\boldx\in U$ can be regarded as the solution of the implicit equation $F(\boldx, \bolda(\boldx)) = 0$, then $\bolda \in \mathcal{C}^{k-1}(U)$ as well.
\end{remark}

\begin{definition}[Laplace-Beltrami operator]
The \emph{Laplace-Beltrami operator} $\Delta_\Gamma u$ of a function $u \in \mathcal{C}^2(\Gamma)$ is defined by
\begin{equation*}
\Delta_\Gamma u(\boldx) := D_xD_x u(\boldx) + D_yD_y u(\boldx), \qquad \boldx\in\Gamma.
\end{equation*}
\end{definition}

\subsection*{A2. Bulk- and surface function spaces}
\begin{definition}[Lebesgue function spaces]
Let $p\in [1,+\infty]$. For $u:\Omega \rightarrow \mathbb{R}$ and $v:\Gamma\rightarrow\mathbb{R}$ the bulk- and surface Lebesgue norms are defined by
\begin{equation*}
\|u\|_{L^p(\Omega)} :=
\begin{cases}
&\displaystyle\left(\intomega |u|^p\right)^{1/p} \quad \text{if}\ p\in [1,+\infty);\\
&\displaystyle\ess\sup_{\boldx\in\Omega} |u(\boldx)| \quad \text{if}\ p=+\infty,
\end{cases}
\qquad
\|v\|_{L^p(\Gamma)} :=
\begin{cases}
&\displaystyle\left(\intgamma |v|^p\right)^{1/p} \quad \text{if}\ p\in [1,+\infty);\\
&\displaystyle\ess\sup_{\boldx\in\Gamma} |v(\boldx)| \quad \text{if}\ p=+\infty,
\end{cases}
\end{equation*}
respectively. The bulk- and surface Lebesgue function spaces are defined as
\begin{align*}
&L^p(\Omega) := \{u:\Omega\rightarrow\mathbb{R}\ |\ \|u\|_{L^p(\Omega)} < +\infty\};\\
&L^p(\Gamma) := \{v:\Gamma\rightarrow\mathbb{R}\ |\ \|v\|_{L^p(\Gamma)} < +\infty\},
\end{align*}
respectively.
\end{definition}

\begin{definition}[Sobolev function spaces]
Let $m\in\mathbb{N}$ and $p\in[1,+\infty]$. For $u:\Omega \rightarrow \mathbb{R}$ and $v:\Gamma\rightarrow\mathbb{R}$, the bulk- and surface Sobolev norms are defined by
\begin{align}
\label{integer_sobolev_norm_bulk}
\|u\|_{W^{m,p}(\Omega)} := 
\begin{cases}
&\displaystyle\left(\sum_{|\boldalpha| \leq m} \|D_{\boldalpha} u\|_{L^p(\Omega)}^p\right)^{1/p} \quad \text{if}\ p\in [1,+\infty);\\
&\displaystyle\max_{|\boldalpha| \leq m} \|D_{\boldalpha} u\|_{L^\infty(\Omega)} \quad \text{if}\ p = +\infty,
\end{cases}\\
\label{integer_sobolev_norm_surf}
\|v\|_{W^{m,p}(\Gamma)} := 
\begin{cases}
&\displaystyle\left(\sum_{|\boldalpha| \leq m} \|D_{\boldalpha} v\|_{L^p(\Gamma)}^p\right)^{1/p} \quad \text{if}\ p\in [1,+\infty);\\
&\displaystyle\max_{|\boldalpha| \leq m} \|D_{\boldalpha} v\|_{L^\infty(\Gamma)} \quad \text{if}\ p = +\infty,
\end{cases}
\end{align}
respectively. The seminorms $|u|_{W^{m,p}(\Omega)}$ and $|v|_{W^{m,p}(\Gamma)}$ are defined by replacing $|\boldalpha| \leq m$ by $|\boldalpha| = m$ in \eqref{integer_sobolev_norm_bulk}-\eqref{integer_sobolev_norm_surf}. The bulk- and surface Sobolev function spaces are defined as
\begin{align*}
&W^{m,p}(\Omega) := \{u:\Omega\rightarrow\mathbb{R}\ |\ \|u\|_{W^{m,p}(\Omega)} < +\infty\};\\
&W^{m,p}(\Gamma) := \{v:\Gamma\rightarrow\mathbb{R}\ |\ \|v\|_{W^{m,p}(\Gamma)} < +\infty\},
\end{align*}
respectively. For $p=2$ we will adopt the standard notations $H^m(\Omega) := W^{m,2}(\Omega)$ and $H^m(\Gamma) := W^{m,2}(\Gamma)$.
\end{definition}

The following definition can be found in \cite{Di_Nezza_2012}.
\begin{definition}[Fractional Sobolev function spaces]
Let $s \in (0,1)$ and $p\in[1,+\infty)$. For $u:\Omega \rightarrow \mathbb{R}$, the bulk- fractional Sobolev norm is defined by
\begin{equation}
\|u\|_{W^{s,p}(\Omega)} := \left(\intomega |u(\boldx)|^p\mathrm{d}\boldx + \intomega\intomega \frac{|u(\boldx)-u(\boldy)|^p}{\|\boldx-\boldy\|^{2+sp}}\mathrm{d}\boldx\mathrm{d}\boldy\right)^{1/p},
\end{equation}
with $\|\cdot\|$ being the Euclidean norm. If $s>1$, $s\notin\mathbb{N}$, by decomposing $s$ as $s = m + \sigma$, where $m\in\mathbb{N}$ and $\sigma \in (0,1)$, the fractional Sobolev norm is defined as
\begin{equation}
\label{fractional_sobolev_norm_bulk}
\|u\|_{W^{s,p}(\Omega)} := \left(\|u\|_{W^{m,p}(\Omega)}^p + \sum_{|\boldalpha| = m}\|D_{\boldalpha}u\|_{W^{\sigma,p}(\Omega)}^p \right)^{1/p}.
\end{equation}
For any (integer or non-integer) $s\in [0,+\infty)$, the \emph{Sobolev-Slobodeckij} space $W^{s,p}(\Omega)$ is defined as
\begin{equation}
W^{s,p}(\Omega) := \left\{u:\Omega\rightarrow\mathbb{R}\ |\ \|u\|_{W^{s,p}(\Omega)} < +\infty\right\},
\end{equation}
where the norm $\|\cdot\|_{W^{s,p}(\Omega)}$ is understood as \eqref{integer_sobolev_norm_bulk} or \eqref{fractional_sobolev_norm_bulk} according to whether $s$ is an integer or a non-integer number.
\end{definition}

\begin{lemma}[Inclusion between fractional Sobolev spaces]
Let $\Omega \subset\mathbb{R}^2$ be a bounded domain with a $\mathcal{C}^1$ boundary $\Gamma$, let $p \in [1,+\infty)$ and $s,s' \in [0,+\infty)$ such that $s < s'$. Then there exists a constant $C>0$ depending on $\Omega$ and $s$ such that
\begin{equation}
\|u\|_{W^{s,p}(\Omega)} \leq C \|u\|_{W^{s',p}(\Omega)},
\end{equation}
for all $u\in W^{s',p}(\Omega)$. Hence, $W^{s,p}(\Omega) \subset W^{s',p}(\Omega)$.
\begin{proof}
See \cite{Di_Nezza_2012}.
\end{proof}
\end{lemma}

\subsection*{A4. Fundamental results in bulk- and surface calculus}
\begin{theorem}[Co-area formula]
\label{thm:co_area_formula}
Let $\delta >0$ such that the \emph{narrow band} $U_\delta := \{\boldx\in \Omega | -\delta \leq d(\boldx) \leq 0\}$ is contained in the Fermi stripe $U$ (see Fig. \ref{fig:exact_domain} for an illustration). For $-\delta \leq s \leq 0$ let $\Gamma_s$ be the \emph{parallel surface} defined by $\Gamma_s := \{\boldx\in \Omega | d(\boldx) = s\}$. For any $0< \varepsilon \leq \delta$ it holds that
\begin{equation}
\label{co_area_formula}
\int_{U_\varepsilon} f(\boldx)\mathrm{d}\boldx = \int_{-\varepsilon}^0 \mathrm{d}s\int_{\Gamma_s} f(\sigma)\mathrm{d}\sigma.
\end{equation}
\begin{proof}
See \cite{dziuk2013finite}.
\end{proof}
\end{theorem}

\begin{theorem}[Narrow band trace inequality]
With the notations of the previous theorem, there exists $C>0$ depending on $\Omega$ such that any $u \in H^1(\Omega)$ fulfils
\begin{equation}
\label{narrow_band_inequality}
\|u\|_{L^2(U_\varepsilon)} \leq C\varepsilon^{\frac{1}{2}} \|u\|_{H^1(\Omega)}.
\end{equation}
\begin{proof}
See \cite{elliottranner2013finite}.
\end{proof}
\end{theorem}

\begin{theorem}[Trace theorem and inverse trace theorem]
\label{thm:trace_theorem}
Let $k\in\mathbb{N}$, $\frac{1}{2} < s \leq k$ and assume that the boundary $\Gamma$ is a $\mathcal{C}^{k}$ curve.\footnote{The original assumption is that $\Gamma$ be a $\mathcal{C}^{k-1,1}$ curve, meaning that its derivatives up to order $k-1$ are Lipschitz continuous. For simplicity, we make the stronger assumption that $\Gamma \in\mathcal{C}^k$.} Then there exists a bounded operator $\Tr: H^s(\Omega) \rightarrow H^{s-\frac{1}{2}}(\Gamma)$, called the \emph{trace operator}, such that $\Tr(u) = u_{|\Gamma}$. The trace operator fulfils
\begin{equation}
\label{trace_inequality}
\|\Tr(u)\|_{H^{s-\frac{1}{2}}(\Gamma)} \leq C\|u\|_{H^s(\Omega)}, \qquad \forall\ u\in H^s(\Omega).
\end{equation}
The trace operator has a continuous inverse operator $\Tr^{-1}:H^{s-\frac{1}{2}}(\Gamma) \rightarrow H^s(\Omega)$ called \emph{Babi\v{c} inverse} such that
\begin{equation}
\|\Tr^{-1}(v)\|_{H^s(\Omega)} \leq C\|v\|_{H^{s-\frac{1}{2}}(\Gamma)}, \qquad \forall\ v\in H^{s-\frac{1}{2}}(\Gamma).
\end{equation}
\begin{proof}
See \cite{sobolev1964} or \cite{Stein_1971}.
\end{proof}
\end{theorem}

A simple consequence of Theorem \ref{thm:trace_theorem} is the following
\begin{corollary}[Normal trace theorem]
\label{thm:normal_trace_theorem}
Let $k\in\mathbb{N}\setminus\{1\}$, $\frac{3}{2} < s \leq k$ and assume that the boundary $\Gamma$ is a $\mathcal{C}^{k}$ curve. There exists a bounded operator $Tn: H^s(\Omega) \rightarrow L^2(\Gamma)$, called the \emph{normal trace operator}, such that $Tn(u) = \nabla u \cdot \boldn$, with $\boldn$ being the outward unit normal vector field on $\Gamma$. The normal trace fulfils
\begin{equation}
\label{normal_trace_inequality}
\|Tn(u)\|_{L^2(\Gamma)} \leq C\|u\|_{H^s(\Omega)}.
\end{equation}
\begin{proof}
Let $u\in H^s(\Omega)$, then $\nabla u \in H^{s-1}(\Omega)$. From Theorem \ref{thm:trace_theorem} we have that $\Tr(\nabla u) \in H^{s-3/2}(\Gamma)$.
If $\|\cdot\|$ denotes the Euclidean norm in $\mathbb{R}^2$, we have that
\begin{equation*}
\begin{split}
\|Tn(u)\|_{L^2(\Gamma)}^2 &= \intgamma (\Tr(\nabla u) \cdot \boldn)^2 \leq \intgamma \|\Tr(\nabla u)\|^2 = \|\Tr(\nabla u)\|_{L^2(\Gamma)}^2 \leq \|\Tr(\nabla u)\|_{H^{s-3/2}(\Gamma)}^2\\
&\leq C\|\nabla u\|_{H^{s-1}(\Omega)}^2 \leq C\|u\|_{H^s(\Omega)}.
\end{split}
\end{equation*}
\end{proof} 
\end{corollary}

\begin{theorem}[Sobolev embeddings]
\label{thm:sobolev_embedding}
Assume that $\Omega$ has a Lipschitz boundary.
\begin{itemize}
\item If $0 < \gamma < 1$, then $H^{1+\gamma}(\Omega) \hookrightarrow \mathcal{C}^{0,\gamma}(\Omega)$ is a continuous embedding, hence $\|u\|_{\mathcal{C}^{0,\gamma}(\Omega)} \leq C_\gamma\|u\|_{H^{1+\gamma}(\Omega)}$. From the definition of the H\"{o}lder space $\mathcal{C}^{0,\gamma}(\Omega)$ we have that
\begin{equation}
\label{sobolev_inequality_holder}
\|u(\boldx) - u(\boldy)\| \leq C_\gamma\|u\|_{H^{1+\gamma}(\Omega)}\|\boldx - \boldy\|^\gamma, \qquad \text{a.e.} \ (\boldx,\boldy) \in \Omega^2.
\end{equation}
\item If $\varepsilon > 0$, then $H^{1+\varepsilon}(\Omega) \hookrightarrow \mathcal{C}(\Omega)$ is a continuous embedding.
\end{itemize}
\begin{proof}
See \cite{adams2003sobolev} for the case of integer-order Sobolev spaces and \cite{Di_Nezza_2012} for the fractional case.
\end{proof}
\end{theorem}

\section*{Appendix B: proofs of main theorems}
\setcounter{equation}{0}
\renewcommand{\theequation}{B.\arabic{equation}}

\subsection*{Proof of Theorem \ref{thm:ritz_proj_L2}}
By assumption we have that $\Tr(u) \in H^2(\Gamma)$. Since $\Gamma$ is a $\mathcal{C}^3$ curve, then for all $F\in\mathcal{F}_h$, the normal projection $\bolda:F \rightarrow \bolda(F)$ is a $\mathcal{C}^2$ function, see Remark \ref{rmk:regularity_normal_projection}. Hence, $(u^{-\ell})_{|F} \in H^2(F)$ and
\begin{equation}
\label{aubin_nitsche_boundary_1}
\sum_{F\in\mathcal{F}_h} \|u^{-\ell}\|_{H^2(F)} \leq C\|u\|_{H^2(\Gamma)}.
\end{equation}
From \eqref{equivalence_gamma_l2}, \eqref{interpolation_error_surface}, the second equation in \eqref{ritz_projection_definition}, and \eqref{aubin_nitsche_boundary_1} we have that
\begin{align}
\label{aubin_nitsche_boundary_L2}
&\|u-(\mathcal{R}u)^\ell\|_{L^2(\Gamma)} = \|u - I_\Gamma(u^{-\ell})^\ell\|_{L^2(\Gamma)} \leq C\|u^{-\ell} - I_\Gamma(u^{-\ell})\|_{L^2(\Gamma_h)} \leq Ch^2\|u\|_{H^2(\Gamma)}\\
\label{aubin_nitsche_boundary_H1}
&|u-(\mathcal{R}u)^\ell|_{H^1(\Gamma)} = |u - I_\Gamma(u^{-\ell})^\ell|_{H^1(\Gamma)} \leq C|u^{-\ell} - I_\Gamma(u^{-\ell})|_{H^1(\Gamma_h)} \leq Ch\|u\|_{H^2(\Gamma)}
\end{align}
From Lemma \ref{lmm:L2_preserving_inverse_trace} there exists $e_B \in H^1(\Omega)$ such that $\Tr(e_B) = \Tr(u-(\mathcal{R}u)^\ell)$ and
\begin{align}
\label{eb_L2}
&\|e_B\|_{L^2(\Omega)} \leq \|\Tr(u-(\mathcal{R}u)^\ell)\|_{L^2(\Gamma)} \leq Ch^{2}\|u\|_{H^2(\Gamma)},\\
\label{eb_H1}
&|e_B|_{H^1(\Omega)} \leq \|\Tr(u-(\mathcal{R}u)^\ell)\|_{H^1(\Gamma)} \leq Ch\|u\|_{H^2(\Gamma)}.
\end{align}
We will use an adapted Aubin-Nitsche duality method. Consider the variational problem: find $\eta \in H^1_0(\Omega)$ such that
\begin{equation}
\label{aubin_nitsche_problem}
\intomega \nabla \eta \cdot \nabla \varphi = \intomega (u-(\mathcal{R}u)^\ell-e_B) \varphi, \qquad \forall\ \varphi \in H^1_0(\Omega).
\end{equation}
Since $u-(\mathcal{R}u)^\ell-e_B \in H^1_0(\Omega)$, by elliptic regularity, the variational problem \eqref{aubin_nitsche_problem} has a unique solution $\eta \in H^3_0(\Omega)$ that fulfils
\begin{align}
\label{elliptic_regularity_1}
&\|\eta\|_{H^2(\Omega)} \leq C\|(u-(\mathcal{R}u)^\ell-e_B)\|_{L^2(\Omega)};\\
&\label{elliptic_regularity_2}
\|\eta\|_{H^3(\Omega)} \leq C\|(u-(\mathcal{R}u)^\ell-e_B)\|_{H^1(\Omega)} \leq C|(u-(\mathcal{R}u)^\ell-e_B)|_{H^1(\Omega)},
\end{align}
thanks to Poincaré's inequality. By combining \eqref{bulk_ritz_estimate_bulk_H1}, \eqref{eb_H1} and \eqref{elliptic_regularity_2} we have that
\begin{equation}
\label{elliptic_regularity_4}
\|\eta\|_{H^3(\Omega)} \leq Ch\|u\|_{H^2(\Omega)} + Ch\|u\|_{H^{2+1/4}(\Omega)}.
\end{equation}
Since $u - (\mathcal{R}u)^\ell - e_B \in H^1_0(\Omega)$ we can choose $\varphi = u - (\mathcal{R}u)^\ell - u_B$ in \eqref{aubin_nitsche_problem} and we get
\begin{equation}
\label{L2_bound_ritz_1}
\begin{split}
\|u - (\mathcal{R} u)^\ell - e_B\|_{L^2(\Omega)}^2 &= \intomega \varphi^2 = \intomega \nabla \eta \cdot \nabla \varphi = \intomega \nabla \eta\cdot \nabla \left(u - (\mathcal{R}u)^\ell - e_B\right)\\
&=  - \intomega \nabla \eta\cdot \nabla e_B + \intomega \nabla \eta\cdot \nabla \left(u - (\mathcal{R}u)^\ell\right)
\end{split}
\end{equation}
We estimate the first term on the right-hand-side of \eqref{L2_bound_ritz_1} by using Green's identity, \eqref{normal_trace_inequality}, \eqref{aubin_nitsche_boundary_L2}, \eqref{eb_L2} and \eqref{elliptic_regularity_1}:
\begin{equation}
\label{L2_bound_ritz_1.5}
\begin{split}
&\left|\intomega \nabla \eta\cdot \nabla e_B\right| \leq \left|\intomega e_b \Delta \eta\right| + \left|\intgamma e_b \nabla\eta\cdot\boldn\right| \leq \|e_B\|_{L^2(\Omega)}\|\eta\|_{H^2(\Omega)} + \|e_B\|_{L^2(\Gamma)}\|Tn(\eta)\|_{L^2(\Gamma)}\\
\leq &\|e_B\|_{L^2(\Omega)}\|\eta\|_{H^2(\Omega)} + \|e_B\|_{L^2(\Gamma)}\|\eta\|_{H^2(\Omega)} \leq Ch^2\|(u-(\mathcal{R}u)^\ell-e_B)\|_{L^2(\Omega)}.
\end{split}
\end{equation}
By using \eqref{lifting_error_bulk_nabla}, \eqref{ritz_projection_definition} and the triangle inequality we estimate the the second term of the right-hand-side of \eqref{L2_bound_ritz_1} as follows
\begin{equation}
\label{L2_bound_ritz_2}
\begin{split}
&\intomega \nabla \eta\cdot \nabla \left(u - (\mathcal{R}u)^\ell\right)\\
= &\intomega \nabla \left(u - (\mathcal{R}u)^\ell\right)\cdot \nabla\left(\eta-I_\Omega(\tilde{\eta})^\ell\right) - \intomega \nabla (\mathcal{R}u)^\ell \cdot \nabla I_\Omega(\tilde{\eta})^\ell + a_h(\mathcal{R}u, I_\Omega(\tilde{\eta}))\\
\leq &|u - (\mathcal{R}u)^\ell|_{H^1(\Omega)}|\eta-I_\Omega(\tilde{\eta})^\ell|_{H^1(\Omega)} - \intomega \nabla (\mathcal{R}u)^\ell \cdot \nabla I_\Omega(\tilde{\eta})^\ell + a_h(\mathcal{R}u, I_\Omega(\tilde{\eta}))\\
\leq &|u - (\mathcal{R}u)^\ell|_{H^1(\Omega)}|\eta-I_\Omega(\tilde{\eta})^\ell|_{H^1(\Omega)} + Ch|(\mathcal{R}u)^\ell|_{H^1(\Omega_B^\ell)}|I_\Omega(\tilde{\eta})|_{H^1(\Omega_B^\ell)}\\
-&\intomegah \nabla \mathcal{R}u \cdot \nabla I_\Omega(\tilde{\eta}) + a_h(\mathcal{R}u, I_\Omega(\tilde{\eta}))\\
\leq &C\left(|u - (\mathcal{R}u)^\ell|_{H^1(\Omega)}  + h^{1/2}|u|_{H^1(\Omega_B^\ell)}\right)\left(|\eta-I_\Omega(\tilde{\eta})^\ell|_{H^1(\Omega)} + h^{1/2}|\eta|_{H^1(\Omega_B^\ell)}\right)\\
-&\intomegah \nabla \mathcal{R}u\cdot \nabla I_\Omega(\tilde{\eta}) + a_h(\mathcal{R}u, I_\Omega(\tilde{\eta})),\\
\end{split}
\end{equation}
where we have used $h < h_0$ in the last inequality. We are left to estimate the right-hand-side of \eqref{L2_bound_ritz_2} piecewise. First, from \eqref{narrow_band_inequality} and \eqref{bulk_ritz_estimate_bulk_H1} we have that
\begin{equation}
\label{L2_bound_ritz_3}
\begin{split}
|u - (\mathcal{R}u)^\ell|_{H^1(\Omega)}  + h^{1/2}|u|_{H^1(\Omega_B^\ell)} \leq &Ch\|u\|_{H^2(\Omega)} + Ch\|u\|_{H^{2+1/4}(\Omega)}.
\end{split}
\end{equation}
Moreover, from \eqref{sobolev_extension}, \eqref{narrow_band_inequality}, \eqref{error_between_lift_and_extension}, \eqref{interpolation_error_bulk}, \eqref{elliptic_regularity_1} and \eqref{elliptic_regularity_4} we have that
\begin{equation}
\label{L2_bound_ritz_4}
\begin{split}
&|\eta-I_\Gamma(\tilde{\eta})^\ell|_{H^1(\Omega)} + h^{1/2}|\eta|_{H^1(\Omega_B^\ell)} \leq C|\eta^{-\ell}-I_\Gamma(\tilde{\eta})|_{H^1(\Omega_h)} + Ch\|\eta\|_{H^2(\Omega)}\\
\leq &C|\eta^{-\ell}-\tilde{\eta}|_{H^1(\Omega_h)} + C|\tilde{\eta}-I_\Gamma(\tilde{\eta})|_{H^1(\Omega_h)} + Ch\|\eta\|_{H^2(\Omega)}\\
\leq &Ch^{2}\|\eta\|_{H^3(\Omega)} + Ch\|\tilde{\eta}\|_{H^2(\Omega_h)} + Ch\|\eta\|_{H^2(\Omega)} \leq Ch\|\eta\|_{H^2(\Omega)} + Ch^{2}\|\eta\|_{H^3(\Omega)}\\
\leq &Ch\|u-(\mathcal{R}u)^\ell - e_B\|_{L^2(\Omega)} + Ch^3\|u\|_{H^2(\Omega)} + Ch^3\|u\|_{H^{2+1/4}(\Omega)}.
\end{split}
\end{equation}
Finally, we estimate the last two terms in \eqref{L2_bound_ritz_2} by adapting the approach used in \cite[Lemma 3.1]{vacca2015virtual}: from \eqref{error_between_lift_and_extension}, \eqref{projection_error}, \eqref{consistency}, \eqref{interpolation_error_bulk} and \eqref{elliptic_regularity_1} we have
\begin{equation}
\label{L2_bound_ritz_5}
\begin{split}
&a_h(\mathcal{R}u, I_\Gamma(\tilde{\eta})) - \intomegah \nabla \mathcal{R}u \cdot \nabla I_\Gamma(\tilde{\eta})\\
= &\sum_{E\in\mathcal{E}_h} \left(\int_E \nabla (\mathcal{R}u - \tilde{u}_\pi)\cdot\nabla(I_\Gamma(\tilde{\eta}) - \tilde{\eta}_\pi) - a_E(\mathcal{R}u - \tilde{u}_\pi, I_\Gamma(\tilde{\eta}) - \tilde{\eta}_\pi)\right)\\
\leq &\sum_{E\in\mathcal{E}_h} |\mathcal{R}u - \tilde{u}_\pi|_{H^1(E)}|I_\Gamma(\tilde{\eta}) - \tilde{\eta}_\pi|_{H^1(E)} \leq \left(Ch\|u\|_{H^2(\Omega)} + Ch\|u\|_{H^{2+1/4}(\Omega)}\right)Ch\|\eta\|_{H^2(\Omega)}\\
= &Ch^2\left(\|u\|_{H^2(\Omega)} + \|u\|_{H^{2+1/4}(\Omega)}\right)\|u-(\mathcal{R}u)^\ell - e_B\|_{L^2(\Omega)}.
\end{split}
\end{equation}
By combining \eqref{L2_bound_ritz_1}-\eqref{L2_bound_ritz_5} and using \eqref{eb_L2} we get
\begin{equation}
\begin{split}
\|u-(\mathcal{R}u)^\ell\|_{L^2(\Omega)}^2 \leq &Ch^2\left(\|u\|_{H^2(\Omega)} + \|u\|_{H^{2+1/4}(\Omega)}\right) \|u-(\mathcal{R}u)^\ell\|_{L^2(\Omega)}\\
+ &Ch^4\left(\|u\|_{H^2(\Omega)} + \|u\|_{H^{2+1/4}(\Omega)}\right)^2,
\end{split}
\end{equation}
with the terms in $H^{2+1/4}(\Omega)$ norm arising only in the simultaneous presence of curved boundaries and non-triangular boundary elements, which proves \eqref{bulk_ritz_estimate_bulk_L2}. \qed

\subsection*{Proof of Theorem \ref{thm:convergence}}
We split the error as
\begin{equation}
\label{split_error}
(u,v)-(U^\ell, V^\ell) = ((u,v)-(\mathcal{R}u, \mathcal{R}v)^\ell) + ((\mathcal{R}u, \mathcal{R}v) - (U,V))^\ell =: (\rho_u, \rho_v) + (\theta_u, \theta_v).
\end{equation}
From \eqref{surface_ritz_estimate_L2_and_H1}, \eqref{bulk_ritz_estimate_bulk_H1}, \eqref{bulk_ritz_estimate_bulk_L2} and \eqref{bulk_ritz_estimate_boundary_L2_and_H1}, since the Ritz projections swap with time derivatives, we have that
\begin{align}
\label{rho_estimate_surf}
&\|(\Tr(\rho_u), \Tr(\rho_{u,t}), \rho_v, \rho_{v,t})\|_{L^2(\Gamma)} \leq Ch^2\|(\Tr(u), \Tr(u_t), v, v_t)\|_{H^2(\Gamma)};\\
\label{rho_estimate_bulk}
&\|(\rho_u, \rho_{u,t})\|_{L^2(\Omega)} \leq Ch^{2}\left(\|(u,u_t)\|_{H^2(\Omega)} + \|(u,u_t)\|_{H^{2+1/4}(\Omega)}\right).
\end{align}
We are left to estimate the the norms of $\theta_u$ and $\theta_v$. By using \eqref{parabolic_problem}, \eqref{parabolic_problem_BSVEM_projection}, \eqref{ritz_projection_definition} and \eqref{split_error} we have the following error equation
\begin{equation}
\label{error_equation}
\begin{split}
& m_h\left(\frac{\mathrm{d}}{\mathrm{d}t}\theta_u^{-\ell}, \varphi\right) + d_u a_h\left(\theta_u^{-\ell}, \varphi\right) + \intgammah \frac{\mathrm{d}}{\mathrm{d}t}\theta_v^{-\ell}\psi + d_v \intgammah \nablagammah \theta_v^{-\ell} \cdot \nablagammah \psi\\
= &m_h\left(\frac{\mathrm{d}}{\mathrm{d}t}U, \varphi\right) - m_h\left(\frac{\mathrm{d}}{\mathrm{d}t}\mathcal{R}u, \varphi\right) + d_u a_h(U,\varphi) - d_u a_h(\mathcal{R}u, \varphi)\\
&+\intgammah \frac{\mathrm{d}}{\mathrm{d}t} V \psi - \intgammah \frac{\mathrm{d}}{\mathrm{d}t}\mathcal{R}v\psi + d_v \intgammah\nablagammah V\cdot\nablagammah\psi - d_v \intgammah\nablagammah\mathcal{R}v\cdot\nablagammah\psi\\
= &\ m_h(I_\Omega q(\Pi^0 U), \varphi) - m_h\left(\frac{\mathrm{d}}{\mathrm{d}t}\mathcal{R}u, \varphi\right) - d_u a_h(\mathcal{R}u, \varphi) + \intgammah I_\Gamma (s(U,V))\varphi_{|\Gamma}\\
& +\intgammah I_\Gamma(r(U,V))\psi - \intgammah \frac{\mathrm{d}}{\mathrm{d}t}\mathcal{R}v\psi - d_v \intgammah\nablagammah\mathcal{R}v\cdot\nablagammah\psi - \intgammah I_\Gamma (s(U,V))\psi\\
= &\underbrace{\intomega \frac{\mathrm{d}}{\mathrm{d}t} u\varphi^\ell - m_h\left(\frac{\mathrm{d}}{\mathrm{d}t}\mathcal{R}u, \varphi\right)}_{T_1}\ + \underbrace{m_h(\Pi^0 I_\Omega q(\Pi^0 U), \varphi)  - \intomega q(u)\varphi^\ell}_{T_2}\\
& +\underbrace{\intgamma \frac{\mathrm{d}}{\mathrm{d}t}v\psi^\ell - \intgammah \frac{\mathrm{d}}{\mathrm{d}t}\mathcal{R}v\psi}_{T_3}\ + \underbrace{\intgammah I_\Gamma(r(U, V))\psi - \intgamma r(u, v)\psi^\ell}_{T_4}\\
& + \underbrace{\intgammah I_\Gamma (s(U,V))(\Tr(\varphi) - \psi) - \intgamma (s(u,v))(\Tr(\varphi) - \psi)^\ell}_{T_5}.
\end{split}
\end{equation}
We now analyse terms $T_1-T_5$ on the right hand side of \eqref{error_equation} separately. For $T_1$ we use \eqref{narrow_band_inequality}, \eqref{lifting_error_bulk_mass} and \eqref{rho_estimate_bulk}:
\begin{equation}
\label{T1_0}
\begin{split}
T_1 = &\intomega (u_t - (\mathcal{R}u_t)^\ell)\varphi^\ell + \intomega (\mathcal{R}u_t)^\ell\varphi^\ell - \intomegah \mathcal{R}u_t\varphi+\intomegah \mathcal{R}u_t\varphi - m_h(\mathcal{R}u_t, \varphi)\\
\leq &\|u_t - (\mathcal{R}u_t)^\ell\|_{L^2(\Omega)}\|\varphi^\ell\|_{L^2(\Omega)} + Ch\|(\mathcal{R}u_t)^\ell\|_{L^2(\Omega_B^\ell)}\|\varphi^\ell\|_{L^2(\Omega_h^\ell)} + \intomegah \mathcal{R}u_t\varphi - m_h(\mathcal{R}u_t, \varphi)\\
\leq &\|u_t - (\mathcal{R}u_t)^\ell\|_{L^2(\Omega)}\|\varphi^\ell\|_{L^2(\Omega)} + Ch^2\|(\mathcal{R}u_t)^\ell\|_{H^1(\Omega)}\|\varphi^\ell\|_{H^1(\Omega)} + \intomegah \mathcal{R}u_t\varphi - m_h(\mathcal{R}u_t, \varphi)\\
\leq &Ch^{2}\left(\|u_t\|_{H^2(\Omega)} + \|u_t\|_{H^{2+1/4}(\Omega)}\right)\|\varphi^\ell\|_{H^1(\Omega)} + \intomegah \mathcal{R}u_t\varphi - m_h(\mathcal{R}u_t, \varphi).
\end{split}
\end{equation}
We estimate the last term in \eqref{T1_0} by using \eqref{sobolev_extension}, \eqref{narrow_band_inequality}, \eqref{equivalence_omega_l2}, \eqref{error_between_lift_and_extension_l2}, \eqref{projection_error}, \eqref{consistency} and \eqref{rho_estimate_bulk}:
\begin{equation}
\label{T1_1}
\begin{split}
&\intomegah \mathcal{R}u_t\varphi - m_h(\mathcal{R}u_t, \varphi)\\
= &\sum_{E\in\mathcal{E}_h} \left(\int_E (\mathcal{R}u_t - \tilde{u}_{t,\pi})\varphi - m_E(\mathcal{R}u_t - \tilde{u}_{t,\pi}, \varphi)\right) \leq \sum_{E\in\mathcal{E}_h} \|\mathcal{R}u_t - \tilde{u}_{t,\pi}\|_{L^2(E)}\|\varphi\|_{L^2(E)}\\
\leq &C\left(\|\rho_{u,t}\|_{L^2(\Omega)} + \|u_t^{-\ell} - \tilde{u}_t\|_{L^2(\Omega_h)} + \sum_{E\in\mathcal{E}_h}\|\tilde{u}_t - \tilde{u}_{t,\pi}\|_{L^2(E)}  \right)\|\varphi^\ell\|_{L^2(\Omega)}\\
\leq &Ch^2\left(\|u_t\|_{H^2(\Omega)} + \|u_t\|_{H^{2+1/4}(\Omega)}\right)\|\varphi^\ell\|_{L^2(\Omega)}.
\end{split}
\end{equation}
By combining \eqref{T1_0} and \eqref{T1_1} we get
\begin{equation}
\label{T1}
T_1 \leq Ch^2\left(\|u_t\|_{H^2(\Omega)} + \|u_t\|_{H^{2+1/4}(\Omega)}\right)\|\varphi^\ell\|_{H^1(\Omega)}.
\end{equation}
We estimate the second term by adapting the approach used in \cite[Theorem 4.2]{adak2019convergence}. From \eqref{narrow_band_inequality}, \eqref{lifting_error_bulk_mass}, \eqref{error_between_lift_and_extension_l2}, \eqref{projection_error}, \eqref{L2_projector_boundedness}, \eqref{interpolation_error_bulk}-\eqref{interpolant_bulk_quasi_bounded} and \eqref{rho_estimate_bulk}, since $q$ is $\mathcal{C}^2$ and Lipschitz continuous, we have that
\begin{equation}
\label{T2}
\begin{split}
T_2 =\ &\intomegah \Pi^0 I_\Omega (q(\Pi^0 U) - q(\tilde{u}))\varphi + \intomegah \Pi^0 (I_\Omega q(\tilde{u}) - q(\tilde{u}))\varphi\\
 &+ \intomegah (\Pi^0 q(\tilde{u}) - q(\tilde{u}))\varphi + \intomegah (q(\tilde{u}) - q(u^{-\ell}))\varphi + \intomegah q(u^{-\ell})\varphi - \intomega q(u)\varphi^\ell\\
 \leq &C\|q(\Pi^0 U) - q(\tilde{u})\|_{L^2(\Omega_h)}\|\varphi\|_{L^2(\Omega_h)} + Ch^2|q(\Pi^0 U) - q(\tilde{u})|_{2,\Omega,h}\|\varphi\|_{L^2(\Omega_h)}\\
 & +Ch^2\|q(\tilde{u})\|_{H^2(\Omega_h)}\|\varphi\|_{L^2(\Omega_h)} + Ch^2\|u\|_{H^2(\Omega)}\|\varphi\|_{L^2(\Omega_h)} + Ch\|q(u)\|_{L^2(\Omega_B^\ell)}\|\varphi^\ell\|_{L^2(\Omega_B^\ell)}\\
 \leq &C\|U -\tilde{u}\|_{L^2(\Omega_h)}\|\varphi\|_{L^2(\Omega_h)} + Ch^2\left(|U|_{H^1(\Omega)} + |\tilde{u}|_{H^1(\Omega_h)} + |\tilde{u}|_{H^2(\Omega_h)}\right)\|\varphi\|_{L^2(\Omega_h)}\\
 & +Ch^2\left(\|u\|_{H^2(\Omega)} + \|u\|_{H^{2+1/4}(\Omega)}\right)\|\varphi^\ell\|_{L^2(\Omega)} + Ch^2\|u\|_{H^1(\Omega)}\|\varphi^\ell\|_{H^1(\Omega)}\\
\leq &C\left(\|\theta_u\|_{L^2(\Omega)} + Ch^2\|U\|_{H^1(\Omega_h)} + Ch^2\|u\|_{H^2(\Omega)} + Ch^2\|u\|_{H^{2+1/4}(\Omega)}\right)\|\varphi^\ell\|_{H^1(\Omega)}.
\end{split}
\end{equation}
For the third term, from \eqref{narrow_band_inequality}, \eqref{equivalence_gamma_l2}, \eqref{lifting_error_surf_mass}, and \eqref{rho_estimate_surf} we have that
\begin{equation}
\label{T3}
\begin{split}
T_3 = &\intgamma \frac{\mathrm{d}}{\mathrm{d}t}v\psi^\ell - \intgammah \frac{\mathrm{d}}{\mathrm{d}t}v^{-\ell}\psi + \intgammah \frac{\mathrm{d}}{\mathrm{d}t}v^{-\ell}\psi - \intgammah \frac{\mathrm{d}}{\mathrm{d}t}\mathcal{R}v\psi\\
\leq &Ch^2\|v_t\|_{L^2(\Gamma)}\|\psi^\ell\|_{L^2(\Gamma)} + C\|\rho_{v,t}\|_{L^2(\Gamma)}\|\psi^\ell\|_{L^2(\Gamma)} \leq Ch^{2}\|v_t\|_{H^2(\Gamma)}\|\psi^\ell\|_{L^2(\Gamma)}.
\end{split}
\end{equation}
To estimate the fourth term we proceed as in \eqref{T2}. By using \eqref{lifting_error_surf_mass}, \eqref{interpolation_error_surface}-\eqref{interpolation_surface_monotonic}, \eqref{aubin_nitsche_boundary_1}, \eqref{rho_estimate_surf}, the $\mathcal{C}^2$ regularity and the Lipschitz continuity of $r$, we have that
\begin{equation}
\label{T4}
\begin{split}
T_4 = &\intgammah I_\Gamma\bigg(r(U,V)-r(u^{-\ell}, v^{-\ell})\bigg)\psi + \intgammah \bigg(I_\Gamma (r(u^{-\ell}, v^{-\ell})) - r(u^{-\ell}, v^{-\ell})\bigg)\psi\\
&+ \intgammah r(u^{-\ell}, v^{-\ell})\psi - \intgamma r(u, v)\psi^\ell \leq C\|(U,V)-(I_\Gamma(u^{-\ell}), I_\Gamma(v^{-\ell}))\|_{L^2(\Gamma_h)}\|\psi\|_{L^2(\Gamma_h)}\\
& + Ch^2|r(u^{-\ell}, v^{-\ell})|_{2,\Gamma,h}\|\psi\|_{L^2(\Gamma_h)} + Ch^2\|r(u,v)\|_{L^2(\Gamma)}\|\psi\|_{L^2(\Gamma)}\\
\leq &C\big(\|(\theta_u,\theta_v)\|_{L^2(\Gamma)} + \|(\rho_u,\rho_v)\|_{L^2(\Gamma)} + h^2|(U,V)|_{H^1(\Gamma)} + h^2\|(u,v)\|_{H^2(\Gamma)}\big) \|\psi^\ell\|_{L^2(\Gamma)}\\
\leq & C \|(\theta_u,\theta_v)\|_{L^2(\Gamma)} \|\psi^\ell\|_{L^2(\Gamma)} +Ch^2\|(u,v)\|_{H^2(\Gamma)} \|\psi^\ell\|_{L^2(\Gamma)}.
\end{split}
\end{equation}
We estimate the fifth term as in \eqref{T4} by also using \eqref{trace_inequality}:
\begin{equation}
\label{T5}
\begin{split}
T_5 \leq &C \big( \|(\theta_u,\theta_v)\|_{L^2(\Gamma)} +Ch^2\|(u,v)\|_{H^2(\Gamma)} \big) \|\Tr(\varphi^\ell) - \psi^\ell\|_{L^2(\Gamma)}\\
\leq &C \big( \|(\theta_u,\theta_v)\|_{L^2(\Gamma)} +Ch^2\|(u,v)\|_{H^2(\Gamma)} \big) \big(\|\varphi^\ell\|_{H^1(\Omega)} + \|\psi^\ell\|_{L^2(\Gamma)} \big).
\end{split}
\end{equation}
By substituting \eqref{T1}-\eqref{T5} into the error equation \eqref{error_equation}, using Young's inequality and choosing $\varphi^\ell = \theta_u$, $\psi^\ell = \theta_v$ we have
\begin{equation}
\begin{split}
\frac{\mathrm{d}}{\mathrm{d}t} \|\theta_u\|_{L^2(\Omega)}^2 + d_u |\theta_u|_{H^1(\Omega)}^2 + \frac{\mathrm{d}}{\mathrm{d}t} \|\theta_v\|_{L^2(\Gamma)}^2 + d_v |\theta_v|_{H^1(\Gamma)}^2\\
\leq C(u,u_t,v,v_t, d_u)\bigg(h^4 + \|\theta_u\|_{L^2(\Omega)}^2 + \|\theta_v\|_{L^2(\Gamma)}^2\bigg) + d_u |\theta_u|_{H^1(\Omega)}^2,
\end{split}
\end{equation}
which yields
\begin{equation}
\begin{split}
\frac{\mathrm{d}}{\mathrm{d}t} \|\theta_u\|_{L^2(\Omega)}^2 + \frac{\mathrm{d}}{\mathrm{d}t} \|\theta_v\|_{L^2(\Gamma)}^2 \leq C(u,u_t,v,v_t, d_u) \bigg(h^4 + \|\theta_u\|_{L^2(\Omega)}^2 + \|\theta_v\|_{L^2(\Gamma)}^2\bigg).
\end{split}
\end{equation}
By applying Gr\"{o}nwall's lemma and accounting for the $h^2$-accuracy of the initial conditions \eqref{parabolic_problem_BSVEM_IC}, we get the desired estimate. \qed

\end{document}

%% file: exact_domain.tex
\begin{tikzpicture}[scale = 1.3, label distance=.05]

\def\a{1.5} 
\def\b{0.2} 
\def\c{0.18} 
\def\f{0.25} 

\tikzset{declare function={gx(\t)      = \a*cos(180*\t/pi);}}
\tikzset{declare function={gy(\t)      = \a*sin(180*\t/pi) + \a*\b*sin(3*180*\t/pi);}}
\tikzset{declare function={vx(\t)      = -\a*sin(180*\t/pi);}}
\tikzset{declare function={vy(\t)      = \a*cos(180*\t/pi) + \a*\b*3*cos(3*180*\t/pi);}}
\tikzset{declare function={nx(\t)      = -vy(\t);}}
\tikzset{declare function={ny(\t)      =  vx(\t);}}
\tikzset{declare function={norm(\t)  = sqrt(pow(nx(\t),2) + pow(ny(\t),2));}}
\tikzset{declare function={ggx(\t)    = gx(\t) + \a*\c*nx(\t)/norm(\t);}}
\tikzset{declare function={ggy(\t)    = gy(\t) + \a*\c*ny(\t)/norm(\t);}}
\tikzset{declare function={fox(\t)    = gx(\t) - \a*\f*nx(\t)/norm(\t);}}
\tikzset{declare function={foy(\t)    = gy(\t) - \a*\f*ny(\t)/norm(\t);}}
\tikzset{declare function={fix(\t)    = gx(\t) + \a*\f*nx(\t)/norm(\t);}}
\tikzset{declare function={fiy(\t)    = gy(\t) + \a*\f*ny(\t)/norm(\t);}}

\coordinate (CC) at (0, 0);

\draw [red, very thick,  domain=0:2*pi, samples=100, fill = lightGreen] plot ({gx(\x)}, {gy(\x)});

\fill [domain=0:2*pi, samples=200, white] plot ({ggx(\x)}, {ggy(\x)});

\fill [domain=0:2*pi, samples=200, blue, opacity=.1] plot ({fox(\x)}, {foy(\x)});

\fill [domain=0:2*pi, samples=200, white] plot ({fix(\x)}, {fiy(\x)});

\node[label=above right:{$\Omega$}] at (CC) {};
\node[label=above:{$\blue{U}$}] at ({gx(.2*pi)}, {gy(.2*pi)})  {};
\node[label=above:{$\green{U_\delta}$}] at ({gx(.5*pi)}, {gy(.5*pi)})  {};
\node[label=above:{$\red{\Gamma}$}] at ({gx(.8*pi)}, {gy(.8*pi)})  {};

\end{tikzpicture}

%% file: mesh_portions.tex
\begin{tikzpicture}[scale = 1.3]

\def\a{1.5} 
\def\b{0.2} 
\def\c{0.18} 
\def\d{2*pi/16} 
\def\f{0.25} 

\tikzset{declare function={gx(\t)      = \a*cos(180*\t/pi);}}
\tikzset{declare function={gy(\t)      = \a*sin(180*\t/pi) + \a*\b*sin(3*180*\t/pi);}}
\tikzset{declare function={vx(\t)      = -\a*sin(180*\t/pi);}}
\tikzset{declare function={vy(\t)      = \a*cos(180*\t/pi) + \a*\b*3*cos(3*180*\t/pi);}}
\tikzset{declare function={nx(\t)      = -vy(\t);}}
\tikzset{declare function={ny(\t)      =  vx(\t);}}
\tikzset{declare function={norm(\t)  = sqrt(pow(nx(\t),2) + pow(ny(\t),2));}}
\tikzset{declare function={ggx(\t)    = gx(\t) + \a*\c*nx(\t)/norm(\t);}}
\tikzset{declare function={ggy(\t)    = gy(\t) + \a*\c*ny(\t)/norm(\t);}}
\tikzset{declare function={fox(\t)    = gx(\t) - \a*\f*nx(\t)/norm(\t);}}
\tikzset{declare function={foy(\t)    = gy(\t) - \a*\f*ny(\t)/norm(\t);}}
\tikzset{declare function={fix(\t)    = gx(\t) + \a*\f*nx(\t)/norm(\t);}}
\tikzset{declare function={fiy(\t)    = gy(\t) + \a*\f*ny(\t)/norm(\t);}}

\coordinate (CC) at (0, 0);

\coordinate (A) at ({gx(0)}, {gy(0)});
\coordinate (B) at ({gx(1*\d)}, {gy(1*\d)});
\coordinate (C) at ({gx(2*\d)}, {gy(2*\d)});
\coordinate (D) at ({gx(3*\d)}, {gy(3*\d)});
\coordinate (E) at ({gx(4*\d)}, {gy(4*\d)});
\coordinate (F) at ({gx(5*\d)}, {gy(5*\d)});
\coordinate (G) at ({gx(6*\d)}, {gy(6*\d)});
\coordinate (H) at ({gx(7*\d)}, {gy(7*\d)});
\coordinate (I) at ({gx(8*\d)}, {gy(8*\d)});
\coordinate (J) at ({gx(9*\d)}, {gy(9*\d)});
\coordinate (K) at ({gx(10*\d)}, {gy(10*\d)});
\coordinate (L) at ({gx(11*\d)}, {gy(11*\d)});
\coordinate (M) at ({gx(12*\d)}, {gy(12*\d)});
\coordinate (N) at ({gx(13*\d)}, {gy(13*\d)});
\coordinate (O) at ({gx(14*\d)}, {gy(14*\d)});
\coordinate (P) at ({gx(15*\d)}, {gy(15*\d)});

\coordinate (A2) at ({ggx(0)}, {ggy(0)});
\coordinate (B2) at ({ggx(1*\d)}, {ggy(1*\d)});
\coordinate (C2) at ({ggx(2*\d)}, {ggy(2*\d)});
\coordinate (D2) at ({ggx(3*\d)}, {ggy(3*\d)});
\coordinate (E2) at ({ggx(4*\d)}, {ggy(4*\d)});
\coordinate (F2) at ({ggx(5*\d)}, {ggy(5*\d)});
\coordinate (G2) at ({ggx(6*\d)}, {ggy(6*\d)});
\coordinate (H2) at ({ggx(7*\d)}, {ggy(7*\d)});
\coordinate (I2) at ({ggx(8*\d)}, {ggy(8*\d)});
\coordinate (J2) at ({ggx(9*\d)}, {ggy(9*\d)});
\coordinate (K2) at ({ggx(10*\d)}, {ggy(10*\d)});
\coordinate (L2) at ({ggx(11*\d)}, {ggy(11*\d)});
\coordinate (M2) at ({ggx(12*\d)}, {ggy(12*\d)});
\coordinate (N2) at ({ggx(13*\d)}, {ggy(13*\d)});
\coordinate (O2) at ({ggx(14*\d)}, {ggy(14*\d)});
\coordinate (P2) at ({ggx(15*\d)}, {ggy(15*\d)});

\coordinate (A3) at ($(A2)!0.5!(CC)$);
\coordinate (B3) at ($(B2)!0.5!(CC)$);
\coordinate (C3) at ($(C2)!0.5!(CC)$);
\coordinate (D3) at ($(D2)!0.5!(CC)$);
\coordinate (E3) at ($(E2)!0.5!(CC)$);
\coordinate (F3) at ($(F2)!0.5!(CC)$);
\coordinate (G3) at ($(G2)!0.5!(CC)$);
\coordinate (H3) at ($(H2)!0.5!(CC)$);
\coordinate (I3) at ($(I2)!0.5!(CC)$);
\coordinate (J3) at ($(J2)!0.5!(CC)$);
\coordinate (K3) at ($(K2)!0.5!(CC)$);
\coordinate (L3) at ($(L2)!0.5!(CC)$);
\coordinate (M3) at ($(M2)!0.5!(CC)$);
\coordinate (N3) at ($(N2)!0.5!(CC)$);
\coordinate (O3) at ($(O2)!0.5!(CC)$);
\coordinate (P3) at ($(P2)!0.5!(CC)$);

\draw[fill=lightGreen] (A) -- (B) -- (B2) -- (A2) -- cycle;
\draw[fill=lightGreen]  (B) -- (C) -- (C2) -- (B2) -- cycle;
\draw[fill=lightGreen]  (C) -- (D) -- (D2) -- (C2) -- cycle;
\draw[fill=lightGreen]  (D) -- (E) -- (E2) -- (D2) -- cycle;
\draw[fill=lightGreen]  (E) -- (F) -- (F2) -- (E2) -- cycle;
\draw[fill=lightGreen]  (F) -- (G) -- (G2) -- (F2) -- cycle;
\draw[fill=lightGreen]  (G) -- (H) -- (H2) -- (G2) -- cycle;
\draw[fill=lightGreen]  (H) -- (I) -- (I2) -- (H2) -- cycle;
\draw[fill=lightGreen]  (I) -- (J) -- (J2) -- (I2) -- cycle;
\draw[fill=lightGreen]  (J) -- (K) -- (K2) -- (J2) -- cycle;
\draw[fill=lightGreen]  (K) -- (L) -- (L2) -- (K2) -- cycle;
\draw[fill=lightGreen]  (L) -- (M) -- (M2) -- (L2) -- cycle;
\draw[fill=lightGreen]  (M) -- (N) -- (N2) -- (M2) -- cycle;
\draw[fill=lightGreen]  (N) -- (O) -- (O2) -- (N2) -- cycle;
\draw[fill=lightGreen]  (O) -- (P) -- (P2) -- (O2) -- cycle;
\draw[fill=lightGreen]  (P) -- (A) -- (A2) -- (P2) -- cycle;

\draw[ultra thick, red] (A) -- (B) -- (C) -- (D) --  (E) -- (F) -- (G) -- (H) -- (I) -- (J) -- (K) -- (L) -- (M) -- (N) -- (O) -- (P) -- cycle;

\fill [domain=0:2*pi, samples=200, blue, opacity=.1] plot ({fox(\x)}, {foy(\x)});

\fill [domain=0:2*pi, samples=200, white] plot ({fix(\x)}, {fiy(\x)});

\draw (A2) -- (B2) -- (B3) -- (A3) -- cycle;
\draw (B2) -- (C2) -- (C3) -- (B3) -- cycle;
\draw (C2) -- (D2) -- (D3) -- (C3) -- cycle;
\draw (D2) -- (E2) -- (E3) -- (D3) -- cycle;
\draw (E2) -- (F2) -- (F3) -- (E3) -- cycle;
\draw (F2) -- (G2) -- (G3) -- (F3) -- cycle;
\draw (G2) -- (H2) -- (H3) -- (G3) -- cycle;
\draw (H2) -- (I2) -- (I3) -- (H3) -- cycle;
\draw (I2) -- (J2) -- (J3) -- (I3) -- cycle;
\draw (J2) -- (K2) -- (K3) -- (J3) -- cycle;
\draw (K2) -- (L2) -- (L3) -- (K3) -- cycle;
\draw (L2) -- (M2) -- (M3) -- (L3) -- cycle;
\draw (M2) -- (N2) -- (N3) -- (M3) -- cycle;
\draw (N2) -- (O2) -- (O3) -- (N3) -- cycle;
\draw (O2) -- (P2) -- (P3) -- (O3) -- cycle;
\draw (P2) -- (A2) -- (A3) -- (P3) -- cycle;

\draw (A3)-- (CC);
\draw (E3)-- (CC);
\draw (I3)-- (CC);
\draw (M3)-- (CC);

\node[label=above right:{$\Omega_h$}] at (CC) {};
\node[label=above:{$\blue{U}$}] at ({gx(.2*pi)}, {gy(.2*pi)})  {};
\node[label=above:{$\green{\Omega_B}$}] at ({gx(.5*pi)}, {gy(.5*pi)})  {};
\node[label=above:{$\red{\Gamma_h}$}] at ({gx(.8*pi)}, {gy(.8*pi)})  {};

\end{tikzpicture}

%% file: element.tex
\begin{tikzpicture}[scale=0.7]

// POLIGONO
\draw[fill = lightGreen] (0,0) 
  -- (2,1)
  -- (4,0) 
  -- (4,4.5)
  -- (0,4.5)
  -- cycle;
\node at (2,{4-sqrt(2)}) {\green{$E$}};  

// FERMI STRIPE
\fill [blue, opacity=.1, domain=-1:5] plot (\x, {3.9-0.5*sin(45*\x)}) -- plot [domain = 5:-1] (\x, {5.1-0.5*sin(45*\x)}) -- cycle;
\node at (4.5,4.95) {$\blue{U}$};
 
 // BOUNDARY EDGE
\draw[dotted] (3,4.5)  -- (3.5, 5.4);
\node at (3.6, 5.5) {$\bar{e}$};

// SURFACE
\draw [red, very thick,domain=-1:5] plot (\x, {4.5-0.5*sin(45*\x)});
\node at (2,5.3) {$\red{\Gamma}$};

\end{tikzpicture}

%% file: exact_element.tex
\begin{tikzpicture}[scale=0.7]

// POLIGONO
\draw[fill = lightGreen] (0,4.5)
  -- (0,0) 
  -- (2,1)
  -- (4,0) 
  -- plot [domain=4:0] (\x, {4.5-0.5*sin(45*\x)});
 \node at (2,{4-sqrt(2)}) {\green{$\breve{E}$}};

// SURFACE
\draw [red, very thick,domain=0:4] plot (\x, {4.5-0.5*sin(45*\x)});

\node at (2,5.3) {$\red{\bolda(\bar{e})}$};

\end{tikzpicture}

%% file: element_intersection.tex
\begin{tikzpicture}[scale=0.9]

// POLIGONO
\draw[fill = lightGreen] (2,3.5) 
  -- (4,4.5)
  -- (0,4.5)
  -- cycle;
\node at (2,{4-sqrt(2)}) {\green{$E$}};  

// FERMI STRIPE
\fill [blue, opacity=.1, domain=-0.3:4.3] plot (\x, {3.5-0.75*sin(90*\x)}) -- plot [domain = 4.3:-0.3] (\x, {5.5-0.75*sin(90*\x)}) -- cycle;
\node at (4.5,4.95) {$\blue{U}$};
 
 // BOUNDARY EDGE
\draw[dotted] (3,4.5)  -- (3.5, 5.4);
\node at (3.6, 5.5) {$\bar{e}$};

// SURFACE
\draw [red, very thick,domain=-0.3:4.3] plot (\x, {4.5-0.75*sin(90*\x)});
\node at (2,5.3) {$\red{\Gamma}$};

\end{tikzpicture}

%% file: element_smooth.tex
\begin{tikzpicture}[scale=0.7]

// POLIGONO
\draw[fill = lightGreen] (2,1) 
  -- (4,4.5)
  -- (0,4.5)
  -- cycle;
\node at (2,{4-sqrt(2)}) {\green{$E$}};  

// FERMI STRIPE
\fill [blue, opacity=.1, domain=-0.5:4.5] plot (\x, {3.9+0.5*sin(90*\x)}) -- plot [domain = 4.5:-0.5] (\x, {5.1+0.5*sin(90*\x)}) -- cycle;
\node at (4.5,4.5) {$\blue{U}$};
 
 // BOUNDARY EDGE
\draw[dotted] (3,4.5)  -- (3.5, 5.4);
\node at (3.6, 5.5) {$\bar{e}$};

// SURFACE
\draw [red, very thick,domain=-0.5:4.5] plot (\x, {4.5+0.5*sin(90*\x)});
\node at (2,5.3) {$\red{\Gamma}$};

\end{tikzpicture}

%% file: element_smooth_exact.tex
\begin{tikzpicture}[scale=0.7]

// POLIGONO
\draw[fill = lightGreen, draw=none] (2,1) 
  -- (4,4.5)
  -- (0,4.5)
  -- cycle;
 \node at (2,{4-sqrt(2)}) {\green{$\breve{E}$}};  

// BORDO NERO DEL POLIGONO
\draw (4,4.5)
  -- (2,1) 
  -- (0,4.5);

// SURFACE
\fill [color=lightGreen, draw=lightGreen, domain=0:2] plot (\x, {4.5+0.5*sin(90*\x)}) -- cycle;
\fill [color=white, draw=white, domain=2:4] plot (\x, {4.5+0.5*sin(90*\x)}) -- cycle;
\draw [red, very thick,domain=0:4] plot (\x, {4.5+0.5*sin(90*\x)});

\node at (2,5.3) {$\red{\bolda(\bar{e})}$};

\end{tikzpicture}

%% file: element_corner.tex
\begin{tikzpicture}[scale=0.7]

// POLIGONO
\draw[fill = lightGreen] (4,1)
  -- (0,1)
  -- (2,5)
  -- cycle;
\node at (2,{4-sqrt(2)}) {\green{$E$}};  

// FERMI STRIPE
\fill [blue, opacity=.1, domain=0:4] plot (\x, {max(4*(sin(45*\x))^(3),1)}) -- plot [domain=4:0] (\x, {2+4*(sin(45*\x))^(1/3)}) -- cycle;
\node at (4,4.5) {$\blue{U}$};
 
// BOUNDARY EDGE e1
\draw[dotted] (1.4,3.8)  -- (0.5, 5.4);
\node at (0.4, 5.5) {$\bar{e}_1$};

// BOUNDARY EDGE e2
\draw[dotted] (2.6,3.8)  -- (3.5, 5.4);
\node at (3.6, 5.5) {$\bar{e}_2$};

// SURFACE
\draw [red, very thick,domain=0:4] plot (\x, {1+4*sin(45*\x)});
\node at (2,5.3) {$\red{\Gamma}$};

\end{tikzpicture}

%% file: element_corner_exact.tex
\begin{tikzpicture}[scale=0.7]

// POLIGONO
\draw[fill = lightGreen] (4,1)
  -- (0,1)
  -- plot [domain=0:4] (\x, {1+4*sin(45*\x)})
  -- cycle;
\node at (2,{4-sqrt(2)}) {\green{$E$}};  
 
// CURVED EDGE e1
\draw[dotted, color=blue] (1,3.6)  -- (0.5, 5.4);
\node at (0.4, 5.5) {$\blue{\bolda(\bar{e}_1)}$};
\draw [blue, very thick,domain=0:2] plot (\x, {1+4*sin(45*\x)});
 
// CURVED EDGE e2
\draw[dotted, color=red] (3,3.6)  -- (3.5, 5.4);
\node at (3.6, 5.5) {$\red{\bolda(\bar{e}_2)}$};
\draw [red, very thick,domain=2:4] plot (\x, {1+4*sin(45*\x)});

// SURFACE

\end{tikzpicture}

%% file: lemma_proof_1.tex
\begin{tikzpicture}[scale = 0.65]

// POLIGONO
\draw[fill = lightGreen] (0,0) 
  -- (2,1)
  -- (4,0) 
  -- (4,4.5)
  -- (0,4.5)
  -- cycle;
 \node at (3.7,{4-sqrt(2)}) {\green{$E$}};
  
 \draw[dotted] (0,2) -- (2,1) -- (4,2);
  
 // STAR-SHAPED-BALL
 \draw [black, fill = lightBlue] (2,{4-sqrt(2)}) circle [radius={sqrt(2)}];
 \node at (2.6, 3.3) {\blue{$\mathcal{B}_E$}};
 
 // RADIUS OF BALL
 \draw (2,{4-sqrt(2)}) -- ({2-sqrt(2)},{4-sqrt(2)});
 \draw[dotted] ({2-sqrt(2)/2}, {4-sqrt(2)}) -- (2, 5.2);
 \node at (2.1, 5.3) {$R_E$};
 
 // CENTER OF BALL
 \fill (2,{4-sqrt(2)}) circle [radius = 2pt];
 \node at (2.6,{4-sqrt(2)}) {$\boldx_E$};

\end{tikzpicture}

%% file: lemma_proof_2.tex
\begin{tikzpicture}[scale=0.65]

// POLIGONO
\draw[fill = lightGreen] (0,0) 
  -- (2,1)
  -- (4,0) 
  -- (4,4.5)
  -- (0,4.5)
  -- cycle;
 \node at (3.7,{4-sqrt(2)}) {\green{$E$}};  
 
 // TRIANGLE
 \draw[fill = lightRed] (0,4.5) -- (2,{4-sqrt(2)}) -- (4,4.5) -- cycle;
 \node at (2, 3.6) {\red{$T_{\bar{e}}$}};
 
 // BOUNDARY EDGE
 \draw[dotted] (3,4.5)  -- (3.5, 5.2);
 \node at (3.6, 5.3) {$\bar{e}$};
 
 // CENTER OF BALL
 \fill (2,{4-sqrt(2)}) circle [radius = 2pt];
 \node at (2.6,{4-sqrt(2)}) {$\boldx_E$};

// SURFACE
\draw [red,thick,domain=-1:5] plot (\x, {4.5+0.3*sin(90*\x)});
\node at (2,5.3) {$\red{\Gamma}$};

\end{tikzpicture}

%% file: lemma_proof_3.tex
\begin{tikzpicture}[scale=0.65]

// POLIGONO
\draw[fill = lightGreen] (0,0) 
  -- (2,1)
  -- (4,0) 
  -- (4,4.5)
  -- (0,4.5)
  -- cycle;
 \node at (3.7,{4-sqrt(2)}) {\green{$E$}};
  
 \draw[dotted] (0,2) -- (2,1) -- (4,2);
  
 // STAR-SHAPED-BALL
 \fill [lightBlue] (2,{4-sqrt(2)}) circle [radius={sqrt(2)}];
 
 // TRIANGLE
 \draw[fill = lightRed] (0,4.5) -- (2,{4-sqrt(2)}) -- (4,4.5) -- cycle;
 
 // BOUNDARY EDGE OF TRIANGLE
 \draw[dotted] (3,4.5)  -- (3.5, 5.2);
 \node at (3.6, 5.3) {$\bar{e}$};

 // HEIGHT OF TRIANGLE
 \draw (2,{4-sqrt(2)}) -- (2,4.5); 
 \draw[dotted] (2,{4-sqrt(2)/2})  -- (2.8, 5.2);
 \node at (2.9, 5.3) {$h_E$};

 // RADIUS OF BALL
 \draw (2,{4-sqrt(2)}) -- ({2-sqrt(2)},{4-sqrt(2)});
 \draw[dotted] ({2-sqrt(2)/2}, {4-sqrt(2)}) -- (2, 5.2);
 \node at (2.1, 5.3) {$R_E$};
 
 // BOUNDARY OF BALL
\draw [black] (2,{4-sqrt(2)}) circle [radius={sqrt(2)}];

\end{tikzpicture}

%% file: lemma_proof_4.tex
\begin{tikzpicture}[scale=0.65]

// POLIGONO
\draw[fill = lightGreen] (0,4.5)
  -- (0,0) 
  -- (2,1)
  -- (4,0) 
  -- (4,4.5);
 \node at (3.7,{4-sqrt(2)}) {\green{$E$}};

 // TRIANGLE
 \draw[fill = lightRed] (0,4.5) -- (2,{4-sqrt(2)}) -- (4,4.5);
 \node at (2, 3.6) {\red{$\breve{T}_{\bar{e}}$}};

// SURFACE
\draw [white,thick,domain=-1:5] plot (\x, {4.5+0.3*sin(90*\x)}); // serve solo a dare la giusta larghezza alla figure
\fill [color=lightRed, draw=lightRed, domain=0:2] plot (\x, {4.5+0.3*sin(90*\x)}) -- cycle;
\fill [color=white, draw=white, domain=2:4] plot (\x, {4.5+0.3*sin(90*\x)}) -- cycle;
\draw [red,thick,domain=0:4] plot (\x, {4.5+0.3*sin(90*\x)});
\node at (2,5.3) {$\red{\bolda(\bar{e})}$};

\end{tikzpicture}

%% file: post_processing_1_a.tex
\begin{tikzpicture}[scale=0.7]

// POLIGONO
\draw[fill=purple] (0,0) 
  -- (2,-2)
  -- (4,-0.3)
  -- (4,0)
  -- cycle;

// SURFACE
\draw [red, very thick,domain=-0.3:4.3] plot (\x, {-0.5*sin(45*\x)});

\fill (4,0) circle [radius = 3pt];
\fill (4,-0.3) circle [radius = 3pt];

\end{tikzpicture}

%% file: post_processing_1_b.tex
\begin{tikzpicture}[scale=0.7]

// POLIGONO
\draw[fill = purple] (0,0) 
  -- (2,-2)
  -- (4,0)
  -- cycle;

// SURFACE
\draw [red, very thick,domain=-0.3:4.3] plot (\x, {-0.5*sin(45*\x)});

\fill (4,0) circle [radius = 3pt];

\end{tikzpicture}

%% file: post_processing_2_a.tex
\begin{tikzpicture}[scale=0.7]

// POLIGONO
\draw[fill = purple] (0,0) 
  -- (1.3,-0.3)
  -- (4,-0.3)
  -- (2.6,0)
  -- cycle;
  
  // SURFACE
\draw [red, very thick, domain=-0.3:4.3] plot (\x, {(-2/10)*\x*\x*\x + (-3/56+33/25)*\x*\x + (39/280-52/25)*\x});

\fill (0,0) circle [radius = 3pt];
\fill (1.3,-0.3) circle [radius = 3pt];
\fill (4,-0.3) circle [radius = 3pt];
\fill (2.6,0) circle [radius = 3pt];

\end{tikzpicture}

%% file: post_processing_2_b.tex
\begin{tikzpicture}[scale=0.7]

// POLIGONO
\draw (0,0) 
  -- (1.3,0)
  -- (2.6,0)
  -- (4,-0.3);
  
  // SURFACE
\draw [red, very thick, domain=-0.3:4.3] plot (\x, {(-2/10)*\x*\x*\x + (-3/56+33/25)*\x*\x + (39/280-52/25)*\x});

\fill (0,0) circle [radius = 3pt];
\fill (1.3,0) circle [radius = 3pt];
\fill (4,-0.3) circle [radius = 3pt];
\fill (2.6,0) circle [radius = 3pt];

\end{tikzpicture}

%% file: elliptic_table.tex
\begin{tabular}{llllllll}
\hline\noalign{\smallskip}
$h$ & $L^2$ error & EOC & $L^\infty$ error & EOC & $N_\Omega$ & $N_\Gamma$ & $\text{cond}_{ell}$\\
\noalign{\smallskip}\hline\noalign{\smallskip}
$7.0711$e-1 & 5.1214e-02 & -      & 5.6622e-02 & -      & 16 & 12 & $8.3781$e+01\\
$3.5355$e-1 & 1.3589e-02 & 1.9141 & 1.4793e-02 & 1.9364 & 60 & 28 & $4.9466$e+02\\
$1.7678$e-1 & 3.6711e-03 & 1.8881 & 4.4904e-03 & 1.7200 & 224 & 60 & $3.0805$e+03\\
$8.8388$e-2 & 9.5200e-04 & 1.9472 & 1.1073e-04 & 2.0198 & 856 & 124 & $2.4991$e+04\\
$4.7611$e-2 & 2.4681e-04 & 2.1820 & 3.1169e-04 & 2.0491 & 3308 & 252 & $2.4124$e+05\\
\noalign{\smallskip}\hline
\end{tabular}   

%% file: elliptic_table_bsfem.tex
\begin{tabular}{llllllll}
\hline\noalign{\smallskip}
$h$ & $L^2$ error & EOC & $L^\infty$ error & EOC & $\text{cond}_{ell}$\\
\noalign{\smallskip}\hline\noalign{\smallskip}
$6.9133$e-1 & 1.1599e-01 & -      & 1.4298e-01 & -      &   $3.8803$e+01\\
$3.6409$e-1 & 2.6858e-02 & 2.2816 & 3.2995e-02 & 2.2869 & $2.3482$e+02 \\
$1.7787$e-1 & 7.4133e-03 & 1.7970 & 9.5269e-03 & 1.7340 & $1.7330$e+03\\
$8.7143$e-2 & 2.4630e-03 & 1.5444 & 3.1075e-03 & 1.5702 & $8.5300$e+03\\
$4.7612$e-2 & 6.8094e-04 & 2.1262 & 8.0268e-04 & 2.2386 & $5.6138$e+04\\
\noalign{\smallskip}\hline
\end{tabular}   

%% file: parabolic_table.tex
\begin{tabular}{llllll}
\hline\noalign{\smallskip}
$h$ & $\tau$ & $L^\infty(L^2)$ error & EOC & $L^\infty(L^\infty)$ error & EOC\\
\noalign{\smallskip}\hline\noalign{\smallskip}
$7.0711$e-1 & $1.000e$-2 & $3.2969$e-02 & -      & $3.4312$e-02 & -      \\
$3.5355$e-1 & $2.500e$-3 & $4.9408$e-03 & 2.7383 & $6.3263$e-03 & 2.4393 \\
$1.7678$e-1 & $6.250e$-4 & $5.3036$e-04 & 3.2197 & $1.3819$e-03 & 2.1947 \\
$8.8388$e-2 & $1.563e$-4 & $6.0452$e-05 & 3.1331 & $2.1639$e-04 & 2.6750 \\
$4.7611$e-2 & $3.906e$-5 & $9.8712$e-06 & 2.9292 & $6.7848$e-05 & 1.8277 \\
\noalign{\smallskip}\hline
\end{tabular}   

%% file: parabolic_table_bsfem.tex
\begin{tabular}{llllll}
\hline\noalign{\smallskip}
$h$ & $\tau$ & $L^\infty(L^2)$ error & EOC & $L^\infty(L^\infty)$ error & EOC\\
\noalign{\smallskip}\hline\noalign{\smallskip}
$6.9133$e-1 & $1.000e$-2 & $9.8785$e-02 & -      & $1.1235$e-01 & -      \\
$3.6409$e-1 & $2.500e$-3 & $1.4697$e-02 & 2.9715 & $2.0306$e-02 & 2.6679 \\
$1.7787$e-1 & $6.250e$-4 & $1.3004$e-03 & 3.3851 & $2.5841$e-03 & 2.8777 \\
$8.7143$e-2 & $1.563e$-4 & $1.5086$e-04 & 3.0190 & $4.6415$e-04 & 2.4064 \\
$4.7602$e-2 & $3.906e$-5 & $1.6979$e-05 & 3.6126 & $1.0210$e-04 & 2.5042 \\
\noalign{\smallskip}\hline
\end{tabular}   